\documentclass[mnsc,nonblindrev]{informs3_WP}

\OneAndAHalfSpacedXI

\usepackage{natbib}
\usepackage[]{bibunits}
 \bibpunct[, ]{(}{)}{,}{a}{}{,}%
 \def\bibfont{\small}%
 \def\bibsep{\smallskipamount}%

\TheoremsNumberedThrough
\ECRepeatTheorems
\EquationsNumberedThrough
\MANUSCRIPTNO{}

\usepackage{float}
\usepackage{nicefrac}
\usepackage{amsmath}
\usepackage{algorithm}
\usepackage[linesnumbered,vlined,ruled,commentsnumbered,algo2e]{algorithm2e}
\SetKwRepeat{Do}{do}{while}
\usepackage{dsfont}
\usepackage{hyperref}
\hypersetup{
   colorlinks   = true,
   urlcolor     = blue,
   linkcolor    = black,
   citecolor    = black, 
   naturalnames = true,
}
\usepackage[show]{ed}
\usepackage{fdsymbol}

\usepackage[position=bottom,caption=false]{subfig}
\usepackage[export]{adjustbox}
\usepackage{makecell}
\usepackage{csquotes}
\usepackage{multirow}
\makeatletter
\g@addto@macro{\UrlBreaks}{\UrlOrds}
\makeatother
\expandafter\def\expandafter\UrlBreaks\expandafter{\UrlBreaks
\do\a\do\b\do\c\do\d\do\e\do\f\do\g\do\h\do\i\do\j%
\do\k\do\l\do\m\do\n\do\o\do\p\do\q\do\r\do\s\do\t%
\do\u\do\v\do\w\do\x\do\y\do\z\do\A\do\B\do\C\do\D%
\do\E\do\F\do\G\do\H\do\I\do\J\do\K\do\L\do\M\do\N%
\do\O\do\P\do\Q\do\R\do\S\do\T\do\U\do\V\do\W\do\X%
\do\Y\do\Z}

\newcommand{\mdash}
           {\discretionary{}{}{\kern 0.1em}---\discretionary{}{}{\kern 0.1em}}

\begin{document}

\RUNAUTHOR{Blaettchen, Calmon, and Hall}
\RUNTITLE{Traceability Technology Adoption}
\TITLE{Traceability Technology Adoption \\ in Supply Chain Networks}

\ARTICLEAUTHORS{%
\AUTHOR{Philippe Blaettchen}
\AFF{Bayes Business School (formerly Cass), City, University of London, London EC1Y 8TZ, United Kingdom,  \EMAIL{philippe.blaettchen@city.ac.uk}} 
\AUTHOR{Andre P. Calmon}
\AFF{Scheller College of Business, Georgia Institute of Technology, Atlanta, GA 30308, \EMAIL{andre.calmon@gatech.edu}}
\AUTHOR{Georgina Hall}
\AFF{Decision Sciences, INSEAD, 77305 Fontainebleau, France, \EMAIL{georgina.hall@insead.edu}}
} 

\ABSTRACT{
Modern traceability technologies promise to improve supply chain management by simplifying recalls, increasing visibility, or verifying sustainable supplier practices. Initiatives leading the implementation of traceability technologies must choose the least-costly set of firms\mdash or \emph{seed set}\mdash  to target for early adoption. Choosing this seed set is challenging because firms are part of supply chains interlinked in complex networks, yielding an inherent \emph{supply chain effect}: benefits obtained from traceability are conditional on technology adoption by a subset of firms in a product's supply chain.
We prove that the problem of selecting the least-costly seed set in a supply chain network is hard to solve and even approximate within a polylogarithmic factor. Nevertheless, we provide a novel linear programming-based algorithm to identify the least-costly seed set. The algorithm is fixed-parameter tractable in the supply chain network's treewidth, which we show to be low in real-world supply chain networks. The algorithm also enables us to derive easily-computable bounds on the cost of selecting an optimal seed set.
Finally, we leverage our algorithms to conduct large-scale numerical experiments that provide insights into how the supply chain network structure influences diffusion. These insights can help managers optimize their technology diffusion strategy.
}

\HISTORY{This paper has been accepted for publication by \emph{Management Science} on August 27, 2023.}

\begin{bibunit}[informs2014]

\maketitle

\vspace{-0.8cm}
\section{Introduction} \label{sec:introduction}

Modern consumer goods supply chains form complex networks spanning dozens of countries and actors. As a result, most firms cannot reliably trace the products they produce and source beyond a few upstream and downstream supply chain tiers. This limited traceability\mdash or ability to trace the processing history, origin of materials, and final destination of products \citep{ISO9000}\mdash has several negative consequences.
First, a lack of traceability is a major barrier to building sustainable, disruption-resilient supply chains \citep{WH2022}, limits supply chain coordination, and increases transaction costs \citep{Wilson2014}. Second, firms lacking traceability are prone to extensive recalls \citep{Wowak2016}, with adverse consequences \citep{Lee2022}. Finally, customers increasingly value traceability, so a lack of it can negatively affect demand \citep{RL2016}.

To counter this traceability deficiency, many firms are leading the development and deployment of new traceability technologies, protocols, standards, and initiatives (henceforth technologies).\footnote{In a recent survey of 150 senior supply chain leaders, 68\% identify traceability as a ``very or extremely" important issue \citep{WEF2021}.
The value of the traceability technology industry is estimated to grow to US\$23 bn.\ by 2025 \citep{Bhandalkar2019}.}
In the food industry, for instance, these traceability initiative leaders are typically large retailers or processors (e.g., Walmart or Tyson Foods), often in collaboration with IT companies (e.g., IBM and its ``Food Trust" traceability solution) or  industry consortia \citep[see][for more examples]{Naidu2017, Youngdahl2018,Haig2020,Hiba2023}. Similarly, in the fashion industry, fast-growing IT companies (such as Textile Genesis) engage with large retail companies (such as H\&M, Lenzing, or Bestseller) to lead traceability initiatives \citep{ahmed2021blockchain}.
Traceability initiative leaders aim to have all players in target supply chains adopt their traceability technology or, more ambitiously, have their technology become the industry standard for specific product categories. However, they often find disseminating their technology across supply chains a daunting and ``painfully complex" task \citep{saenz2022traceability}.

This complexity is due to two key reasons. The first is a network effect unique to supply chain technologies called the \emph{supply chain effect}. To illustrate this effect, consider a producer of chocolate-based products who wishes to use a traceability technology to trace the origins of its cocoa and the final destination of its products.
While traceability may benefit all players in the chocolate supply chains (e.g., by improving the detection of supplier malpractice, increasing demand visibility, or enabling sustainability certifications), benefits are only obtained if products are traceable throughout their entire supply chains.
If a subset of firms in the supply chain does not adopt the technology, the product produced by that supply chain is no longer fully traceable, and technology adoption benefits are drastically diminished for all firms. Thus, the supply chain effect requires most or even all firms involved in the product's supply chain to adopt the traceability technology for firms to benefit \citep[see, e.g.,][]{Behnke2020,Sternberg2020}. The second reason is the complex structure of modern supply chains. Traceability initiative leaders interact with thousands of firms in hundreds of partially overlapping supply chains, thus forming a \emph{supply chain network}. While these overlapping supply chains help alleviate the supply chain effect, they complicate the design of traceability technology dissemination strategies.

Traceability initiative leaders often proactively engage with a set of early adopter firms in the supply chain network to jumpstart technology diffusion. We refer to this set as the network's \emph{seed set}. Engagement with the seed set is costly and usually includes pilot programs, subsidies, or cost-sharing incentives. The leader's goal is to have seed set firms adopt the technology and then influence other firms into adopting it, triggering broad technology diffusion \citep{WEF2021, saenz2022traceability}. Engaging with the ``best" seed set is a critical decision for traceability leaders, and, to build an effective technology dissemination strategy, they must answer a few vital managerial questions: \emph{(i) What is the lowest-cost seed set that ensures the whole network eventually adopts the technology?} \emph{(ii) How does the size of the seed set depend on the network structure?}  and \emph{(iii) What are the different roles that seed set firms play in the diffusion process?}

The existing network diffusion literature addresses these questions through various mathematical models and frameworks  \citep[e.g.][]{rogers2010diffusion}. 
However, none of these models are tailored to supply chain networks, nor the specificities of traceability technology. In particular, there is no research on how the supply chain effect influences technology diffusion. Our paper aims to fill this gap by introducing a new model which incorporates the supply chain effect and can be used to optimize the dissemination of traceability technology in supply chain networks. This model enables us to answer the strategic questions above and can guide the design of traceability technology diffusion strategies. 
From a theoretical perspective, our framework extends and applies recent results from Integer Programming to technology diffusion, building a new bridge between these two fields.

More specifically, our theoretical contributions are as follows. In Section~\ref{sec:SCTM}, we introduce our new technology diffusion model, the Supply Chain Traceability Model (SCTM), and formalize the seed set selection problem ($MIN$-$SCTM$). In Section~\ref{sec:Comp.Comp},  we prove that $MIN$-$SCTM$ is hard to solve and approximate and that the supply chain effect drives this complexity. In light of this result, any exact solution algorithm for $MIN$-$SCTM$ must be parametrized by a structural parameter of the network. We propose such an algorithm in Section \ref{sec:exact.sols}, more specifically, a fixed-parameter tractable (FPT) linear programming-based algorithm with parameter \emph{treewidth} of the supply chain network. This parameter measures how ``tree-like'' the network is and is low for real-world supply chains. We further provide two approximation schemes for $MIN$-$SCTM$. One is a principled heuristic that returns upper and lower bounds on the optimal cost, explicitly trading off accuracy with computational time (Section \ref{sec:exact.sols}). The other is a simple heuristic based on our managerial insights (Section~\ref{sec:heuristics}).
These algorithms and heuristics collectively answer question \emph{(i)} above.

We then conduct a series of large-scale numerical experiments using our optimization framework to answer questions \emph{(ii)} and \emph{(iii)}. In Section \ref{sec:insight1}, we address \emph{(ii)} and find that the \emph{Jaccard clustering} of a supply chain network is a crucial predictor of the seed set size. This measure can thus estimate the effort required to disseminate a traceability technology.
 Section \ref{sec:insight2} examines \emph{(iii)}. We observe two types of early adopter firms in the seed set: \emph{starter} and \emph{helper} firms. Starter firms are positioned within supply chains that are made traceable early in the diffusion process and help ``jumpstart" diffusion. Conversely, helper firms are part of supply chains that become traceable at later stages of the diffusion process. These firms help circumvent the supply chain effect and ``transfer" diffusion across different network parts. We show that the ratio of starter-to-helper nodes in the seed set has a non-linear relationship with a network's \emph{modularity}. Our insights can help managers tailor their diffusion strategy to a supply chain network's structure.

\section{Literature Review} \label{sec:literature.review}

Most operations management papers on traceability technologies focus on the tools and IT infrastructure supporting traceability, ranging from RFID \citep{dutta2007rfid,heese2007inventory,whang2010timing} to data management systems such as blockchain technology \citep{babich2020om,chod2020financing,cui2020values}. In contrast, our paper abstracts away specific technological details and focuses on diffusion across supply chains, an issue that affects all traceability technologies. 

Our approach contributes to the network diffusion literature by building on the Linear Threshold Model (LTM) \citep{granovetter1978threshold}. In the LTM, a node in a network adopts an innovation (such as a new technology) after a certain fraction of its neighbors have adopted the same innovation. The model we develop, the SCTM, relates to the LTM in two ways. First, the SCTM is a generalization of the LTM, allowing for interactions through hyperedges, not only direct neighbors. Thus, the SCTM encodes the supply chain effect\mdash where \emph{sets}  of firms in a supply chain must adopt the technology for adoption benefits to become available\mdash which the LTM cannot do directly. Second, as described in Section \ref{sec:aux_graph}, by introducing an auxiliary graph to our supply chain network hypergraph, one can view the SCTM as a weighted generalization of the LTM on a highly structured graph. The auxiliary graph enables us to relate our results to existing results for the target set selection (TSS) problem in the LTM. For example, our hardness and inapproximability results for $MIN$-$SCTM$ add to the results on the hardness of TSS under additional structural assumptions \citep{Kempe2003,chen2009approximability}. We also provide an extension of the dynamic programming-based fixed-parameter tractable algorithm for TSS by \citet{ben2011treewidth} to our setting.

This paper also builds on recent results from the integer programming literature.
Specifically, \citet{laurent2009sums} and \citet{bienstock2018lp} show that binary linear programs are amenable to linear programming reformulations that are FPT in the treewidth of a graph related to the original formulation. We employ these results to derive an LP-based FPT algorithm for $MIN$-$SCTM$. By doing so, we open up a new application area for these integer programming techniques in network diffusion while providing an innovative approach to TSS in the LTM.\footnote{While we do not apply our algorithm to this specific problem, we could easily extend our approach.} Namely, our LP formulation overcomes the implementation difficulties of the existing dynamic programming approach and provides a principled way for designing heuristics and obtaining managerial insights.

\section{The Supply Chain Traceability Model (SCTM)} \label{sec:SCTM}

We introduce our model in Section \ref{subsec:model.description} and show how it describes different supply chain relationships in Section \ref{subsec:case_studies}. In Section  \ref{sec:aux_graph},  we present an auxiliary graph that is key to solving $MIN$-$SCTM$.

\subsection{Model Description}\label{subsec:model.description}

\begin{figure}[t]
    \centering
    \includegraphics{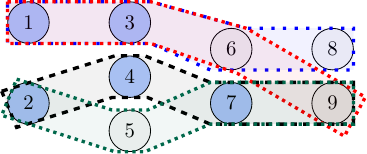}
    \caption{Example of a hypergraph $G$ with $n=9$ and $m=4$. Here, $N_F=\{1,\ldots,9\}$ and $E=\{e_{blue}, e_{red}, e_{green}, e_{black}\}$ with $e_{blue}=\{1,3,6,8\},~ e_{red}=\{1,3,6,9\},~ e_{green}=\{2,5,7,9\},~ e_{black}=\{2,4,7,9\}$.}
    \label{fig:adoption_example}
\end{figure}

Consider a hypergraph $G$ with nodes $N_F=\{1,\ldots,n\}$ and hyperedges\footnote{A hyperedge $e_j, j=1,\ldots,m$ in hypergraph $G$ is a subset of nodes in $N_F$.} $E=\{e_1,\ldots,e_m\}$. Figure \ref{fig:adoption_example} provides an example. Nodes could represent firms and hyperedges subsets of firms collaborating to produce a product (see Section \ref{subsec:case_studies}).
The size of $e_j$ is $k_j$, and $k_j$ is upper-bounded by a constant $k$. Without loss of generality, we assume $G$ is connected\mdash otherwise, we repeat our analysis on each connected component independently.

\subsubsection*{State of the Network.} Each node $i \in N_F$ has a state $x_{it} \in \{0,1\}$ at the end of period $t \in \{0,1,2,\ldots\}$. A node in state $1$ is said to be \emph{active} or have \emph{adopted the technology}. Once active, a node remains in state 1 for all future periods. The \emph{state of the network} at time $t$, given by $S_t$, is the set of active nodes at the end of period $t$, that is $S_t = \{i \in N_F\ :\ x_{it} = 1\}$.

\subsubsection*{Parameters of the Network.} For each node $i$ and hyperedge $e_j$, there is a \emph{benefit} $r_{ji}$  provided by $e_j$ to $i$. For each node $i$, there is also an \emph{adoption} cost, $c_i$,  and a \emph{seeding} cost, $w_i$. Each hyperedge $e_j$ has an \emph{adoption threshold}  $\theta_j$, which is the minimum number of nodes in $e_j$ that need to be active for $e_j$ to become \emph{active} or \emph{traceable}.

\subsubsection*{Activation Process and State Equation.}
In periods $t \in \{1,2,\ldots\}$, a node $i \in N_F$ in state 0 decides whether to become active and switch to state 1. Node $i$ switches states if the adoption cost $c_i$ is outweighed by the adoption benefit $b_i(S_t)$, computed in the following way. Consider the set of all hyperedges $e_j$ that $i$ belongs to and which are traceable \emph{after node $i$ decides to become active}. We let the set of indexes of these hyperedges be $\mathcal{B}_i(S_{t})$. Then,  $\mathcal{B}_i(S_{t}) = \{ j\ :\ i \in e_j,\  |S_{t} \cap e_j | \geq  \theta_j - 1\}$.\footnote{ This definition comes from the following argument: Let $j \in \mathcal{B}_i(S_{t}).$ Then, $i \in e_j$ and $\sum_{i' \in e_j} x_{i't} \geq \theta_j$. If $i$ decides to activate, then $x_{it}=1$. Therefore, $\sum_{i' \in e_j} x_{i't}=\sum_{i'\neq i, i'\in e_j} x_{i't}+1 \geq \theta_j$, which is equivalent to $|S_{t} \cap e_j | \geq  \theta_j - 1.$ }

Each hyperedge $e_j$ with $ j \in \mathcal{B}_i(S_t)$  generates traceability benefit $r_{ji} \geq 0$ for node $i$ once the node is active. Thus, the benefit $i$ obtains from becoming active is $b_i(S_{t}) = \sum_{j \in \mathcal{B}_i(S_{t})} r_{ji}$.
We define the activation process's state equations, which we call the \emph{SCTM activation process}, as
$$S_{t+1} = S_{t} \cup \{ i \in N_F \setminus S_{t}\ :\ b_i(S_{t}) \geq c_i\},\ \forall t = 0, 1, \ldots.$$
Given a hypergraph $G$ and some initial set $S_0 \subseteq N_F$, this process is well-defined: A unique final set of adopters $S_\infty \subseteq N_F$ exists and can be attained in a finite number of steps.

\subsubsection*{Activation Process Example.}
Figure \ref{fig:adoption_example} illustrates the activation process. We focus on the network structure's impact and set $r_{ji} \geq c_i\ \forall i \in N_F, e_j \in E$. Let the initial set of active nodes that have adopted the technology be $S_0 = \{1,2,3,4,7\}$ and assume we are in period $t = 1$. We further assume that traceability benefits only kick in when all nodes in a hyperedge have adopted the technology, that is, $\theta_j~=~4\ \forall e_j \in E$. In the case of 
$e_{blue}$, $e_{red}$, and $e_{green}$, there are two active nodes, so no node can obtain a benefit from these hyperedges by becoming active. However, hyperedge $e_{black}$ has three active nodes, with only Node 9 inactive. By becoming active, Node 9 ensures that $e_{black}$ is active or traceable and obtains benefits. In particular, $\mathcal{B}_9(S_0) = \{ j\ :\ 9 \in e_j,\  |S_0 \cap e_j | \geq  3\} = \{ black\}$ and $b_9(S_0) = r_{black,9}$. Because $r_{black,9} \geq c_9$, Node 9 becomes active, and $S_1 = \{1,2,3,4,7,9\}$. 
In period $t = 2$, hyperedges $e_{red}$ and $e_{green}$ are one node away from becoming active via Nodes 6, respectively 5. Consider Node 6: $\mathcal{B}_6(S_1) = \{ j\ :\ 6 \in e_j,\  |S_1 \cap e_j | \geq  3\} = \{ red\}$, so its benefit is $b_6(S_1) = r_{red,6} \geq c_6$. Similarly for Node 5: $\mathcal{B}_5(S_1) = \{ j\ :\ 5 \in e_j,\  |S_1 \cap e_j | \geq  3\} = \{ green\}$, so the benefit is $b_5(S_1) = r_{green,5} \geq c_5$. Both nodes become active, and $S_2 = \{1,2,3,4,5,6,7,9\}$. This leaves Node 8 as the only inactive node at the end of period $t = 2$. Since $\mathcal{B}_8(S_2) = \{ j\ :\ 8 \in e_j,\  |S_2 \cap e_j | \geq  3\} = \{ blue\}$, we have $b_8(S_2) = r_{blue,8}\geq c_8$ and Node 8 becomes active in period 3. Thus, $S_3 = S_\infty = N_F$.

\subsubsection*{Decision.}
At time 0, the decision-maker chooses the initial set of adopters, or \emph{seed set} $S_0 \subseteq N_F$. Other nodes are in state 0. The decision-maker could be, for example, a retail chain interested in tracing the origins of a set of products it sells. We provide a detailed discussion in Section \ref{subsec:case_studies}.

\subsubsection*{Problem Formulation.}
The decision-maker chooses the lowest-cost seed set $S_0$ so all nodes $i \in N_F$ eventually become active. We call this problem $MIN$-$SCTM$:
\begin{equation}\label{eq:tss.sctm}
\begin{aligned}
    OPT = \min_{S_0 \subseteq N_F}\ \ \ & \sum_{i \in S_0} w_i  \\
    \text{s.t.} \ \ \ & S_{\infty} = N_F  \\
    & S_{t+1} = S_{t} \cup \{ i \in N_F \setminus S_{t}\ :\ b_i(S_{t}) \geq c_i\},\ \forall t=0,1,\ldots  
\end{aligned}
\end{equation}
Here, given a budget constraint, we minimize the seed set cost to achieve a final set of adopters rather than maximize the size of the final set of adopters. This is because the former is more appealing in practice: network effects and economies of scale inherent to traceability technology will eventually lead to a single industry-wide standard. For a decision-maker to reap benefits in the long term, broad adoption is needed, even at a high initial cost. In addition, we constrain this final set of adopters to be $N_F$, as it leads to an upper bound on the seeding cost of any subset of $N_F$. Our results can be directly extended to any subset of $N_F$.

\subsubsection*{Assumptions on the Parameters.} We assume that all nodes can become active, that is, $\sum_{j : i \in e_j} r_{ji} \geq c_i\ \forall i \in N_F$. We also assume that all values $r_{ji}$ and $c_i$ are integers\mdash values can be scaled if they are initially rational. For a given $i \in N_F$, if the greatest common divisor $g_i$ of $\{r_{ji}\}_j$ and $c_i$ is larger than one, we consider benefits $\{\nicefrac{r_{ji}}{g_i}\}_j$ and cost $\nicefrac{c_i}{g_i}$. This is a useful preprocessing step, as the runtime of our subsequent algorithm will scale with the maximum of costs and benefits. Moreover, we assume  $\theta_j \in \{2,3,\ldots,k_j\}$. If $\theta_j > k_j$, hyperedge $e_j$ never provides traceability benefits, and we remove it from $G$. If $\theta_j =1$, traceability benefits prevail when just one node is active, so we can remove the hyperedge and reduce the adoption cost of each node $i \in e_j$ by the node's benefit, that is, $c_i := c_i - r_{ji}$.
Finally, we assume that $c_i > 0$: If $c_i \leq 0$ for a node $i$, we can remove the node and let $\theta_j := \theta_j - 1$ for all $j$ with $i \in e_j$. If this leads to $\theta_j < 2$, we repeat the previous simplification.

\subsection{Model Discussion and Examples}\label{subsec:case_studies} 

We assume that $G$ is a \emph{supply chain network}. 
A node $i \in N_F$ represents a \emph{firm}, and a hyperedge $e_j \in E$ a \emph{supply chain}, that is, a subset of firms collaborating to produce product $j$. Suppose $G$ is the supply chain network of a given product category (such as fruits, vegetables, or dairy). Then, nodes may not be entire firms but rather divisions handling this product category. Thus, our model allows for \emph{parts of firms} to adopt traceability technology rather than firms in their entirety.

The benefit $r_{ji}$ obtained by firm $i$ if the supply chain of product $j$ is traceable can correspond to a guarantee of continuing or expanding demand for the product or other more intangible benefits such as improved supply chain resiliency, better coordination, or lower quality costs. As traceability benefits are manifold and take different forms at different supply chain stages, they can be challenging to estimate, and traceability initiatives spend significant efforts on this task. Absent detailed information on the value of $r_{ji}$ for each firm, an initial estimate for $r_{ji}$ can be obtained by considering the total traceability benefits to a supply chain $j$ (which are simpler to estimate than the individual benefits) and prorating them to $r_{ji}$ according to the value added by each firm.

The adoption threshold $\theta_j$ is the number of firms in hyperedge $e_j$ that must adopt the traceability technology for product $j$ to become traceable. We assume that all nodes in $e_j$ contribute equally to this threshold, though, in practice, contributions may be unequal. We can extend our model to this setting, but the notation and analysis become more complex, so we leave this to future research. If $\theta_j$ is unknown, one can set $\theta_j=k_j$, which implies that \emph{all} firms in the supply chain of product $j$ must be active to obtain traceability benefits. This worst-case scenario is a natural assumption for several traceability applications and provides an upper bound on the optimal seeding cost.

The adoption cost $c_i$ of the traceability technology by firm $i$ can represent IT, auditing, and training costs and all discounted future costs of operating the technology. If firms receive technology adoption benefits independently of other firms' adoption, we can consider $c_i$ as the adoption cost net of such benefits. Our adoption process assumes that this cost is a one-time cost for each firm and that, once the technology has been adopted by firm $i$, the technology is available for all other products (if any) processed by firm $i$. Without detailed firm-level information on $c_i$, one can set it to reflect the size of $i$ plus any fixed technology adoption costs.

\subsubsection*{Decision-Makers.} The most common decision-maker for \eqref{eq:tss.sctm} is a large firm wishing to make products it processes traceable, i.e., it is a node of $G$. We refer to such a node as the network's \emph{lead firm}.\footnote{We assume we only have one lead firm, but our discussion can easily be adapted for multiple lead firms.} Alternatively, the decision-maker can be a traceability initiative external to the supply chain network, such as a large IT company or certification organization, or even a mixture of internal and external players \citep{Naidu2017,Youngdahl2018}.

Our model can deal with both of these scenarios. If the decision-maker is an external initiative, it can solve \eqref{eq:tss.sctm} directly to obtain the optimal seed set. If the decision-maker is internal to the network, i.e., it is a lead firm, it can also solve \eqref{eq:tss.sctm} to obtain the optimal seed set. However, the hypergraph $G$ over which \eqref{eq:tss.sctm} is solved will be slightly different. The set $N_F$ of nodes will include an additional node, corresponding to the lead firm, and this node will have both its seeding and adoption costs set to zero.\footnote{Through its role, the lead firm would have adopted the technology at the start of the diffusion process.} To avoid having to distinguish between both cases for the remainder of the paper, we note that solving \eqref{eq:tss.sctm} over $G$ as defined this way is equivalent to solving \eqref{eq:tss.sctm} over a slightly modified graph obtained by deleting the lead firm node from $G$ and setting $\theta_j\mathrel{\mathop{:}}=\theta_j-1$ for all hyperedges $e_j$ to which this node belongs (see Section \ref{subsec:model.description}). In other words, in settings where the decision-maker is internal, one can modify the underlying hypergraph $G$ to revert to the external decision-maker case, so we do not need to distinguish between the two cases.

\subsubsection*{Assembly Networks and Aggregators.} Our model can represent an assembly network, where multiple inputs are needed to produce an output. An example is given in Figure \ref{fig:power_a}, with six nodes and two hyperedges, $e_{blue}=\{1,2,4,6\}$ and $e_{red}=\{3,5,6\}.$ Here, Node 4 converts the inputs from Nodes 1 and 2 into a single output. Thus, traceability benefits from $e_{blue}$ only kick in if \emph{both} Nodes 1 and 2 (as well as 4) adopt, accounted for by $\theta_{blue}=4$: if only Nodes 1 and 4 have adopted, then Node 6 will not get any benefits. This example can be generalized to more complex settings.

\subsubsection*{Power Dynamics in Supply Chain Networks.} Hyperedges can also encode more intangible relationships, such as power dynamics, where different firms within a supply chain may have asymmetric abilities to influence technology adoption. For example, a large retailer such as Walmart can have considerable influence over technology adopted by its suppliers \citep{Nash2018}.
Our model can encode such dynamics through additional hyperedges:

\begin{figure}[t]
    \centering
    \subfloat[Assembly network.]{
        \includegraphics{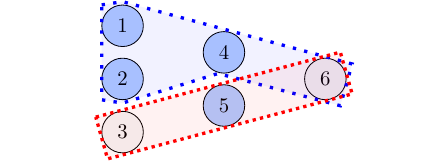}
        \label{fig:power_a}
    }
    \hspace{0.1cm}
    \subfloat[Additional hyperedge for power dynamics.]{
        \includegraphics{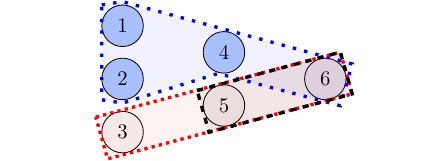}
        \label{fig:power_b}
    }
    \caption{Representing assembly networks and power dynamics using our model.}\label{fig:power}
\end{figure}

Consider first the hypergraph in Figure \ref{fig:power_a}. Assuming that $w_i = w_{i'}\ \forall i,i' \in N_F$, that $r_{ji} \geq c_i$, $\forall i \in N_F, e_j \in E$, and that $\theta_{blue}= 4$ and $\theta_{red}=3$, then $S_0=\{1,2,4,5\}$ is an optimal seed set.
Now, say that Node 6, upon adoption, forces Node 5 to adopt the same technology. We can model this effect by introducing an additional hyperedge $e_{black} = \{5,6\}$ in Figure \ref{fig:power_b}, with threshold $\theta_{black}= 2$, benefits to Node 6 of zero ($r_{black,6}=0$), and benefits to Node 5 corresponding to its adoption costs ($r_{black,5} = c_5$). If, for example, Node 5 adopts before Node 6, then the additional hyperedge has no effect. If, however, Node 6 adopts at time $t$, while Node 5 has not adopted, then $black \in \mathcal{B}_5(S_t)$. But $r_{black,5} \geq c_5$, so Node 5 will adopt. Thus, the seed set $S_0 = \{1,2,4\}$ with one less node is optimal.

\subsubsection*{Model Limitations.} While $SCTM$ can describe rich supply chain relationships, it has a few limitations. For example, the model does not describe potential changes in supply chain relationships and firms' sourcing strategies due to  technology adoption decisions.  Furthermore, it assumes firms act myopically in each period (a common assumption in network diffusion models) and does not describe more sophisticated strategic games between firms. Enriching $SCTM$ to model evolving supply chain dynamics and strategic behavior could be a fruitful source of research and insights. 

\subsection{The Auxiliary Graph} \label{sec:aux_graph}

A crucial construct for solving $MIN$-$SCTM$ is the \emph{auxiliary graph} $G'=(N',E')$ of hypergraph $G$, which will help connect the SCTM to the linear threshold model (LTM). As mentioned in Section~\ref{sec:literature.review}, the LTM is a popular model of diffusion in networks where a node becomes active if the number of active neighbors exceeds the node's threshold \citep{Kempe2003}.

\begin{figure}[t]
    \centering
    \subfloat[Hypergraph $G$.]{
        \includegraphics{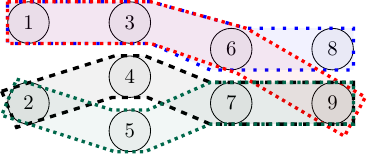}
        \label{fig:aux_graph_orig} 
    }   
    \hspace{0.5cm}
    \subfloat[Auxiliary graph $G'$.]{
        \includegraphics{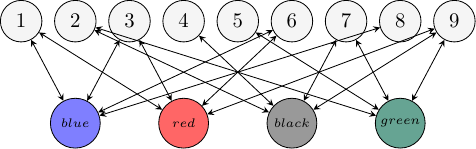}
        \label{fig:aux_graph_aux} 
    }
    \caption{The supply chain network from Figure \ref{fig:adoption_example}, as well as the associated auxiliary graph $G'$.}\label{fig:aux_graph}
\end{figure}

Define $G'$ to be a bipartite graph with a node for each firm $i \in N_F$ (``firm-node'') and a node for each supply chain $e_j \in E$ (``SC-node'', denoted with $j$). We let $N_{SC}$ be the set of SC-nodes, so $N'=N_F \cup N_{SC}.$ Each node $i \in N_{F}$ (resp.\ $j \in N_{SC}$) has a threshold $c'_i=c_i$ (resp.\ $c'_j = \theta_j - 1$). The graph is \emph{weighted} and \emph{directed}. Edges in $E'$ are added as follows: if $i \in e_j$ in $G$, then $(i,j) \in E'$ with weight $w'_{ij}=1$ and $(j,i) \in E'$ with weight $w'_{ji}=r_{ji}$. Figure \ref{fig:aux_graph} recalls graph $G$ from Figure \ref{fig:adoption_example} and displays the corresponding auxiliary graph $G'$ (without weights for legibility).

We next define the activation process on $G'$: A (firm or SC) node $i$ becomes active at time $t$ ($x'_{it}=1$) if the sum of the incoming edge weights from active nodes exceeds its threshold. In other words, if $i \in N_F$, $i$ becomes active if $\sum_{\{(j,i) \in E'~:~j \text{ is active} \} } r_{ji} \geq c_i$. If $j \in N_{SC}$, $j$ becomes active if $|\{(i,j) \in E'~:~i \text{ is active} \}|  \geq \theta_j - 1$. The activation process on the auxiliary graph can be viewed as a weighted-edge version of the LTM, so we refer to it as the \emph{LTM activation process}. 
In an analogous way to the SCTM activation process on $G$, we  define state equations, given a set $S_0' \subseteq N'$:
$$S'_{t+1}=S'_t \cup \left\lbrace i \in N'\backslash S'_t~:~ \sum_{j: (j,i) \in E'} w'_{ji}x_{jt} \geq c_i' \right\rbrace,\ \forall t=0,1,\ldots,$$
where $S'_t$ corresponds to the set of nodes active at time $t$ in $G'$. As before, a unique final set of adopters $S_{\infty}'$ is attained in a finite number of steps. Our first result shows that any SCTM activation process on $G$ can be replicated via an LTM activation process on $G'$:

\begin{proposition} \label{prop:sctm.eq.ltm}
Let $S_0 = S_0' \subseteq N_F$. Then, $S_t=S'_{2t} \cap N_F\ \forall t=0,1,2,\ldots$, and $S_{\infty}=S'_{\infty} \cap N_F$.
\end{proposition}

Proofs for this section are in Appendix \ref{appendix:model}. Figure \ref{fig:scn_ltm_equi} exemplifies the processes' equivalence. In both graphs, we start with the seed set $S_0 = S_0' = \{1,2,3,4,7\}$. Under the SCTM, Node 9 becomes active in $t=1$ because it can make $e_{black}$ active, Node 5 (resp.\ 6) in $t=2$, because it can make $e_{green}$ (resp.\ $e_{red}$) active, and Node 8 in $t=3$, because it can make $e_{blue}$ active. Consider now the LTM: At time $t=1$, Node ${black}$ corresponding to $e_{black}$ has three incoming active neighbors, while its threshold is three. Hence, it becomes active. With Node ${black}$ active, Node 9 has one incoming active neighbor. Because $w'_{black,9}=r_{black,9} \geq c_9 = c'_9$, it becomes active, and $S'_2 \cap N_F = \{1,2,3,4,7\} = S_1$. Now, both of the Nodes ${red}$ and ${green}$ have three incoming active neighbors, and, again, their thresholds are three, so they become active at time $t=3$, followed by Nodes 5 and 6 at time $t=4$. Again, $S'_4 \cap N_F = \{1,2,3,4,5,6,7\} = S_2$. Repeating this one last time, we see that $S'_6 \cap N_F = N_F = S_3$.

\begin{figure}[t]
    \centering
    \subfloat[The SCTM activation process.]{
        \includegraphics{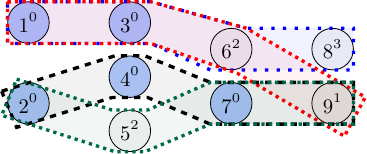}
        \label{fig:sctm_ltm_1} 
    }   
    \hspace{0.5cm}
    \subfloat[The LTM activation process.]{
        \includegraphics{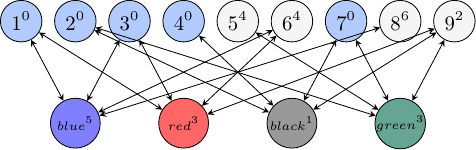}
        \label{fig:sctm_ltm_2} 
    }
    \caption{Equivalence between the SCTM and LTM activation processes. Superscripts represent the period in which a node becomes active.}\label{fig:scn_ltm_equi}
\end{figure}

Our next result formalizes the minimum cost seed set problem on the auxiliary graph.

\begin{corollary} \label{cor:sctm.ltm}
Optimization problem (\ref{eq:tss.sctm}) is equivalent to the optimization problem
\begin{equation}\label{eq:tss.ltm}
\begin{aligned}
    OPT = \min_{S'_0 \subseteq N_F}\ \ \ & \sum_{i \in S'_0} w_i  \\
    \text{s.t.} \ \ \ & S_{\infty}' = N'  \\
    & S'_{t+1}=S'_t \cup \left\{i \in N'\backslash S'_t\ :\ \sum_{j: (j,i) \in E'} w'_{ji}x_{jt} \geq c_i' \right\},\ \forall t=0,1,\ldots  
\end{aligned}
\end{equation}
\end{corollary}
Problem \eqref{eq:tss.ltm} has similarities with the \emph{target set selection (TSS)} problem under the LTM, but also several differences (see Section \ref{sec:literature.review}). In a nutshell, the LTM activation process we consider is a weighted version of the regular LTM activation process. Moreover, the underlying graph in the TSS problem is generic. Here, $G'$ has some additional structure:
It is bipartite and has constraints on the nodes' degrees and thresholds. Hence, one cannot leverage existing results from the literature to show the problem's hardness (see Section \ref{sec:Comp.Comp}). However, one can typically adapt algorithms for solving the TSS to \eqref{eq:tss.ltm} though they can be improved by leveraging the additional structure.

\section{Computational Complexity} \label{sec:Comp.Comp}

We turn to the computational complexity of $MIN$-$SCTM$. To isolate the effects of the graph's structure on the complexity of $MIN$-$SCTM$, we consider the simplest setting with trivial cost-benefit analysis. Namely, we let $r_{ji}=c_i=1\ \forall i=0,\ldots,n,\ j=1,\ldots,m$. Thus, any node whose activation would make a supply chain traceable will activate as the benefits will automatically outweigh the costs. This simplified setting can also be of interest in its own right, for example, in traceability systems that aim to restrict counterfeit drugs \citep[see, e.g.,][]{Lock2019}. We further assume that $k_j=k$ and $\theta_j=\theta\ \forall j=1,\ldots,m$. We define the decision version of Problem \eqref{eq:tss.sctm} under these assumptions and provide a \emph{full characterization} of the difficulty of answering it.

\begin{definition} \label{def:dec.sctm}
$DEC$-$SCTM$ is the decision version of $MIN$-$SCTM$, with simplified $G$: \\
\textsc{Input:} Integer $h$; hypergraph $G$ as defined in Section \ref{sec:SCTM} with benefits $r_{ji}=1$, adoption costs $c_i=1$, edges with $k_j=k$ and $\theta_j=\theta$, and rational seeding costs $w_i$, for all $i=0,\ldots,n,\ j=1,\ldots,m$.\\
\textsc{Question:} Is there a seed set $S_0$ of cost $\sum_{i \in S_0} w_i \leq h$ leading to full (SCTM) activation of $G$?
\end{definition}

\begin{theorem} \label{thm:comp}
The hardness of answering $DEC$-$SCTM$ depends on $k$ and $\theta$ as follows:
\vspace{0.1cm}

\begin{center}
    \begin{tabular}{ccccc}
    \hline
         \up\down $\theta$ / $k$ & $k=1$ & $k=2$ & $k=3$ & $k\geq 4$ \\
         \hline
         \up
         $\theta=1$ & in P & in P & in P & in P\\
         $\theta=2$ &  & in P & in P & in P\\
         $\theta=3$ &  &    & NP-hard & NP-hard\\
         \down $\theta \geq 4$ & & & & NP-hard\\
         \hline
    \end{tabular}
\end{center}
\end{theorem}
\vspace{0.3cm}

This section's proofs are in Appendix \ref{appendix:comp}. The theorem states that if supply chains in $G$ contain three or more firms, and at least three are needed for traceability benefits, $DEC$-$SCTM$ is hard to answer. Then, if $P\neq NP$, there is no hope of a polynomial-time algorithm for $MIN$-$SCTM$. 

From a managerial perspective, Theorem \ref{thm:comp} shows that optimizing technology dissemination when the benefits to a firm only depend on one supplier's or buyer's adoption decision ($\theta \leq 2$) is ``easy'' (in P). However, once the supply chain effect is in place and the benefits of a traceability technology depend on adoption by second-tier suppliers and buyers ($\theta \geq 3$), there is a phase shift, and optimization becomes complex (NP-hard).
 
\subsubsection*{Relationship to Hardness of Target Set Selection under the LTM.} 
Due to the specific weights and thresholds used in Theorem \ref{thm:comp} and Corollary \ref{cor:sctm.ltm}, Theorem \ref{thm:comp} is equivalent to the following result: The decision version of target set selection under the LTM is NP-hard if the underlying graph $G'$ is bipartite, with one set of nodes, $N_{F}$, having threshold 1, and the other set of nodes, $N_{SC}$, having threshold $\theta-1 \leq k-1$ and degree $k$. TSS under the LTM remains NP-hard under various assumptions on the graph and its thresholds \citep{chen2009approximability,ben2011treewidth,centeno2011irreversible, nichterlein2013tractable,chopin2014constant}. However, none of the existing results cover our specific case. Indeed, the proofs of these results rely on reductions from vertex cover and require the construction of an instance where each node's threshold and degree are equal. This reduction is not feasible in our case due to our structural assumptions (nodes in $N_{SC}$ have degree $k$ and threshold $\theta-1 \leq k-1$), so we must resort to a more complex proof. Thus, Theorem \ref{thm:comp} also contributes to the LTM literature, more specifically, to understanding which problem structures make TSS hard to solve. This is difficult to determine a priori. For example, TSS is fixed-parameter tractable with respect to treewidth, which tends to be low for sparse graphs, \emph{and} with respect to cluster edge deletion number, which tends to be low for dense graphs \citep{chen2009approximability,nichterlein2013tractable}. 

\vspace{0.2cm}

Not only can we not answer $DEC$-$SCTM$ in polynomial time if $k \geq \theta \geq 3$, we cannot provide a meaningful approximation of the true solution under a slightly stronger assumption than $P \neq NP$:

\begin{proposition} \label{prop:eps-inapprox}
For any $k \geq \theta \geq 3$, there exists an $\alpha > 1$ such that, unless $NP \subseteq DTIME(n^{polylog(n)})$, the optimal value to \eqref{eq:tss.sctm} with $r_{ji}=c_i=1\ \forall i=0,\ldots,n,\ j=1,\ldots,m$ cannot be approximated in polynomial time within the ratio of $O\left( \alpha^{\log^{1-\xi} n} \right)$ for any fixed constant $\xi > 0$.
\end{proposition}

We have shown that $MIN$-$SCTM$ is hard to solve or approximate when $k \geq \theta \geq 3$, even if the adoption costs are negligible compared to the traceability benefits. These results suggest that complexity is not mainly driven by costs and benefits but rather by the structure of $G$ and the supply chain effect. Any hope of obtaining a polynomial-time algorithm for $MIN$-$SCTM$ must rely on assuming additional structure. We discuss this next.

\section{An Exact Solution Algorithm for \texorpdfstring{$MIN$-$SCTM$}{MIN-SCTM}} \label{sec:exact.sols}

We now provide a linear programming-based \emph{fixed-parameter tractable} algorithm for $MIN$-$SCTM$ with respect to the \emph{treewidth} of $G$. Precise definitions of these concepts are in Section \ref{subsec:tw.fpt}, but at a high level, this means that the complexity of solving $MIN$-$SCTM$ is significantly reduced when the treewidth of the supply chain network is small. In light of this, Section~\ref{subsec:tw.fpt} provides evidence that real-world supply chain networks have small treewidth. Thus, solving $MIN$-$SCTM$ on real-world networks is not as hopeless as one may be led to believe by the results in Section \ref{sec:Comp.Comp}.

Moving forward, Section \ref{subsec:BiLP.Reform} builds an integer programming formulation of $MIN$-$SCTM$ that exploits tree decompositions of $G'$. While the formulation allows for solving the problem in many practical instances, we go further and leverage the formulation to derive an LP-based FPT algorithm in Section \ref{subsec:LP.FPT}. The techniques we use come from the integer programming literature and, to our knowledge, have never been applied in the context of network diffusion. Finally, Section \ref{sec:lb} introduces a hierarchy of LPs that directly trades off computational complexity and approximation quality to obtain lower bounds on the optimal value of $MIN$-$SCTM$. We also use this hierarchy of LPs to obtain upper bounds (and corresponding feasible sets).

\subsection{Treewidth and FPT algorithms} \label{subsec:tw.fpt}

We define the concept of treewidth \citep{bodlaender1994tourist}, central to the rest of the paper.

\begin{definition} \label{def:tw}
Let $G$ be a (hyper)graph with nodes $N$ and (hyper)edges $E$. A \emph{tree decomposition} of $G$ is a pair $\mathcal{T}=(T,\{X_z\}_{z \in T})$, with tree $T$ and bags $X_z$ for each node $z\in T$, such that:
\begin{enumerate}[(a)]
    \item $\bigcup_{z \in T}X_z=N$.
    \item If $\{i_1,\ldots,i_h\} \in N$ belong to (hyper)edge $e \in E$, there must be a set $X_{z}$ with $\{i_1,\ldots,i_h\} \in X_z$.
    \item If a node $i \in N$ appears in two distinct bags $X_x$ and $X_y$, then it appears in all bags $X_z$ such that $z$ is on the (unique) path between $x$ and $y$ in $T$.
\end{enumerate}
The tree decomposition's \emph{width} is $\max_{z \in T} |X_z|-1$. The \emph{treewidth} $tw(G)$ of $G$ is simply the minimum width over all tree decompositions of $G$.
\end{definition}
A (hyper)graph $G$ always admits a trivial tree decomposition with a single node $T$ containing $N$, so $tw(G) \leq n-1$. However, $tw(G)$ is much smaller when $G$ is ``tree-like''\mdash trees have treewidth~1.

\begin{figure}[t]
    \centering
    \includegraphics{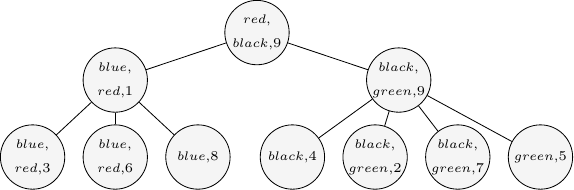}
    \vspace{0.2cm}
    \caption{Tree decomposition of the auxiliary graph in Figure \ref{fig:aux_graph_aux}. Each bag in the tree corresponds to a set of nodes from the auxiliary graph.}\label{fig:treedec}
\end{figure}

Two graphs play important roles in our model: the original hypergraph $G$ and the auxiliary graph $G'$. Figure \ref{fig:treedec} displays a tree decomposition of $G'$ from Figure \ref{fig:aux_graph_aux} with treewidth two. The uppermost bag contains Firm-node $9$ and SC-nodes ${red}$ and ${black}$. Following (b), as SC-node ${red}$ also appears in a bag at the bottom left, it must appear at the intermediate level on the left-hand side. A natural question is how the treewidth of $G$ and $G'$ relate: this is the focus of our next result.

\begin{proposition} \label{prop:courcelle}
Let $G$ be a hypergraph as defined in Section \ref{subsec:model.description} and let $G'$ be its auxiliary graph as defined in Section \ref{sec:aux_graph}. Then, $tw(G') \leq tw(G)+1.$
\end{proposition}
The proof of Proposition \ref{prop:courcelle} is in Appendix \ref{appendix:LP.ackerman}. Our analysis uses a tree decomposition of $G'$ with treewidth $tw(G')=\omega'$, so we related our results to $tw(G)=\omega$ using the upper bound on $\omega'$.

Another central concept to our work is \emph{fixed-parameter tractability} \citep{downey2012parameterized}:

\begin{definition}\label{def:fpt}
A problem parameterized by $\theta$ is \emph{fixed-parameter tractable (FPT)} with respect to $\theta$ if it can be solved in $f(\theta) n^{O(1)}$ time, where the function $f$ does not depend on $n$. 
\end{definition}

The algorithms we develop are fixed-parameter tractable in the treewidth of $G$. Thus, they are particularly valuable in settings where the treewidth of the supply chain network is small.
 Data from \citet{Willems2008} provides evidence that supply chain networks have small treewidth. The data represents 38 acyclic network structures gathered from companies in 22 industries, from which we randomly generate 657 supply chain networks. The generation process details are in Appendix~\ref{appendix:data}. Figure \ref{fig:hists} displays upper bounds on the auxiliary graph treewidths\footnote{To compute a minimal tree decomposition, we use the Flow Cutter algorithm from the \emph{PACE 2017 Parameterized Algorithms and Computational Experiments Challenge} \citep{dell2018pace}.} relative to the networks' sizes  for this dataset. We focus on the auxiliary graph treewidth $\omega'$, rather than the possibly larger hypergraph treewidth $\omega$, as our results (including Corollary \ref{cor:lp}) hold for $\omega'$. We observe that $\omega'$ is much smaller than the supply chain network size $n+m$. The mean (resp.\ median) treewidth is 9 (resp.\ 6), which is less than 3\% of the mean size of 310 (resp.\ 151). More broadly, we note that the supply chain management literature frequently assumes supply chain network graphs to be trees \citep[see, e.g.,][]{graves2000optimizing}, which have a treewidth of one.

\begin{figure}[t]
    \centering
    \begin{adjustbox}{clip,trim=0cm 0.75cm 0cm 0cm,max width=\textwidth}
    \includegraphics{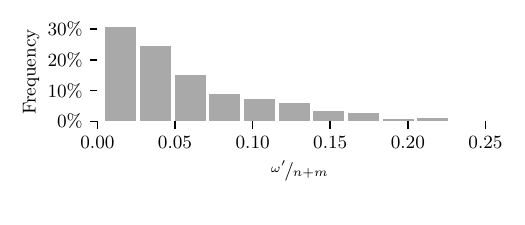}
    \label{fig:node_hist}
    \end{adjustbox}
    \caption{Histogram of $\frac{\omega'}{n+m}$ for 657 randomly generated supply chain networks based on \citet{Willems2008}.}\label{fig:hists}
\end{figure}

\subsection{Binary Linear Programming Reformulations of \texorpdfstring{$MIN$-$SCTM$}{MIN-SCTM}} \label{subsec:BiLP.Reform}

The goal of this section is to reformulate $MIN$-$SCTM$ as a binary linear program (BiLP), that is, an optimization problem of the following form:
\begin{equation}
\begin{aligned} 
\min_{x \in \{0,1\}^n} &c^Tx\\
\text{s.t. } & Ax \geq b,
\end{aligned}
\label{eq:binLP}
\end{equation}
where $c \in \mathbb{R}^n, A \in \mathbb{R}^{m \times n}$, and $b \in \mathbb{R}^m$. The next result clarifies why this could be of interest.

\subsubsection*{Linear Programming Formulations of BiLPs.} We first define an \emph{intersection graph}:

\begin{definition}
The intersection graph of \eqref{eq:binLP} is an undirected graph with a vertex for each variable $x_i,~i=1,\ldots,n$ and an edge for each pair $(x_i,x_j)$ that feature in the same constraint.
\end{definition}

\begin{proposition}[\citealt{laurent2009sums, bienstock2018lp}]\label{prop:binLP}
If the intersection graph of (\ref{eq:binLP}) has treewidth $\omega$, then there is an equivalent reformulation of (\ref{eq:binLP}) as a linear program with $O(2^\omega n)$ variables and constraints.
\end{proposition}
Since LPs with $h$ variables and constraints can be solved in $O(h^{2.5})$ time \citep{jiang2020faster}, solving this LP is an FPT algorithm for (\ref{eq:binLP}) with parameter treewidth of the BiLP's intersection graph. Thus, our goal moving forward is to reformulate $MIN$-$SCTM$ as a BiLP, whose intersection graph's treewidth is upper-bounded by a function involving $\omega'$. From Proposition \ref{prop:binLP}, one can construct an equivalent linear program, which constitutes our FPT algorithm for $MIN$-$SCTM$. As mentioned, an LP approach of this type is new to the network diffusion literature. 

\subsubsection*{A First BiLP Formulation of $MIN$-$SCTM$.}  \citet{ackerman2010combinatorial} provide a BiLP formulation of the target set selection problem in the linear threshold model, which\mdash to the best of our knowledge\mdash is the only such formulation in the literature. By virtue of Corollary~\ref{cor:sctm.ltm}, we adapt the formulation in \citet{ackerman2010combinatorial} to our setting as follows. Let 
$$E_{F,SC} = \{ (i,j) \in E': i \in N_F, j \in N_{SC} \}  \text{ and } E_{SC,F} = \{ (j,i) \in E': i \in N_F, j \in N_{SC} \}.$$

\begin{proposition}[Adapted from \citealt{ackerman2010combinatorial}]\label{prop:ackerman}
Let $G'=(N',E')$ represent the auxiliary graph of hypergraph $G$. Consider the following BiLP:
\begin{subequations} \label{eq:MILP}
\begin{alignat}{2}
&\min & \quad  &\sum_{i \in N_F} w_i s_i  \notag \\
&\text{s.t.} &  &\sum_{\{j~|~(j,i) \in E_{SC,F} \} } r_{ji} \ell_{ji} \geq c_i (1 - s_i),  ~\forall i \in N_F,  \label{eq:milp.const1} \\
& & &\sum_{\{i~|~(i,j) \in E_{F,SC}\}} \ell_{ij} \geq \theta_j - 1, ~\forall j \in N_{SC}, \label{eq:milp.constr2}\\
& & &\ell_{ij} + \ell_{ji} = 1,  ~\forall i \in N_F, \forall j \in N_{SC}, \label{eq:milp.constr3}\\
&  & &\ell_{i_1j_1}+\ell_{j_1i_2}+\ell_{i_2j_2}+\ell_{j_2i_1} \leq 3,~\forall i_1,i_2 \in N_F, ~\forall j_1,j_2 \in N_{SC}\label{eq:milp.constr4},\\
& & &s_i \in \{0,1\},~\forall i\in N_F,~ \ell_{ij}, \ell_{ji} \in \{0,1\},  ~\forall i \in N_F,~\forall j \in N_{SC}.\notag
\end{alignat}
\end{subequations}
The set $S_0^* = \{ i \in N_F : s_i^* = 1\}$ is a solution to the problem $MIN$-$SCTM$ on hypergraph $G$.
\end{proposition}
The proof of this and the following result are in Appendix \ref{appendix:LP.ackerman}. The BiLP constructs a directed acyclic graph (DAG) on $N'$, encoding the activation sequence of $G'$. In other words, there is an edge ($\ell_{ij} = 1$) between node $i,j \in N'$ if node $i$ contributes to node $j$'s activation. The seed set nodes are then simply the sources of the DAG and are encoded via the variables $s_i,~i \in N_F$. Constraints \eqref{eq:milp.const1} and \eqref{eq:milp.constr2} enforce that activation can only proceed as described in Section \ref{sec:SCTM}, whereas constraints \eqref{eq:milp.constr3} and \eqref{eq:milp.constr4} ensure that the final graph is, respectively, directed and acyclic.

Although \eqref{eq:MILP} solves $MIN$-$SCTM$, the treewidth of its intersection graph is at least $nm-1$, which precludes us from using Proposition \ref{prop:binLP} to derive an FPT algorithm from this BiLP. 

\begin{proposition} \label{prop:tw.benzwi}
The treewidth of the intersection graph of (\ref{eq:MILP}) is at least equal to $nm-1$. 
\end{proposition}
Proposition \ref{prop:tw.benzwi} makes it clear that one cannot simply use any BiLP formulation of $MIN$-$SCTM$ to leverage Proposition \ref{prop:binLP}: great care must be taken to formulate the BiLP appropriately.

\subsubsection*{A Second BiLP Formulation of $MIN$-$SCTM$.} The next formulation is closely tied to a tree decomposition $\mathcal{T}'=(T',\{X'_z\}_{z \in T'})$ of $G'$ with treewidth $\omega'$. For the remainder of this section, we assume wlog that $T'$ is binary. That is, each node $z \in T'$ has at most two children, $c_1(z)$ and $c_2(z)$.\footnote{If $T'$ is not binary, we choose an arbitrary node as the root and proceed top-down. If a node $z$ has $\tilde{n}>2$ children $c_1(z),\ldots,c_{\tilde{n}}(z)$, we create a new node $z'$ with bag $X'_{z'}=X_z'$ and add $z'$ and $c_1(z)$ as children of $z$ and $c_2(z), \ldots, c_{\tilde{n}}(z)$ as children of $z'$, and repeat as necessary. This process returns a binary tree whose bags are exactly those in $G'.$} To avoid the issues described in Proposition \ref{prop:tw.benzwi}, we first replace constraints (\ref{eq:milp.constr3}) and (\ref{eq:milp.constr4}) by constraints that are tree decomposition-dependent. We then obtain
\begin{subequations} \label{eq:MILP_Inter2}
\begin{alignat}{2}
&\min & \quad  &\sum_{i \in N_F} w_i s_i  \notag \\
&\text{s.t.} &  &\sum_{(j,i) \in E_{SC,F}} r_{ji} \ell_{ji} \geq c_i (1 - s_i),~\forall i \in N_F , \label{eq:milpInter2.const1} \\
& & &\sum_{(i,j) \in E_{F,SC}} \ell_{ij} \geq \theta_j - 1, ~\forall j \in N_{SC}, \label{eq:milpInter2.constr2}\\
& & &\ell_{ij} + \ell_{ji} = 1 , ~\forall i,j \in N' \cap X_z', \forall z \in T', \label{eq:milpInter2.constr3}\\
&  & &\ell_{ij}+\ell_{jk}+\ell_{ki} \leq 2,~\forall i,j,k \in N' \cap X'_z, \forall z \in T', \label{eq:milpInter2.constrend}\\
& & &s_i \in \{0,1\},~\forall i\in N_F,~ \ell_{ij}, \ell_{ji} \in \{0,1\},  ~\forall i, j \in  N' \cap X_z', \forall z \in T'. \notag
\end{alignat}
\end{subequations}
As it turns out, (\ref{eq:MILP_Inter2}) and (\ref{eq:MILP}) are equivalent, which we formally state next.

\begin{proposition} \label{prop:eq.10}
Problems \eqref{eq:MILP_Inter2} and \eqref{eq:MILP} are equivalent. In particular, let $(s_i^*,\ell_{ji}^*, \ell_{ij}^*)$ be a solution to (\ref{eq:MILP_Inter2}). The set $S_0^*=\{i \in N_F~|~s_i^*=1\}$ is a solution to the problem $MIN$-$SCTM$ on hypergraph $G$.
\end{proposition}
The proof is in Appendix \ref{appendix:LP.fpt}. Although (\ref{eq:MILP}) and (\ref{eq:MILP_Inter2}) are equivalent, (\ref{eq:MILP_Inter2}) can have much fewer constraints than (\ref{eq:MILP}) if the treewidth of $G'$ is small. To see this, note that (\ref{eq:MILP}) has $O\left(n^2 m^2\right)$ constraints, driven by \eqref{eq:milp.constr4}, which makes the activation sequence acyclic. In contrast, since there always exists a (binary) tree decomposition $T'$ of $G'$ with at most $4 (n+m)$ bags  \citep[see Lemma 13.1.2. of][]{kloks1994treewidth} and since $T'$ has bags of size at most $\omega + 1$, Formulation \eqref{eq:MILP_Inter2} has $O\left( (n+m) \cdot \omega^3\right)$ constraints. When $\omega$ is much smaller than $m$ and $n$, (\ref{eq:MILP_Inter2}) will be a much smaller optimization problem than (\ref{eq:MILP}). 

This difference directly translates to a difference in solving time for \eqref{eq:MILP} and \eqref{eq:MILP_Inter2}. To illustrate, we draw a sample of 150 supply chain networks with $m+n \geq 100$ from the set of networks described in Section \ref{subsec:tw.fpt}, for which we can solve \eqref{eq:MILP} and \eqref{eq:MILP_Inter2} within three hours using \citet{gurobi} on a computing cluster with 24 cores. The average time for solving \eqref{eq:MILP_Inter2} on these instances is 50\% less than that for \eqref{eq:MILP}, even when accounting for the calculation time of the tree decomposition of $G'$. When the decomposition is not accounted for, solving \eqref{eq:MILP_Inter2} takes 72\% less time than solving \eqref{eq:MILP} on average. For several instances, we can solve  \eqref{eq:MILP_Inter2} to optimality within three hours but not \eqref{eq:MILP}.

Despite the encouraging computational results, it is not clear that the treewidth of the intersection graph of \eqref{eq:MILP_Inter2} can be upper-bounded by a function of $\omega'$. In fact, the intersection graph's treewidth can be as large as $m+1$. To see this, suppose that $i \in N_F$ belongs to $m$ supply chains. Then, constraint (\ref{eq:milpInter2.constr2}) leads to a clique of size $m+1$ in the intersection graph, containing $\{\ell_{ji}\}_{(j,i) \in E_{SC,F}}$ and $s_i$. This implies that the treewidth is at least $m+1$, and we can still not use Proposition \ref{prop:binLP} to derive an LP-based FPT algorithm. Hence, we next introduce new variables representing partial sums of terms appearing in constraints \eqref{eq:milpInter2.const1} and \eqref{eq:milpInter2.constr2}. The idea is to replace a constraint such as $v+w+y+z \geq x$ with three equivalent constraints $x_1+x_2 \geq x$, $v+w \geq x_1$, $y+z \geq x_2$, leading to an intersection graph of smaller treewidth, at the expense of more constraints and variables. This is similar in spirit to \citet{bienstock2018lp}. However, our formulation explicitly leverages the tree decomposition of $G'$ as well as the specificities of our problem.

\subsubsection*{A Third (and Final) BiLP Formulation of $MIN$-$SCTM$.}
Before proceeding, we introduce some new notation. For $i \in N_{F}$ and $j \in N_{SC}$, we let $J^i=\{j~|~(j,i)\in E_{SC,F}\}$ and $I^j=\{i~|~(i,j) \in E_{F,SC}\}$ be the sets that appear in constraints (\ref{eq:milpInter2.const1}) and (\ref{eq:milpInter2.constr2}). As mentioned above, our goal is to group the variables appearing in these constraints into partial sums in an effective way. We do this by leveraging $T'$, with the idea that variables that appear in the same partial sum should have indices in the same bag in $T'$ so that the resulting intersection graph has low treewidth. We describe how to do this based on a constraint in (\ref{eq:milpInter2.const1}). Let $i \in N_F$: we look for a \emph{partition} of $J^i$ into sets $\{J_z^i\}$ such that $J_z^i \subseteq X'_z$ for $z \in T'$. We can then rewrite $\sum_{j \in J^i} r_{ji} \ell_{ji}= \sum_{z}  \left(\sum_{j \in J_{z}^i} r_{ji} \ell_{ji}\right)$, where $\sum_{j \in J_{z}^i} r_{ji} \ell_{ji}$ are partial sums with the variables' indices all belonging to the same bag.

Let $T'_i$ be the subtree of $T'$ when restricted to bags containing $i$, $z_0^i$ an arbitrarily chosen root node of $T'_i$, and $|T'_i|$ the number of tree nodes in $T'_i.$ By construction, this tree's bags contain all of $J^i$, which we next partition. We split tree nodes of $T_i'$ into two sets: a ``useful'' set $T_i^G$ to build the partition, and a ``useless'' set $T_i^{\tilde{G}}$ to keep track of over-counting. Set $T_i^G$ is obtained by sequentially adding tree nodes from $T_i'$ while ensuring that the associated bag contains at least one $j$ that is not already present in the bags of $T_i^G$. We stop when $T_i^G$ contains all of $J^i$ in its bags. The set $T_i^{\tilde{G}}$ contains the remaining tree nodes of $T'_i$. We are now ready to build our partition $\{J_z^i\}_{z \in T_i^G}$ of $J_i$: we arbitrarily number the bags $\{X_z'\}_{z \in T_i^G}$ and let $J_z^i=X_z' \cap N_{SC}$ for the first tree node. For other tree nodes, we let $J_z^i = X_z' \cap N_{SC} \setminus \{j \in J_i~|~j\text{ present in previous bags} \}$. For constraints \eqref{eq:milpInter2.constr2}, we similarly define $T'_j$ as the subtree of $T'$ when restricted to bags containing $j$, with $z_0^j$ an arbitrary root of $T'_j$ and $|T'_j|$ the number of tree nodes in $T'_j$. Within $T'_j$, we construct analogous concepts $T_j^G$ and $T_j^{\tilde{G}}$, replacing $J^i$ by $I^j$. We also let $\{I_{z}^j\}_{z \in T_j^G}$ be the counterpart of $\{J_z^i\}_{z \in T_j^G}$. Letting 
\begin{align} \label{eq:nu.nv}
n_u=\lfloor \log_2(\max\{c_1,\ldots,c_n\}+1) \rfloor -1 \text{ and } n_v=\lfloor \log_2(\max\{\theta_1,\ldots,\theta_m\}+1) \rfloor -1,
\end{align}
we can formulate our final BiLP. As this BiLP is quite cumbersome to write out, we place it in Appendix \ref{appendix:LP.fpt} as equation \eqref{eq:MILP_New}. Note that \eqref{eq:MILP_New} is not precisely in the form of \eqref{eq:binLP} to keep the formulation legible. However, it can be obtained by appropriately substituting variables using equations \eqref{eq:milpNew.binu}, \eqref{eq:milpNew.binv}, and \eqref{eq:milpNew.equality1} in constraints \eqref{eq:milpNew.binu} and \eqref{eq:milpNew.binv}. This reformulation is a common technique for replacing integer variables with binary ones \citep[see, e.g.,][]{watters1967reduction}. Furthermore, we have not differentiated between nodes with different numbers of children in the constraints. With a slight abuse of notation, we assume that if  $c_1(z)= \emptyset$ and/or $c_2(z)= \emptyset$, the corresponding sum is dropped. We now show two results regarding (\ref{eq:MILP_New}), with proofs in Appendix \ref{appendix:LP.fpt}. These enable us to propose an LP-based FPT algorithm for $MIN$-$SCTM$ in the next section.

\begin{proposition} \label{prop:eq.MIN.SCTM}
Optimization problems (\ref{eq:MILP_Inter2}) and (\ref{eq:MILP_New}) are equivalent. In particular, let $\left(s_i^*,\ell_{ij}^*, u_{iz}^{r*}, \tilde{u}_{iz}^{r*}, v_{jz}^{r*}, \tilde{v}_{jz}^{r*}\right)$ be an optimal solution to (\ref{eq:MILP_New}). The set $S_0^*=\{i \in N_F~|~ s_i^*=1\}$ is a solution to the $MIN$-$SCTM$ problem on hypergraph $G.$
\end{proposition}

\begin{proposition} \label{prop:tw.MILP.New}
Let $\vartheta_{\max} = \max \left\lbrace \max_{i=1,\ldots,n} \{ c_{i}\}, \max_{j=1,\ldots,m} \{ \theta_{j}\} -1 \right\rbrace$.
The treewidth of the intersection graph of (\ref{eq:MILP_New}) is at most $O(\omega'^2+\omega' \log_2(\vartheta_{\max}))$, where $\omega'$ is the treewidth of $G'.$
\end{proposition}

\subsection{A Linear Programming FPT Algorithm for \texorpdfstring{$MIN$-$SCTM$}{MIN-SCTM}} \label{subsec:LP.FPT}
We use Propositions  \ref{prop:courcelle}, \ref{prop:eq.MIN.SCTM}, and \ref{prop:tw.MILP.New} to show one of our main results.

\begin{theorem} \label{thm:eq.refor}
Let $G$ be a hypergraph with $tw(G)=\omega$. Then, there is an equivalent reformulation for $MIN$-$SCTM$ as a linear program with a number of constraints and variables in 
{\small 
\begin{align*}
O\left(2^{\omega^2}\cdot \vartheta_{\max}^{\omega}  \left(n+(n+m)\omega^3+8 \log_2\left(\max \{c_1,\ldots,c_n\}\right) \left(n+m\right) + \log_2 \left(\max \{\theta_1,\ldots,\theta_m\}\right)  \left(n+m\right) \right)\right).
\end{align*}}
\end{theorem}
The proof of Theorem \ref{thm:eq.refor} can be found in Appendix \ref{appendix:LP.fpt}. Theorem \ref{thm:eq.refor} indicates that an LP reformulation of $MIN$-$SCTM$ with $2^{O(\omega^2)} \vartheta_{\max}^{O(\omega)} O(n+m)$ constraints and variables exists. Following \citet{jiang2020faster}, we then have the following:

\begin{corollary} \label{cor:lp}
There is a linear programming-based FPT algorithm for solving $MIN$-$SCTM$ with parameters $\omega$ and $\vartheta_{\max}$, running in time at most
$2^{O(\omega^2)} \vartheta_{\max}^{O(\omega)} O \left((n+m)^{2.5}\right).$
\end{corollary}
Up to this point, we emphasized that our FPT algorithms are with respect to the treewidth $\omega$. Technically speaking, they are FPT algorithms with respect to the parameters $\omega$ and $\vartheta_{\max}$. However, one can reasonably assume a bound, independent of $n$ and $m$, on $c_i$ and $\theta_j$, for $i=1,\ldots,n$ and $j=1,\ldots,m$, and so on $\vartheta_{\max}$. This is because $c_i$ is a firm-dependent parameter that describes the adoption costs. At the same time, $\theta_j$ is a ``local" supply chain-dependent parameter, describing interactions between a fixed set of firms that does not change as the network scales.\footnote{Such an assumption does not preclude $MIN$-$SCTM$ from being a hard problem to solve and approximate as evidenced in Section \ref{sec:Comp.Comp} where $\vartheta_{\max}$ is equal to 2.} Under this assumption, one can remove the dependency of the algorithm's runtime on $\vartheta_{\max}$.

We now present the construction of the LP that Corollary \ref{cor:lp} relies on. This formulation differs from the one in, e.g., \citet{bienstock2018lp}. The main advantage of our formulation is that one can easily use it to derive smaller LPs that provide good-quality lower bounds for $MIN$-$SCTM$ (see below for a detailed discussion). It is based on the tree decomposition of the intersection graph of \eqref{eq:MILP_New}, which we denote by $\mathcal{S}=(S,\{W_z\}_{z \in S})$. We also let $\omega_z=|W_z|$. Each bag $W_z$ contains variables that we rename $x_1^z,\ldots,x_{\omega_z}^z$, which are subsets of the decision variables in \eqref{eq:MILP_New}. Furthermore, we associate a set of $l_z$ constraints from \eqref{eq:MILP_New} with each bag $W_z$: the constraints that only feature variables in $W_z$. As these are linear in the decision variables, we can write them as $g_{\emptyset}^l +\sum_{i=1}^{\omega_z} g_{\{x_{i}^z\}}^l \cdot x_i^z \geq 0$ for $l=1,\ldots,l_z.$ We are ready to state our LP formulation of $MIN$-$SCTM$:
  \begin{equation} \label{eq:LP}
    \begin{aligned}
    &\min_{Y_\mathbb{S}} &&\sum_i w_i Y_{\{s_i\}}\\
    &\text{s.t. } &&Y_{\emptyset}=1, ~~ \sum_{\{\mathbb{S} \in 2^{W_z}~|~ \mathbb{T} \subseteq \mathbb{S}\}} (-1)^{|\mathbb{S}|-|\mathbb{T}|} Y_{\mathbb{S}} \geq 0, ~\forall \mathbb{T} \in 2^{W_z},~\forall z \in S,\\
    & &&\sum_{\{\mathbb{S} \in 2^{W_z}~|~\mathbb{T} \subseteq \mathbb{S}\}} (-1)^{|\mathbb{S}|-|\mathbb{T}|} \left( g_{\emptyset}^l Y_{\mathbb{S}}+ \sum_{i=1}^{\omega_z} g_{\{x_{i}^z\}}^l \cdot Y_{\{x_i^z\} \cup \mathbb{S}}\right) \geq 0,~\forall \mathbb{T} \in 2^{W_z},~\forall l=1,\ldots,l_z,\forall z \in S,
    \end{aligned}
   \end{equation}
   where $2^{W_z}$ corresponds to all possible subsets of variables in $W_z$. For small values of $\omega_z$, \eqref{eq:LP} can be solved exactly as commercial solvers allow for millions of variables and constraints. It can be the case, however, that when $\omega_z$ becomes larger, this problem becomes difficult to solve due to memory constraints. One advantage of the LP approach is that it provides a principled way of deriving less computationally-intense heuristics, as seen in Section \ref{sec:lb}.
   
\subsubsection*{Comparison to a Dynamic Programming-Based FPT Algorithms.} We also propose a dynamic programming (DP)-based FPT algorithm for $MIN$-$SCTM$ in Appendix~\ref{appendix:dp}, which is a non-trivial generalization of an algorithm introduced by \citet{ben2011treewidth} for the target set selection problem in the LTM. We derive an analogous statement to Corollary \ref{cor:lp} for this algorithm and show that its runtime is at most $2^{O(\omega \log_2(\omega))} \vartheta_{\max}^{O(\omega)} (n+m)$. Thus, theoretically, the DP-based algorithm has a slightly lower run time than the LP-based algorithm.\footnote{Recall the definition of an FPT algorithm in Section \ref{subsec:tw.fpt}. The quality of an FPT algorithm is based on how small one can make $f(\omega)$, rather than the exponent of $(n+m)$ \citep[see, e.g.,][]{lokshtanov2011known}. This is important here as our LP-based algorithm is close to its DP counterpart for $f(\omega)$ but not as close for the exponent of $(n+m)$.} However, it has the usual caveat of being difficult to implement, requiring specific coding of the algorithmic procedure. Such an ad-hoc approach typically does not optimize for memory or processing capabilities. In contrast, widely available commercial solvers can solve linear programs with thousands of variables and millions of constraints in a matter of minutes. One can further speed up the solving time by using techniques such as set-up parallelization and warm-starting, which are tried-and-tested techniques for LPs. Last but not least, unlike the DP-based algorithm, the LP-based algorithm allows us to systematically derive good-quality lower and upper bounds on $MIN$-$SCTM$. We see this now.

\subsection{Bounds from a Hierarchy of Linear Programs.} \label{sec:lb}
We now derive principled upper and lower bounds on the optimal value of $MIN$-$SCTM$.

\subsubsection*{Lower Bounds on the Optimal Value of $MIN$-$SCTM$.} Any seed set leading to full activation of the graph directly implies an upper bound on the optimal value of $MIN$-$SCTM$. Obtaining lower bounds, on the other hand, that go above and beyond the simple lower bound obtained by considering the continuous relaxation of (\ref{eq:MILP_New}), or equivalently (\ref{eq:MILP_Inter2}), can be trickier. Problem (\ref{eq:LP}) provides us with a process to generate increasingly powerful lower bounds via a \emph{hierarchy} of LPs, that is, a family $\{LP_{\kappa}\}_{\kappa=1,\ldots,\omega}$, where $LP_{\omega}$ is equal to (\ref{eq:LP}) (i.e., $LP_{\omega}$ solves $MIN$-$SCTM$ exactly) and where the objective value of $LP_{\kappa}$ is an increasingly tight lower bound on the objective value of (\ref{eq:LP}) as $\kappa$ grows. Recall the notation given in Section \ref{subsec:LP.FPT}. The linear program $LP_{\kappa}$ is given thus:
\begin{equation} \label{eq:LP.k}
    \begin{aligned}
    \min_{Y_\mathbb{S}}~~~& \sum_i w_i Y_{\{s_i\}}\\
    \text{s.t.}~~~&Y_{\emptyset}=1,\\
    & \sum_{\{\mathbb{S} \in 2^{W_z} ~|~ \mathbb{T} \subseteq \mathbb{S} \subseteq \mathbb{U}\}} (-1)^{|\mathbb{S}|-|\mathbb{T}|} Y_{\mathbb{S}} \geq 0, ~\forall \mathbb{T}, \mathbb{U} \in 2^{W_z}, \mathbb{T} \subseteq \mathbb{U}, |\mathbb{U}| = \min\{\kappa+1,n\},\forall z \in S, \\
    & \sum_{\{\mathbb{S} \in 2^{W_z}~|~ \mathbb{T} \subseteq \mathbb{S} \subseteq \mathbb{U}\}} (-1)^{|\mathbb{S}|-|\mathbb{T}|} \left( g_{\emptyset}^lY_{\mathbb{S}}+ \sum_{i=1}^{\omega_z} g_{\{x_{i}^z\}}^l \cdot Y_{\{x_i^z\} \cup \mathbb{S}}\right) \geq 0, ~\forall \mathbb{T}, \mathbb{U} \in 2^{W_z},\mathbb{T} \subseteq \mathbb{U}, |\mathbb{U}| =\kappa,\\
    &  \hspace{9.75cm} \forall l=1,\ldots,l_z,\forall z \in S.
\end{aligned}
\end{equation}
The objective value of (\ref{eq:LP.k}) increases with $\kappa$. When $\kappa=\omega$, it equals the solution to (\ref{eq:LP}) \citep[see][]{laurent2003comparison}. The LP's size also increases in $\kappa$, providing an explicit trade-off between accuracy and computation time. Interestingly, when $\kappa=0$, we obtain the simple lower bound mentioned above.

\begin{proposition} \label{prop:kappa}
When $\kappa=0$, \eqref{eq:LP.k} is equivalent to \eqref{eq:MILP_New}, where the binary variables have been replaced by continuous variables on $[0,1]$.
\end{proposition}

Experimentally, we compute lower bounds for a sample of 100 supply chain network instances with $m \geq 50$.\footnote{We use a dataset similar to the \citet{Willems2008} networks introduced in Appendix \ref{appendix:data}. However, to make solving more challenging, we only standardize $r_{ji}=1$ and randomly vary $w_i$, $c_i$, and $\theta_j$ for all $i \in N_F$, $j \in N_{SC}$.} We measure the relative gap between the costs when solving $LP_\kappa$, and an upper bound obtained when attempting to solve \eqref{eq:MILP_Inter2} using Gurobi, as described in Section \ref{subsec:BiLP.Reform}. When $\kappa=0$, the median gap (resp.\ IQR) between our lower bound and the upper bound is -51\% (resp.\ -62\% -- -42\%). When $\kappa=1$, the median gap dramatically improves to -19\% (resp.\ -31\% -- -11\%), indicating that for small values of $\kappa$, we already obtain powerful lower bounds that considerably outperform those that can be obtained from a simple relaxation.

\subsubsection*{Upper Bounds on the Optimal Value of $MIN$-$SCTM$.}

We can further leverage solutions to \eqref{eq:LP.k} to obtain upper bounds on the optimal value of $MIN$-$SCTM$ and corresponding feasible seed sets. The most direct approach would be to sort nodes according to their score $Y_{\{s_i\}}^*$ from the solution. We could then add them to the seed set one-by-one until the entire network activates. However, this does not consider supply chain dependencies and has poor numerical performance. Noting that activation proceeds along supply chains and that some seed nodes are only relevant late in the activation process, our algorithm instead sequentially selects supply chains based on their constituents' scores. This sequential process is motivated by our insights in Section \ref{sec:insight2}.

The heuristic, formalized in Algorithm \ref{alg:ub}, involves three steps. First, for fixed $\kappa$, it solves \eqref{eq:LP.k}. Second, it uses the optimal solution $Y_{\{s_i\}}^* \in [0,1]$, as a ``score" for each firm $i$. The heuristic chooses the (relevant) subset of nodes with the highest average score for each supply chain. Namely, suppose a given supply chain requires $h$ nodes for full activation and has inactive nodes $I$. Then, the heuristic examines $\binom{|I|}{h}$ combinations of nodes in that supply chain and chooses the combination with the highest average score as a candidate for the seed set. Finally, the heuristic compares the average scores of candidate combinations across supply chains, greedily choosing the combination (and, thus, the next supply chain to become active) with the highest average.

\begin{algorithm2e}[!ht]
 \caption{An LP relaxation-based upper bound. \label{alg:ub}}
 \SingleSpacedXI
\KwData{A hypergraph $G=(N_F,E)$ as defined in Section \ref{subsec:model.description}.}
\KwResult{A feasible seed set $S_0$ of $G$.}
Initialize: $S_0 = \emptyset$; $\{Y_{s_i}^*\}_{i \in N_F} \leftarrow$ solution to $LP_\kappa$ applied to $G$\;
\While{$x_{i,\infty} = 0$ for some $i \in N_F$, given $S_0$ as seed set}{
    $s_{\max} \leftarrow 0$;
    $O_{\max} \leftarrow \emptyset$\;
    \For{$e_j \in E$}{
        $I \leftarrow \{i \in e_j : x_{i,\infty} = 0 \text{ given } S_0 \text{ as seed set}\}$\;
        $h \leftarrow \theta_j - 1 - \sum_{i \in e_j} x_{i,\infty}$\;
        \For{$O \in choose(I,h)$}{
            $s \leftarrow \frac{1}{|O|} \sum_{i \in O} Y_{s_i}^*$\;
            \If{$s \geq s_{\max}$}{
                $s_{\max} \leftarrow s$;
                $O_{\max} \leftarrow O$\;
            }
        }
    }
    $S_0 \leftarrow S_0 \cup O_{\max}$\;
}
\end{algorithm2e}

We apply our heuristic to the same 100 instances, using solutions to $LP_0$ and $LP_1$. This time, we measure the relative gap between the resulting upper bound and either the previous lower bound or the one obtained when attempting to solve \eqref{eq:MILP_Inter2} in Gurobi, whichever is higher. When $\kappa = 0$, the median gap (resp.\ IQR) is 28\% (resp.\ 21\% -- 42\%). When $\kappa=1$, the gap dramatically improves to 12\% (resp.\ 8\% -- 20\%). Recall that the heuristic provides a feasible solution, so the choice of $\kappa=1$ already enables high-quality approximations to the optimal seed set. Increasing $\kappa$ allows us to improve upon those even more. In the case where we take $\kappa=\omega$ and (\ref{eq:LP.k}) returns the optimal solution $Y_{s^*_i}=s^*_i$, the heuristic always seems (empirically) to return the optimal solution.

\section{Managerial Insights and a Simple Heuristic for Solving  \texorpdfstring{$MIN$-$SCTM$}{MIN-SCTM}} \label{sec:heuristics}

We now show how our optimization framework sheds light on two critical questions:  \emph{How does the size of the seed set depend on the network structure?}  and \emph{What are the different roles that seed set firms play in the diffusion process?} Answering these questions can help traceability initiative leaders estimate the effort required to disseminate their technology and how to engage with different firms in the supply chain network \citep{WEF2021, Sandoval2022}. 

Section \ref{sec:insight1} addresses the first question and finds that supply chain networks with a higher degree of \emph{Jaccard clustering}\mdash a measure of how much overlap there is between firms' neighborhoods\mdash tend to have larger seed set sizes. Section \ref{sec:insight2} addresses the second question and shows that networks with intermediate levels of \emph{modularity}\mdash a measure of whether there are communities in the network and how tightly knit these communities are \citep{newman2006modularity}\mdash give rise to \emph{helper} nodes in the seed set. These nodes belong to supply chains that become traceable in later stages of the diffusion process and ``help" diffusion move between different network parts. Interestingly,  modularity and clustering frequently appear in supply chain network analyses \citep{perera2017network}.

While the LTM literature has studied the influence of modularity and clustering on network diffusion, there is no consensus on the direction of that influence. For instance, \cite{Acemoglu2011} suggest that network diffusion may be more widespread on networks with a smaller degree of clustering, in contrast to the conclusions given in \cite{centola2007cascade} and \cite{centola2010spread}. Similarly, \cite{shakarian2012large} find that highly modular graphs tend to have a lower diffusion rate, contradicting the conclusions of \cite{nematzadeh2014optimal}. Our numerical experiments will resolve these conflicts for traceability technology diffusion in supply chain networks.

\subsection{The Relationship Between Seed Set Size and Network Structure} \label{sec:insight1}

We examine which \emph{structural characteristics} of supply chain networks influence seed set \emph{size} and ultimately show that the higher the \emph{Jaccard clustering} of a supply chain network (which is defined below), the larger the seed set. Identifying the drivers of the seed set size has important practical implications for managers. First, it can help them evaluate the effort required to disseminate a technology in a particular network. For example, traceability initiative leaders (such as Walmart) interested in increasing the traceability of a set of supply chain networks might, \emph{ceteris paribus}, prefer to focus on product categories or industries whose supply chain networks display low Jaccard clustering. Second, firms interested in increasing their influence over their supply chain network should not only consider costs and risks when expanding their supply chains. They should also examine how their expansion strategy influences the network's degree of clustering.

Jaccard clustering was introduced in \cite{latapy2008basic} and is frequently used in hypergraph-based models \cite[e.g.,][]{klamt2009hypergraphs}. At a high level, this metric reflects whether the nodes in $G$ tend to belong to the same sets of hyperedges or not. Formally, let $E_i$ be the set of hyperedges that a node $i$ is a part of, i.e., $E_i = \{ j \in N_{SC} : i \in e_j \}$ and denote by $|E_i|$ its cardinality. Then, let the \emph{neighborhood similarity} between two nodes $i$ and $k$ of the hypergraph be the relative overlap between $E_i$ and $E_k$. Namely,
$ns(i,k)=\frac{|E_i \cap E_k|}{|E_i \cup  E_k|}.$
This measure of set overlap, also known as the \emph{Jaccard index} \citep{jaccard1912distribution}, is commonly used in network theory, machine learning, and biology. If $E_i$ and $E_k$ overlap perfectly, then $ns(i,k) = 1$. Conversely, if $E_i$ and $E_k$ are very different\mdash for example if $i$ belongs to many hyperedges that do not contain $k$\mdash  then $ns(i,k)$ will be low. We denote the set of nodes that share at least one hyperedge with $i$ as the \emph{neighborhood} of $i$, defined as $\mathcal{N}_i = \{ i'\in N_F : 
| E_i \cap E_k | \geq  1 \}$. The \emph{Jaccard clustering coefficient} of $i$ is the average neighborhood similarity between $i$ and its neighbors. Formally,
\begin{equation} \label{eq:cluster.i.def}
J_i=\frac{\sum_{k \in \mathcal{N}_i} ns(i,k)}{|\mathcal{N}_i|}.
\end{equation}
If all nodes in the neighborhood of $i$ only belong to a single hyperedge, i.e., all firms belong to a single supply chain (to which $i$ also belongs), then $J_i = 1$. Conversely, in a star graph with each node sharing one hyperedge with the central node $i$, $J_i$ converges to zero as the number of nodes grows. To obtain a clustering measure for a hypergraph $G$, we compute the average clustering coefficient of all its nodes: $J = \frac{1}{n} \sum_{i\in N_F} J_i$, which we denote as the \emph{Jaccard clustering} of the network.

\begin{figure}[t]
    \centering
    \subfloat[Higher Jaccard clustering]{
        \includegraphics{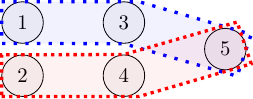}
        \label{fig:clustering1} 
    }   
    \hspace{1cm}
    \subfloat[Lower Jaccard clustering]{
        \includegraphics{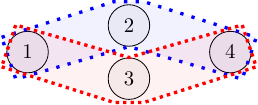}
        \label{fig:clustering2} 
    }  
    \caption{Illustration of Jaccard clustering.}\label{fig:scs}
\end{figure}

We illustrate this definition through the example in Figure \ref{fig:scs}. In Figure \ref{fig:clustering1}, we have  $ns(1,5) = ns(2,5) = ns(3,5) = ns(4,5) = 1/2$ and  $ns(1,3) = ns(2,4) = 1$. Thus, $J_5 = \frac{1}{2}$ and $J_i = \frac{1/2 + 1}{2} = 0.75$ for $i \neq 5$, so $J=0.7$. In Figure \ref{fig:clustering2}, we have $ns(1,2) = ns(1,3) = ns(2,4) = ns(3,4) = 1/2$ while $ns(1,4) = 1$. Thus, $J_1 = J_4 =  \frac{2}{3}$ and $J_2 = J_3 =  \frac{1}{2}$, giving $J=0.5833$.
The Jaccard clustering of the supply chain network in Figure \ref{fig:clustering1} is higher than in Figure \ref{fig:clustering2} since there are fewer overlaps between the two supply chains and more firms that share the same neighborhoods. 

The higher Jaccard clustering, the larger the seed set. Take Figure \ref{fig:scs}, assume that $w_i = c_i = r_{ij}=1$ and $\theta_j = 3$ for all $i \in N_F, j \in N_{SC}$, and recall that Figure \ref{fig:clustering1} has higher clustering than Figure \ref{fig:clustering2}. Correspondingly, Figure \ref{fig:clustering1} requires a larger percentage of firms to be seeded (60\% vs. 50\%). This is because the two supply chains in Figure \ref{fig:clustering1} have little overlap and, due to the supply chain effect, activation from one does not entirely transfer to the other. More generally, the higher the Jaccard clustering, the more likely there will be clusters of supply chains with limited overlap, and the more difficult it will be for the diffusion process to enter regions of the network with non-traceable supply chains without adding additional active firms from within that region. We use problem \eqref{eq:MILP_Inter2} to further formalize this reasoning through an optimization lens in Appendix \ref{subsec:cluster.IP}.

We present several numerical experiments confirming that Jaccard clustering accurately predicts seed set size and has more predictive power than alternative clustering metrics and other structural network characteristics. Our experiments use three sets of supply chain networks, described in detail in Appendix \ref{appendix:data}, with more than 1,600 instances across them. The first set is derived from the real-world supply chain networks in \cite{Willems2008}. The second set consists of randomly generated networks with structural characteristics similar to supply chains found in practice. The third set contains random networks used in Section \ref{sec:insight2}. Our results are consistent across the three datasets.

\begin{figure}[t]
    \centering
    \subfloat[\citet{Willems2008}]{ \label{fig:correlations1}
    \begin{adjustbox}{clip,trim=0cm 0.3cm 0cm 0cm,max width=\textwidth}
        \includegraphics{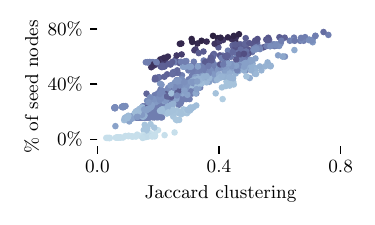}
    \end{adjustbox}
    }   
    \hspace{-0.5cm}
    \subfloat[Random]{ \label{fig:correlations2}
        \begin{adjustbox}{clip,trim=0.7cm 0.3cm 0cm 0cm,max width=\textwidth}
            \includegraphics{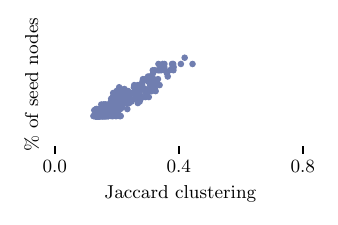}
        \end{adjustbox}
    }  
    \hspace{-0.5cm}
    \subfloat[Modular]{ \label{fig:correlations3}
        \begin{adjustbox}{clip,trim=0.7cm 0.3cm 0cm 0cm,max width=\textwidth}
            \includegraphics{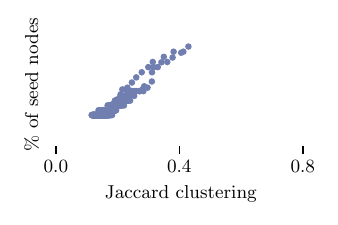}
        \end{adjustbox}
    }  
    \caption{Relationship between Jaccard clustering and the seed set size for different datasets. Colors vary with the average length of supply chains in the networks, from 2 (light) to 8 (dark).}\label{fig:correlations}
\end{figure}

We first analyze the relationship between Jaccard clustering and the percentage of total firms in the seed set, depicted in Figure \ref{fig:correlations}. As one would expect from the supply chain effect, the longer the length of the supply chains in the network, the more seed nodes are required. In addition, Jaccard clustering is highly correlated with the percentage of firms in the seed set, especially when considering supply chain networks with the same average supply chain length. This is clearly visible in Figures \ref{fig:correlations2} and \ref{fig:correlations3}, where all supply chains have length five.

The relationship between Jaccard clustering and the percentage of firms in the seed set is further confirmed by the correlation coefficient, which we compare to that of alternative clustering metrics\footnote{In hypergraphs, there is no consensus on the definition of clustering. Jaccard clustering is perhaps not an immediate choice as it requires two rounds of averaging to obtain a measure over graphs and does not reduce to the most common graph clustering metric when applied to hypergraphs with hyperedges of size 2 (i.e., graphs).} (see Appendix \ref{appendix:measures} for definitions).
In Table \ref{tab:clustering}, we observe that Jaccard clustering is an excellent predictor of seed set size, its performance surpassing other clustering metrics across all datasets.

\begin{table}[t]
\centering
\caption{Correlation coefficients between the percentage of seed nodes in the optimal solution and different clustering metrics. The highest for each dataset is highlighted in bold.}
\label{tab:clustering}
\begin{tabular}{lccc}
\hline 
 \up\down
 & \citet{Willems2008} & Random & Modular \\
\hline
\up
Jaccard clustering  & \textbf{0.7743} & \textbf{0.9114} & \textbf{0.9446} \\
Projection clustering & 0.7196 & 0.4819 & 0.4359 \\
Projection clustering (weighted) & 0.7081 & 0.3729 & 0.4482 \\
Hourglass clustering & 0.3265 & 0.5110 & 0.0950 \\
Repetition of partners & 0.2327 & -0.8301 & -0.8052 \\
\hline
\end{tabular}
\end{table}

Next, we create predictive models of the seed set percentage as a function of a large set of key network measures, including various clustering metrics. Appendix \ref{appendix:measures} states these metrics, and our analysis is in Appendix \ref{appendix:prediction}. When we use random forest models to predict the seed set percentage, we achieve root mean square errors of $0.016$--$0.024$ on unseen test data. Jaccard clustering is the first- or second-most predictive variable on each dataset. We then conduct an experiment using a penalized regression model and vary the penalty parameter. As we increase the penalty term, i.e., as we force the regression model to rely on fewer explanatory variables, the importance of Jaccard clustering increases and consistently reaches the first position among different metrics.

Our results indicate that, for supply chain networks, the effect discussed in \citet{Acemoglu2011}, that higher clustering\footnote{We note that \citet{Acemoglu2011} uses a form of hourglass clustering.} leads to clusters that the diffusion process has trouble ``entering" unless there are active nodes, dominates the effect observed in \citet{centola2007cascade} and \citet{centola2010spread}, that nodes in highly-clustered networks share more neighbors and more opportunities for diffusion. This is a consequence of the supply chain effect.

\subsection{Two Kinds of Seed Set Firms: Starters and Helpers} \label{sec:insight2}

We now investigate the different roles played by seed set nodes in the diffusion process. To illustrate these different roles, consider the example in Figure \ref{fig:adoption_example}. The seed set consists of Firms 1, 2, 3, 4, and 7. While these firms are active at the beginning of the diffusion process, the hyperedges (i.e., supply chains) they belong to become active at different times. Hyperedge $e_{black}$ containing Firms 2, 4, and 7 becomes active in the first period of the diffusion process (as we have $\theta_{black}=4$) while hyperedges $e_{red}$ and $e_{blue}$ containing Firms 1 and 3 become active in the second and third periods, respectively. This example highlights a pattern that we observe across instances of $MIN$-$SCTM$. Namely, firms in the seed set can be of two types: firms that belong to hyperedges that become active at the \emph{start} of the diffusion process, or firms that belong to hyperedges that only become active at \emph{later} stages of the diffusion process. We refer to the former category as \emph{starter firms} and to the latter as \emph{helper firms} because these firms ``help'' diffusion after the first period. Firms 2, 4, and 7 in Figure \ref{fig:adoption_example} are starter firms, while 1 and 3 are helper firms. 

\begin{figure}[t]
    \centering
    \subfloat[High modularity]{
        \includegraphics{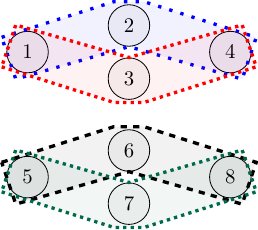}
        \label{fig:modularity1} 
    }    
    \hspace{0.5cm}
    \subfloat[Intermediate modularity]{
        \includegraphics{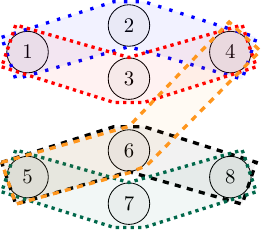}
        \label{fig:modularity2} 
    } 
    \hspace{0.5cm}
    \subfloat[Low modularity]{
        \includegraphics{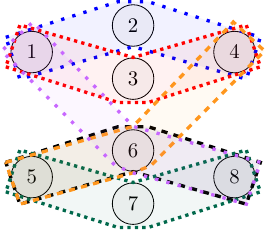}
        \label{fig:modularity3} 
    }  
    \caption{Illustration of modularity. The parameters are $w_i = c_i = r_{ij}=1$ and $\theta_j = 3$ for all $i \in N_F, j \in N_{SC}$.}\label{fig:mod}
\end{figure}

From a managerial perspective, starter and helper firms play very different roles in the diffusion process. Starter firms immediately lead to the traceability of a set of supply chains. Making these supply chains (or, equivalently, the products they produce) traceable ``jumpstarts" diffusion in the network. Conversely, helper firms are targeted without the explicit intent to make their supply chains traceable early in the diffusion process but to help keep the diffusion process ``moving" along the network. As a consequence, helper firms are strategically placed, typically at the juncture between different network parts, as they ``transfer" the activation process from one part to another. Interestingly, the strategic importance of ``helper" firms has been picked up on in the practitioner literature (for example, \citealt{WEF2021} call them ``alliance brokers"). To illustrate the different roles played by the seed set firms, consider the example in Figure \ref{fig:modularity2}. Seeding Firms 1 and 2 at the top of the network will not be enough to ensure that the four firms at the bottom become active. One helper firm, such as 6, must be added to the seed set for diffusion to complete. 

Managers may follow a strategy where starter firms are seeded at the beginning of the diffusion process while helper firms are seeded as the diffusion process ``approaches" their supply chains. For the network in Figure \ref{fig:adoption_example}, this would mean seeding Firms 2, 4, and 7  early on while seeding Firms 1 and 3 after $e_{black}$ and $e_{green}$ are active. Our IP formulation for $MIN$-$SCTM$  outputs the order in which hyperedges become active, which can be used to identify starter and helper firms.

We now examine how the supply chain network structure influences the seed set's proportion of starter and helper firms. 
Our main observation is that the proportion of starter firms in the seed set is \emph{a V-shaped function} of modularity.
Given a set of groups of nodes within a network, modularity measures how connected different groups are. Thus, a very modular graph has ``tightly connected" groups of nodes that are  ``loosely connected" with other groups. 
Introducing a modularity metric for hypergraphs is the first step toward establishing this V-shaped relationship. Appendix \ref{appendix:measures} discusses how we specialize the commonly-used definition in \citet{newman2006modularity} to our setting.

Consider the networks in Figure \ref{fig:mod} to develop intuition around this relationship. When a supply chain network has high modularity, there are near-disconnected groups in the network, and there is little possibility for helper firms to transfer activation across groups. In Figure \ref{fig:modularity1}, modularity is the highest among the three networks (0.5), and we require four starter firms out of four seed firms (e.g., 1, 2 and 5, 6). When modularity is intermediate, groups of nodes are somewhat connected, and helper firms emerge to transfer activation across these groups. In Figure \ref{fig:modularity2}, modularity is intermediate (0.431), and we require two starter firms and one helper firm out of three seed firms (e.g., 5, 6, and 4). Finally, when modularity is low, there are not many loosely-connected groups. As a result, there is less of a need for helper firms to transfer activation across groups. In Figure \ref{fig:modularity3}, modularity is the lowest (0.401), and both seed set firms are starter firms (e.g., 5, 6). In consequence, the percentage of starters in the seed set is a V-shaped function of modularity.

\begin{figure}[t]
    \centering
    \subfloat[$\overline{m} = 20$]{\label{fig:mod_u_40}
        \begin{adjustbox}{clip,trim=0cm 0.3cm 0cm 0cm,max width=\textwidth}
            \includegraphics{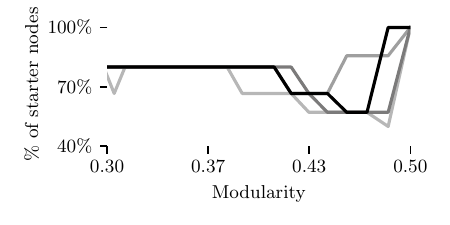}
        \end{adjustbox}
    }
    \hspace{0.5cm}
    \subfloat[$\overline{m} = 30$]{\label{fig:mod_u_60}
        \begin{adjustbox}{clip,trim=0cm 0.3cm 0cm 0cm,max width=\textwidth}
            \includegraphics{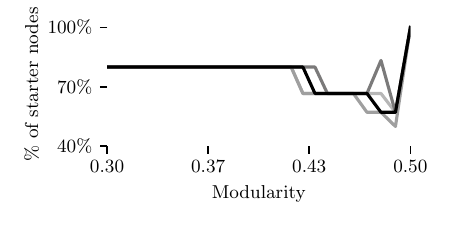}
        \end{adjustbox}
    }
    \caption{Starter firm percentage for varying numbers of connecting hyperedges on otherwise identical two-part networks, as described in Appendix \ref{appendix:data}.
    Each line corresponds to different starting network.}
    \label{fig:mod_u}
\end{figure}

To depict the emergence of helper firms, we randomly generate supply chain networks with varying modularity (see Appendix \ref{appendix:data}). Simply put, we first generate a network with $\bar{m}$ hyperedges, which we duplicate to obtain two identical but disconnected groups of nodes. We then sequentially add $\bar{m}$ hyperedges across the groups, calculating the modularity, starter firms, and helper firms of the resulting network each time we add a  hyperedge. As the number of cross-group hyperedges increases, the modularity of the overall supply chain network decreases. Figure  \ref{fig:mod_u} shows the percentage of starter firms in the seed set as a function of the network's modularity. Each line of the plot corresponds to an instance of this cross-group hyperedge addition process with a randomly generated initial network and with the number of cross-group hyperedges varying from 0 to $\bar{m}$. 

On the left-hand side of the plots, there are $\bar{m}$ cross-group hyperedges, and modularity is low, so activation flows unhindered through the network, and the seed set consists primarily of starter firms. On the right-hand side of the plots, there are zero cross-group hyperedges, and modularity is high, so there are few or no helper firms since there is limited diffusion across the two network parts. In the middle of the plots, there are some cross-group hyperedges, and modularity is intermediate. Helper firms emerge to ``help" transfer the diffusion process across the two network parts. Thus, the percentage of starter firms in the seed set decreases. Note that other network metrics heavily correlated with modularity, such as graph density, would describe the same effect.

Our insights indicate that traceability initiative leaders should be mindful of the modularity of their supply chain networks when searching for helper firms to target. If the network's modularity is intermediate, helper firms at the juncture of loosely connected network parts are likely needed to ensure broad technology diffusion. These firms help
circumvent the supply chain effect and can be targeted later in the diffusion process.

\subsection{A Simple Heuristic} \label{subsec:heuristic}

We use the insights derived in Sections \ref{sec:insight1} and \ref{sec:insight2} to propose a simple seeding heuristic for managers. This heuristic first computes a \emph{score} for each combination of firms that, when active, lead to a supply chain becoming active: the score includes the model parameters and Jaccard clustering.
Once these scores have been computed, the heuristic seeds the set of firms with the highest score and lets diffusion propagate until it stops. It then repeats this process on the remaining non-active parts of the graph. By proceeding sequentially, the heuristic mimics the strategy mentioned in Section \ref{sec:insight2} of first seeding starter firms and then seeding helper firms as diffusion propagates.

The heuristic follows Algorithm \ref{alg:ub}, but instead of scores derived from an LP, it scores a node's Jaccard clustering. Specifically, in the heuristic, Line 8 of the algorithm becomes
$$
s \leftarrow \frac{1}{|O|}  \frac{ \left(\sum_{i \in O} J_i \right) \left( \sum_{i \in O} c_i \right) }{\left(\sum_{i \in O} w_i \right) \left( \sum_{i \in O} r_{ji} \right)}.
$$
Note that computing this score is simple and does not require any optimization software. 

When comparing the heuristic to a lower bound, as in Section \ref{sec:lb}, the median gap (resp.\ IQR) is 19\% (resp.\ 13\% -- 28\%). In comparison, when paths are chosen randomly (but the choice of nodes within paths is still optimized based on $w_i$), the gap is 44\% (resp.\ 25\% -- 76\%).\footnote{The median gap is 94\% (resp.\ 74\% -- 137\%) when nodes within paths are also chosen at random. When supply chain effects are ignored, and random selection is by node instead of by path, the gap is 135\% (resp.\ 95\% -- 245\%).}
The heuristic even outperforms the upper bound from Section \ref{sec:lb} with $\kappa=0$.
However, it performs worse than the upper bound when $\kappa=1$ and, naturally, even worse as $\kappa$ grows further. Hence, the heuristic is a starting point for managers who desire a quick estimate of the cost of seeding or the companies to seed. When more computation time and an optimization solver are available, the principled heuristics providing upper and lower bounds from Section \ref{sec:lb} are preferable.

\section{Conclusion} \label{sec:conclusion}

Modern traceability technologies can alleviate supply chain risks, improve sustainability, reduce transaction costs, and enhance demand. However, the benefits of such technologies are only unlocked when subsets of firms in a supply chain adopt the technologies, a phenomenon we refer to as the supply chain effect. This effect has profound implications for companies leading traceability initiatives that often struggle to design a dissemination strategy for their technology \citep{WEF2021}. Successful traceability technology dissemination strategies must address a few key questions:\emph{(i) What is the lowest-cost seed set that ensures the whole network eventually adopts the technology?} \emph{(ii) How does the size of the seed set depend on the network structure?}  and \emph{(iii) What are the different roles that seed set firms play in the diffusion process?} While the technology diffusion literature offers several network models and optimization frameworks to address these questions, they cannot be readily applied to supply chain networks due to the supply chain effect.

To address the questions above, we contribute to the technology diffusion literature by introducing the Supply Chain Traceability Model as a new framework incorporating the supply chain effect. We then define $MIN$-$SCTM$, the problem of finding the minimum cost seed set of nodes that guarantees diffusion throughout the network, which can be viewed as the optimization formulation of question \emph{(i)}. We prove that $MIN$-$SCTM$ is not just NP-hard; it is inapproximable in polynomial time, even when all supply chains in the network only contain three firms and when the cost-benefit analysis is trivial. This result indicates that the supply chain effect and the network's intricate structure are the two drivers of $MIN$-$SCTM$'s complexity. Thus, any effective procedure for solving $MIN$-$SCTM$ must take advantage of particular network structures.

Therefore, we design an LP-based FPT algorithm for $MIN$-$SCTM$ with parameter treewidth of the supply chain network. The use of treewidth is practice-driven: the treewidth of publicly available supply chain networks is often an order of magnitude smaller than the network size. The approach we use to design our FPT algorithm is based on recent integer programming techniques that are new to the technology diffusion literature. Specifically, we show how to bound the treewidth of the intersection graph of an integer programming formulation of $MIN$-$SCTM$ by the treewidth of the network where technology diffusion occurs.  
The resulting FPT algorithm is an explicit LP formulation of this integer program that can be solved using existing optimization solvers. Our procedure also outputs a hierarchy of approximations of $MIN$-$SCTM$ with an explicit tradeoff between the accuracy of the solution and the time taken to compute it. In short, our optimization framework answers question \emph{(i)} by introducing algorithms that solve large instances of $MIN$-$SCTM$ to near-optimality within reasonable computational time and that are easy to implement.

We further employ our optimization framework to address questions \emph{(ii)} and \emph{(iii)} and obtain several new managerial insights. For the relationship between the network structure and the seed set size, we observe that a supply chain network with high Jaccard clustering tends to have a large optimal seed set. This result contrasts with existing and conflicting results in the network diffusion literature. While higher clustering might facilitate diffusion due to neighborhood overlaps (leading to smaller seed sets), it also makes it more difficult for the diffusion process to enter tightly-knit ``clusters" of firms (leading to larger seed sets). In supply chain networks, the latter phenomenon dominates due to the supply chain effect: If firms tend to be in the same supply chains, more firms are needed in the seed set.
As for the different roles played by seed set firms, we observe that networks with an intermediate degree of modularity require ``helper" firms, i.e., firms that help transfer diffusion between different parts of the network. These firms are part of supply chains that only activate in later stages of the diffusion process and can thus be targeted later. Collectively, our insights can help managers shape their technology diffusion strategy, for example, by indicating which product categories (and corresponding supply chain networks) tend to require smaller seed sets and which networks are more likely to need helper firms to promote diffusion. 

A promising future research direction is to extend $SCTM$ to address the limitations discussed in Section \ref{subsec:case_studies}. In particular, assuming that firms engage in strategic games or adjust sourcing decisions based on their technology adoption could lead to interesting new results. Another exciting research direction is to explore the connections between integer programming techniques and network diffusion problems. For instance, examining LP formulations for seed set selection that depend on other graph parameters (beyond treewidth) might produce new approaches and insights. 

More broadly, the optimization-based tool set and analysis we develop might be helpful for other problems requiring the coordination of multiple firms in supply chains, such as adopting sustainable supply chain practices. Consider, for example, circular economy initiatives. A supply chain can only really be ``circular" if all companies in the supply chain adopt a consistent set of practices and technologies. We speculate that the tools and approaches introduced in this paper can assist with designing new sustainable supply chain strategies.

\ACKNOWLEDGMENT{We would like to thank the Department Editor, the Associate Editor, and two anonymous reviewers for their constructive comments and suggestions throughout the revision process. We would also like to express our gratitude to IBM Research and Ashish Jagmohan, for stimulating the initial research on this project, and to Manpreet Hoora, Dan Iancu, Mihalis Markakis, Gonzalo Mu\~{n}oz, Karthik Ramachandran, and Beril Toktay who provided invaluable reflections at different stages of the revision process. Finally, we would like to acknowledge the research cyberinfratsructure resources and services provided by the Partnership for an Advanced Computing Envrionment (PACE) at Georgia Institute of Technology.}

\SingleSpacedXI
 
\putbib[library]
 
\end{bibunit}

\newpage
\OneAndAHalfSpacedXI

\begin{APPENDICES}
\begin{bibunit}[informs2014]
\section{Results Linked to the Model Definition} \label{appendix:model}

This section contains the proofs of all results contained in Section \ref{sec:SCTM}.

\proof{Proof of Proposition \ref{prop:sctm.eq.ltm}}
As $S_0' \subseteq N_F$ and the graph is bipartite, $S'_{2t+1} \backslash S'_{2t} \subseteq N_{SC}$ and $S'_{2t+2} \backslash S'_{2t+1} \subseteq N_F$. In other words, SC-nodes activate in odd time periods and firm-nodes activate in even time periods. In particular, at time $2t+1$, the nodes added are
$$\left\{j \in N_{SC} \backslash S_{2t} ~:~ |\{(i,j) \in E'~:~i \text{ is active} \}|  \geq \theta_j - 1 \right\}= \bigcup_{i \in S_{2t} \cap N_F} \mathcal{B}_i(S'_{2t} \cap N_F).$$
In light of this, at time $2t+2$, the nodes added are
$$\left\{i \in N_{F} \backslash S_{2t+1} ~:~ \sum_{(j,i) \in E' ~:~ \text{ $j$ is active}} r_{ji} \geq c_i \right\}= \left\{i \in N_{F} \backslash S_{2t+1} ~:~ \sum_{j \in \mathcal{B}_i(S'_{2t} \cap N_F)} r_{ji} \geq c_i \right\}.$$
From this and the fact that no firm nodes are added at time $2t+1$, it follows that
$$S'_{2t+2} \cap N_F = (S'_{2t} \cap N_F) \cup \left\{i \in N_{F} \backslash S_{2t} ~:~ \sum_{j \in \mathcal{B}_i(S'_{2t} \cap N_F)} r_{ji} \geq c_i \right\}.$$
As $S_0=S'_0$, we get the state equations of the SCTM activation process, so $S'_{2t} \cap N_F=S_{t}$, $\forall t$.
\Halmos
\endproof

\proof{Proof of Corollary \ref{cor:sctm.ltm}}
This is an immediate consequence of Proposition \ref{prop:sctm.eq.ltm} and the fact that if $S_{\infty}=N_F$, then it must be the case that $S'_{\infty}=N'$, by virtue of $| \{i \in N_F\ :\ i \in e_j \}| = k_j \geq \theta_j$ for all $j \in N_{SC}$.
\Halmos
\endproof
\section{Computational Complexity Results} \label{appendix:comp}

\subsection{Definition of Building Blocks Used in the Proofs} \label{appx:proof_prelims}

The proofs in this section require similar hypergraph structures with edges of size $k_j=3$ and thresholds $\theta_j = 3$. These ``building blocks'' $B_U$ (resp.\ $C_{U,V}$) consist of 3 (resp.\ 5) nodes and are defined in Figure \ref{fig:NP_hard_proof_construction} in gray (resp.\ blue). 
$C_{U,V} \xrightarrow{L,R} B_V$ denotes that links $C_{U,V} \xrightarrow{L} B_V$ and $C_{U,V} \xrightarrow{R} B_V$ exist.

\begin{figure}[htp!]
    \subfloat[$B_U \leftrightarrow C_{U,V}$ with 7 hyperedges ($B_U$ active $\Leftrightarrow$ $C_{U,V}$ active). \label{fig:NP_hard_proof_construction_outgoing}]{ 
        \includegraphics{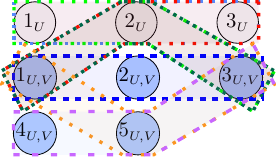}
    }
    \hspace{0.8cm}
        \subfloat[$C_{U,V} \xrightarrow{L} B_V$ with 5 hyperedges ($C_{U,V}$ active $\Rightarrow$ $B_V$ requires $1_V$ or $3_V$ to activate). \label{fig:NP_hard_proof_construction_incoming_L}]{ 
        \includegraphics{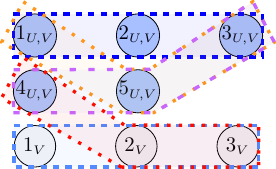}
    }
    \hspace{0.8cm}
        \subfloat[$C_{U,V} \xrightarrow{R} B_V$ with 5 hyperedges ($C_{U,V}$ active $\Rightarrow$ $B_V$ requires $1_V$ or $2_V$ to activate). \label{fig:NP_hard_proof_construction_incoming_R}]{ 
        \includegraphics{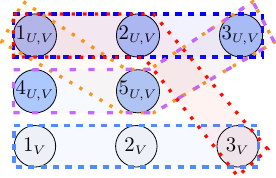}
    }
    \caption{Construction of building blocks. \label{fig:NP_hard_proof_construction}}
\end{figure}

\subsection{NP-Hardness of \texorpdfstring{$MIN$-$SCTM$}{MIN-SCTM(eps)} (Theorem \ref{thm:comp})} \label{appendix:hardness.MIN.SCTM.eps}

\begin{lemma} \label{lem:np.hard.k3.theta3}
$DEC$-$SCTM$ is NP-hard when $k=3$ and $\theta=3$.
\end{lemma}

\proof{Proof.} Let $DEC$-$LTM$ be the following decision problem:\\
\textsc{Input:} A directed graph $G'=(N',E')$ with weights $w'_{ij}\ \forall (i,j) \in E'$ and thresholds $c'_i\ \forall i \in N'$, and $h' \in \mathbb{N}$.\\
\textsc{Question}: Is there a seed set of size less than or equal $h'$ leading to full activation of $G$ in the LTM sense?\footnote{Recall that a node activates in the LTM sense if the sum of its incoming edge weights from active nodes exceeds its threshold; see Section \ref{sec:aux_graph}.}\\
$DEC\text{-}LTM$ is NP-hard, even if $w'_{ij} = 1\ \forall (i,j) \in E'$ and $c'_i=|\mathcal{N}_i|\ \forall i \in N'$ \citep[see][Proof of Theorem 2.7]{Kempe2003}, which we assume. Wlog, we also assume that $|\mathcal{N}_i| > 0$ for $i \in N'$ (if a node has no incoming edges, its benefit is equal to 0, as is its threshold. It will thus always activate by itself and would never be part of a minimal seed set).
We construct a reduction from $DEC$-$LTM$ to $DEC$-$SCTM$.

\paragraph*{Construction of the reduction.} We use the building blocks from Appendix \ref{appx:proof_prelims} to build $G$ from $G'$, indexing blocks $B_U$ and $C_{U,V}$ using sets of nodes in $N'$, i.e., $U,V \subseteq N'$. For any node $i$ in $N'$, add a block $B_{i}$ to $G$. Recursively split its set of incoming neighbors $\mathcal{N}_i$ into two sets of equal size (if $|\mathcal{N}_i|$ is even), or into two sets of size differing by one (if $|\mathcal{N}_i|$ is odd). In other words, $\mathcal{N}_i = \mathcal{N}_i^0 \cup \mathcal{N}_i^1 = \{\mathcal{N}_i^{0,0} \cup \mathcal{N}_i^{0,1} \} \cup \{ \mathcal{N}_i^{1,0} \cup \mathcal{N}_i^{1,1}\} = \ldots$, stopping when the sets in the decomposition contain one element. Add corresponding blocks $B_{\mathcal{N}_i^0},$ $B_{\mathcal{N}_i^1},$ $B_{\mathcal{N}_i^{0,0}},$ $B_{\mathcal{N}_i^{0,1}},$ $B_{\mathcal{N}_i^{1,0}},$ $B_{\mathcal{N}_i^{1,1}}, \ldots$ to $G$, except when $|\mathcal{N}_i^x| = 1$. Further, add blocks $C_{\mathcal{N}_i^0,\{i\}},$ $C_{\mathcal{N}_i^1,\{i\}},$ $C_{\mathcal{N}_i^{0,0},\mathcal{N}_i^0},$ $C_{\mathcal{N}_i^{0,1},\mathcal{N}_i^0},$ $C_{\mathcal{N}_i^{1,0},\mathcal{N}_i^1},$ $C_{\mathcal{N}_i^{1,1},\mathcal{N}_i^1}\ldots$ to $G$ and link the blocks as follows: $C_{\mathcal{N}_i^0,\{i\}} \xrightarrow{L} B_{\{i\}}$, $C_{\mathcal{N}_i^1,\{i\}} \xrightarrow{R} B_{\{i\}}$, $B_{\mathcal{N}_i^0} \leftrightarrow C_{\mathcal{N}_i^0,\{i\}}$, $B_{\mathcal{N}_i^1} \leftrightarrow C_{\mathcal{N}_i^1,\{i\}}, \ldots$. This is illustrated in Figure \ref{fig:NP_hard_proof_construction_incoming_total}. If a block already exists (e.g., if incoming neighbors overlap), use the existing one. Finally, add a block $B_0$ and blocks $C_{\{i\},0}\ \forall i \in N$ with $B_{\{i\}} \leftrightarrow C_{\{i\},0}$ and $C_{\{i\},0} \xrightarrow{L,R} B_0$. As constructed, $G$ is a hypergraph with edges of size $k=3$. Furthermore, we assume that $r_{ji}=1$, $c_i=1$, $\theta_j=3$, and $w_i=1$ for $i=1,\ldots,n$ and $j=1,\ldots,m$, and we let $h=h'+1$.

\begin{figure}[t]
    \centering
    \includegraphics{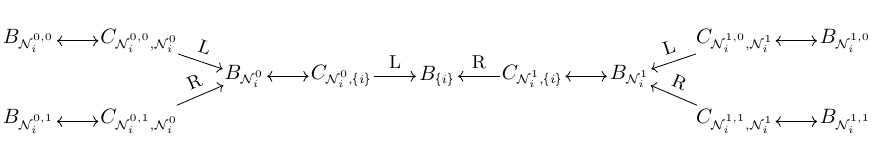}
    \caption{Linkage of building blocks for incoming edges.}  \label{fig:NP_hard_proof_construction_incoming_total}
\end{figure}

This construction is polynomial in $|N'|$ as the number of blocks is polynomial in $|N'|.$ Indeed, as $|\mathcal{N}_i| \leq |N'|-1\ \forall i \in N$, there are at most $\lceil \log_2(|N'|-1) \rceil$ recursive splits on $\mathcal{N}_i$, which implies that the number of sets $\mathcal{N}_{i}^0, \mathcal{N}_i^1, \ldots$ generated for node $i$ is at most equal to $2^1+2^2+\ldots+2^{\lceil \log_2(|N'|-1) \rceil} \leq 8(|N'|-1)$.

Moreover, the LTM activation process in $G'$ can be replicated via SCTM activation in $G$ by equating active nodes $i$ in $G'$ with active blocks $B_{\{i\}}$ in $G$ and assuming that no other nodes are initially active in $G$. We can then show that a node $i_a$ activates in $G'$ if and only if block $B_{\{i_a\}}$ activates in its entirety in $G$. To see this, note that each block $B_{\{i\}}$ only appears in two places in $G$: exactly once at the root of a tree such as the one given in Figure \ref{fig:NP_hard_proof_construction_incoming_total} and possibly many times as a leaf of trees of this type, associated to other blocks $B_{\{j\}}$.  When $B_{\{i\}}$ is at the root of a tree, it activates if all of the leaves of the tree (corresponding to blocks $B_{\{j\}}, j \in \mathcal{N}_i$) activate. If it is a leaf, it will not activate through this particular tree if the root of the tree activates, nor if any other leaves activate. From this, we prove the statement: if node $i_a$ is active in $G'$, then all incoming neighbors $j\in \mathcal{N}_{i_a}$ of $i_a$ are active. This implies that blocks $B_{\{j\}}, j \in \mathcal{N}_{i_a}$ are fully active, which leads from our previous discussion, to activation of $B_{\{i_{a}\}}$. Conversely, if $B_{\{i_a\}}$ activates with only blocks $B_{\{i\}}, i \in N'$ having initially been activated, then it must be the case that all blocks in the tree rooted at $B_{\{i_a\}}$ became active at some point, which can only happen if all of its leaves $B_{\{j\}}, j \in \mathcal{N}_{i_a}$ were fully activated. By equivalence, this means that in $G'$, nodes $j \in \mathcal{N}_{i_a}$ are active, and so $i_a$ would activate.

\paragraph*{\texorpdfstring{$DEC\text{-}LTM$}{} with $(G',h')$ answers YES if and only if \texorpdfstring{$DEC\text{-}SCTM$}{} with $(G,h)$ answers YES.}
``Only if'': Take a seed set $S_0'$ in the LTM of size $h'$ leading to full activation. For each $i \in S_0'$, add the corresponding node $5_{\{i\},0}$ to the seed set $S_0$ of the SCTM. Finally, take any one of the nodes in $S_0'$, say $j$, and add the node $1_{\{j\},0}$ to $S_0$, such that $|S_0| = h+1$. Clearly, the block $C_{\{j\},0}$ fully activates, then, block $B_0$, then each block $C_{\{i\},0}$ with $i \in S_0$. Next, $B_{\{i\}}$ activates for all $i \in S_0$. Each other block $B_{\{j\}}$ with $j \notin S_0$ only activates if all incoming nodes are fully activated. Hence, activation proceeds exactly as under the LTM and full activation eventually occurs under the SCTM.

``If'': Assume NO for $DEC\text{-}LTM$ and YES for $DEC\text{-}SCTM$. Let $S_0$ be an activating seed set for $G$ of size $h$.
If there is a node $j \in S_0 \cap  \left( C_{\mathcal{N}_i^x,\mathcal{N}_i^y} \cup B_{\mathcal{N}_i^x} \right)$ for some $\mathcal{N}_i^x$ that is not a singleton, then there is another activating seed set $\tilde{S}_0$, $|\tilde{S}_0| \leq |S_0|$, without $j$, but with $2_{\{i\}}$ or $3_{\{i\}}$.
Consider a block associated with $i \in N'$ that contains a node $j \in \tilde{S}_0$. Once it is active, the block $B_0$ fully activates, which in turn implies that any other such block can be activated with only the node $5_{\{i\},0}$. The first block needs to have two nodes for any activation to occur and the seed set to be minimal. Assume this block is associated with node $j \in N'$.  We can replace the seed nodes in this block with the nodes $4_{\{j\},0}$ and $5_{\{j\},0}$ wlog. We then know from before that a set $S_0' =\{i | 5_{\{i\},0} \in \tilde{S}_0\}$ (of size $\leq h -1 = h'$) leads to full activation in the LTM, giving a contradiction.
\Halmos
\endproof

\proof{Proof of Theorem \ref{thm:comp}.} 
Lemma \ref{lem:np.hard.k3.theta3} shows that $DEC$-$SCTM$ is hard when $(k=3,\theta=3)$. We first show that the cases $\theta=1$ and $\theta=2$ are in P. If $\theta=1$, then it is easy to see that all firms will adopt at the first time step so the seed set can be taken to be empty. If $\theta=2$, then one need only seed the node of lowest cost. This will lead to all supply chains it belongs to activating, which will in turn lead to all supply chains connected to those supply chains to activate. Connectivity of $G$ thus ensures full activation. \\
We now show that if $DEC$-$SCTM$ is hard for fixed $k$ and $\theta=k$ and $\theta \leq k,$ then $DEC$-$SCTM$ remains hard for $k+1$ and $\theta+1=k+1.$ To do this, let $G$ be the graph over which $DEC$-$SCTM$ is hard with fixed $k$ and fixed $h$. We let $\tilde{G}$ be $G$, except that we add the same additional node $\iota$, to each hyperedge $e_j \in G$. Note that the size of each hyperedge is now $k+1$. We further let the cost to seed of the additional node be $w_{\iota}=1$ for all  $i \in N$, let the edge threshold be $k+1$, and take $\tilde{h}=h+1$. One can check that $S_0 \cup \{ \iota\}$ is a seed set of size less than or equal to $\tilde{h}$ leading to full activation of $\tilde{G}$ if and only if $S_0$ is a seed set of size less than equal to $h$ leading to full activation of $G$, which implies the result. Similarly, we can show that if $DEC$-$SCTM$ is hard for fixed $k$ and $\theta$, then $DEC$-$SCTM$ remains hard for $k+1$ and $\theta$. This involves building the same graph as $\tilde{G}$ above, except that we do not add $\{\iota\}$ to the seed set and we let $\tilde{h}=h$. Thus, for any $(\theta,k)$ such that $\theta  \leq k$, one can proceed from the case $(k=3,\theta=3)$ to $(k,\theta)$ by applying successively the operations $(k,\theta=k) \mapsto (k+1,\theta+1=k+1)$ and $(k,\theta) \mapsto (k+1,\theta)$ as described above.
\Halmos
\endproof
\vspace{0.2cm}

\subsection{Hardness of Approximation of \texorpdfstring{$MIN$-$SCTM$}{MIN-SCTM} (Proposition \ref{prop:eps-inapprox})} \label{appendix:hardness.approx.MIN.SCTM.eps}

\proof{Proof of Proposition \ref{prop:eps-inapprox}.} We proceed as before using a reduction from the LTM on a graph $G'$ where all nodes have a threshold $c_i' \leq 2$ and all edges have weight $w_{ij}'=1$.
Below, we will show that for any $G'$ with $|N'|=n'$ nodes and LTM activation, we can create, in polynomial time, a hypergraph $G$ with $r_{ji}=1, c_i=1$, $k=3,\theta=3,$ activating via the $SCTM$, with the following properties: (i) the number of nodes is $n \leq (n')^\beta$ for a constant $1 \leq \beta < \infty$; (ii) if $OPT_{LTM}$ is the size of the minimum seed set for $G'$, and $S_0^*$ is a minimum size seed set of $G$ with $OPT = |S_0^*|$ , then $OPT  = OPT_{LTM} + 1$.

Given the construction, assume that for all $\alpha > 1$ there is an algorithm approximating $MIN$-$SCTM$  with result $OPT'$ such that $OPT' = OPT \cdot O(\alpha^{\log^{1-\xi} n})$ for some $\xi > 0$. Clearly, $OPT'$ is an upper-bound on $OPT_{LTM}$. Moreover, $OPT' = (OPT_{LTM} + 1) \cdot O(\alpha^{\log^{1-\xi} n}) < OPT_{LTM} \cdot O(\alpha^{\log^{1-\xi} n}) < OPT_{LTM} \cdot O(\alpha^{\beta^{1-\xi} \log^{1-\xi} n'}) = OPT_{LTM} \cdot O\left(\left(\alpha^{\beta^{1-\xi}}\right)^{\log^{1-\xi} n'}\right)$. As $\alpha > 1$, let $\alpha = 2^{\frac{1}{\beta^{1-\xi}}}$. We then have a direct contradiction to the result that there is no polynomial-time approximation algorithm with output $OPT'$ and $OPT' < OPT_{LTM} \cdot O(2^{\log^{1-\xi} n'})$ for any $\xi > 0$  \citep[cf.\ Corollary 4.1 in][]{chen2009approximability}. The case where $k \geq \theta \geq 3$ can be obtained in a similar fashion using the extension operations given in the proof of Theorem \ref{thm:comp}.

\paragraph*{A note on the complexity class.}

Assume instead that $OPT' = OPT \cdot O(2^{\log^{1-\xi} n})$. Then, $OPT' < OPT_{LTM} \cdot O\left(2^{\beta^{1-\xi} \log^{1-\xi} n'}\right)$. To establish a contradiction, we require that there is a constant $M > 0$ such that $M \cdot 2^{\beta^{1-\xi} \log^{1-\xi} n'} \leq 2^{\log^{1-\xi} n'}$. This is equivalent to $M \leq 2^{\log^{1-\xi} n' \left( 1 - \beta^{1-\xi} \right)}$. However, unless $\beta = 1$, $1 - \beta^{1-\xi} < 0$, and the right-hand side tends to zero as $n$ goes to infinity. But because $M>0$, this does not hold.
The same issue arises in the sequence of proofs leading up to Corollary 4.1 in \citet{chen2009approximability}, where a reduction is created to a graph with a higher number of nodes at each step (that is, $\beta > 1$). This is of no consequence for the complexity class that we apply here, however, as we can replace $2$ by some $\alpha \in (1,2)$ in each step.

\paragraph*{Construction of the reduction.}

Let $G'=(N',E')$, $|N'| = n'$ be a directed graph, activated via the LTM. Assume that all nodes have a threshold $c_i' \in \{1,2\}$ and all edges have weight $w'_{ij}=1$. We proceed with a similar construction of a hypergraph $G = (N,E)$ with $k=3$, as in the proof of Theorem \ref{thm:comp}. For each $i \in N'$, define block $B_{\{i\}}$. Either (i) node $i$ has a threshold of 1, then for each (directed) edge $(j,i) \in E'$, construct a block $C_{\{j\},\{i\}}$, as well as the linkages $B_{\{j\}} \leftrightarrow C_{\{j\},\{i\}} \xrightarrow{L,R} B_{\{i\}}$. Alternatively, (ii) it has a threshold of 2. Say the set of incoming neighbors is $\mathcal{N}_i$. Then, construct block $B_{\{j,j'\}}$ for any pair $(j,j') \in \mathcal{N}_i^2$. There are $\binom{|\mathcal{N}_i|}{2} \leq \binom{n-1}{2} \leq n^2$ such pairs. For each pair, construct blocks $C_{\{j\},\{j,j'\}}$, $C_{\{j'\},\{j,j'\}}$, and $C_{\{j,j'\},\{i\}}$, as well as the following linkages: $B_{\{j\}} \leftrightarrow C_{\{j\},\{j,j'\}} \xrightarrow{L} B_{\{j,j'\}}$, $B_{\{j'\}} \leftrightarrow C_{\{j\},\{j,j'\}} \xrightarrow{R} B_{\{j,j'\}}$, and $B_{\{j,j'\}} \leftrightarrow C_{\{j,j'\},\{i\}} \xrightarrow{L,R} B_{\{i\}}$. As before, we add a block $B_0$ and the corresponding blocks $C_{\{i\},0}$ for all $i \in N$ with connections $B_{\{i\}} \leftrightarrow C_{\{i\},0} \xrightarrow{L,R} B_0$. The construction is, again, polynomial in $n'$. Moreover, following the same steps as before, one can show that the solution to $MIN$-$LTM$ on graph $G'$ is $OPT_{LTM}$ if and only if the solution to $MIN$-$SCTM$ on graph $G$ is $OPT = OPT_{LTM} + 1$. Note that $n \leq (n')^\beta$ for a constant $1 \leq \beta < \infty$.
\Halmos
\endproof
\section{Results Regarding the Linear Programming-Based FPT Algorithm} \label{appendix:lp}

\subsection{Proof of Propositions \ref{prop:courcelle}, \ref{prop:ackerman}, and \ref{prop:tw.benzwi}} \label{appendix:LP.ackerman}

\proof{Proof of Proposition \ref{prop:courcelle}}
This follows from \citet[Lemma 14]{courcelle2015clique}, the definition of $G'$, noting that $(i,j) \in E'$ if and only if $(j,i) \in E'$, and considering the undirected version of $G'$.
\Halmos
\endproof

The proof of Proposition \ref{prop:ackerman} follows that given in \cite{ackerman2010combinatorial} with some small modifications linked to the specificities of $G'$ (e.g., its bipartiteness, the fact that its edges are weighted, etc.).

\proof{Proof of Proposition \ref{prop:ackerman}}
Due to Corollary \ref{cor:sctm.ltm}, we show that any feasible solution for (\ref{eq:tss.ltm}) is feasible for (\ref{eq:MILP}) and vice-versa, and that the objective functions are equivalent. 

Let $S_0' \cup N_F$ be a feasible solution to (\ref{eq:tss.ltm}). We construct the following solution to (\ref{eq:MILP}): we let $s_i=1$ if $i \in S_0'$ and $s_i=0$ if not, and we let $\ell_{ij}=1$ (resp. $\ell_{ji}=1$) if $i \in N_F$ precedes $j \in N_{SC}$ (resp. $j \in N_{SC}$ precedes $i \in N_F$) in terms of activation. We show that this solution is feasible, constraint-by-constraint. Constraint (\ref{eq:milp.const1}) trivially holds if $i \in S_0$. If $i \notin S_0$, then $s_i = 0$. As $S_0$ leads to full activation of $G'$, there must be some $t$ such that $\sum_{\{ (j,i) \in E_{SC,F} \text{ s.t.\ } j \in S_t \}} r_{i,j} \geq c_i$. By construction, the relevant variables $\ell_{ji}$ must be set to $1$. Hence, $\sum_{ \{(j,i) \in E_{SC,F} \}} r_{i,j} \ell_{ji} \geq c_i$ and (\ref{eq:milp.const1}) holds. A similar argument can be used to show that (\ref{eq:milp.constr2}) holds. Constraint (\ref{eq:milp.constr3}) also holds as $i$ cannot simultaneously activate $j$ and $j$ activate $i$. Likewise, if (\ref{eq:milp.constr4}) were violated, then $\ell_{i_1j_1}+\ell_{j_1i_2}+\ell_{i_2j_2}+\ell_{j_2i_1}=4$, which would imply that $i_1$ activates $j_1$, which activates $i_2$, which activates $j_2$, which activates $i_1$, once again. This is not possible as $i_1$ is already activated. Thus (\ref{eq:milp.constr4}) holds. Finally, the objective function of (\ref{eq:MILP}) is equal to the weight of the seed set, which is identical to that of (\ref{eq:tss.ltm}).

We now consider a feasible solution to (\ref{eq:MILP}). To set $\{\ell_{ij},\ell_{ji}\}$, we can associate a directed acyclic graph on $N'$ with an edge from $i\in N_F$ (resp. $j \in N_{SC}$) to $j \in N_{SC}$ (resp. $i \in N_F$) if $\ell_{ij}=1$ (resp. $\ell_{ji}=1$). The graph is directed (constraint (\ref{eq:milp.constr3})) and acyclic (constraint (\ref{eq:milp.constr4})), as $G'$ is bipartite. Thus, we are able to define a topological ordering on the nodes in $N'$. 
We let $S_0'=\{i \in N_F~|~s_i=1\}$. Consider $t \geq 1.$ We define:
$$S'_{2t}=S'_{2t-1} \cup \{ i \in N_F \backslash S'_{2t-1}~|~\ell_{ji}=1 \Rightarrow j \in S'_{2t-1}, \forall j \in \{j~|~(j,i) \in N_{SC,F}\}\}$$
and
$$S'_{2t+1}=S'_{2t} \cup \{ j \in N_{SC}\backslash S'_{2t}~|~\ell_{ij}=1 \Rightarrow i \in S'_{2t}, \forall i \in \{i~|~(i,j) \in N_{F,SC}\}\}.$$
As $\{\ell_{ij}, \ell_{ji}\}$ define a topological ordering on $N'$, we have $S'_{\infty}=N'$. Furthermore, for $i \in S'_{2t}\backslash S'_{2t-1}$, we have:
$$c_i \leq \sum_{\{j|(j,i)\in E_{SC,F}\}} r_{ji}\ell_{ji}=\sum_{\{j|(j,i)\in E_{SC,F} \cap \ell_{ji=1}\}} r_{ji}\ell_{ji}=\sum_{\{j|(j,i)\in E_{SC,F} \cap S'_{2t-1}\}} r_{ji}=\sum_{\{j~|~(j,i)\in E_{SC,F}\}} r_{ji}x_{j(2t)},$$
where the first inequality is due to (\ref{eq:milp.const1}) and the second equality is due to the definition of $S'_{2t}$. A similar set of inequalities can be derived for $j \in S'_{2t+1} \backslash S'_{2t}$, and thus the activation process defined in this way is exactly that given in (\ref{eq:tss.ltm}). As the objectives of (\ref{eq:tss.ltm}) and (\ref{eq:MILP}) are equivalent, this concludes the proof.
\Halmos
\endproof

\proof{Proof of Proposition \ref{prop:tw.benzwi}}
Constraint (\ref{eq:milp.constr4}) enforces that $\ell_{i_1j_1}$ and $\ell_{i_2j_2}$ are connected in the intersection graph of (\ref{eq:MILP}) for any $i_1,i_2 \in N_F$ and $j_1,j_2 \in N_{SC}.$ There are $n \times m$ pairs $(i,j) \in N_F \times N_{SC}$ and all of the corresponding $nm$ variables $\ell_{ij}$ are connected: this creates a clique in the intersection graph of (\ref{eq:MILP}) of size $nm$. As the treewidth of a clique is $nm-1$ and the treewidth of any subgraph of the intersection graph lower-bounds the treewidth of the intersection graph, the result follows.
\Halmos
\endproof

\subsection{Formulation of the Final BiLP and Proofs of Propositions \ref{prop:eq.10}--\ref{prop:tw.MILP.New} and Theorem \ref{thm:eq.refor}} \label{appendix:LP.fpt}

\proof{Proof of Proposition \ref{prop:eq.10}.}
The two objective functions are identical, so it suffices to show that (\ref{eq:MILP}) and (\ref{eq:MILP_Inter2}) are equivalent. First, let $\big(\{\tilde{\ell}_{ij}, \tilde{\ell}_{ji}\}_{i \in N_F,j \in N_{SC}}, \{\tilde{s}_i\}_{i \in N_F} \big)$ be a feasible solution to (\ref{eq:MILP}). Define a solution to (\ref{eq:MILP_Inter2}): set $s_i=\tilde{s}_i$ for $i=1,\ldots,N_F$ and $\ell_{ij}=\tilde{\ell}_{ij}$ (resp. $\ell_{ji}=\tilde{\ell}_{ji}$) for any $\ell_{ij}$ (resp. $\ell_{ji}$) present in (\ref{eq:MILP_Inter2}) with $i \in N_F$ and $j \in N_{SC}.$ Recall that the variables $\{\tilde{\ell}_{ij}\}_{ij}$ from (\ref{eq:MILP}) define a DAG and that this DAG is connected. For any pair $(i,j) \in X_z'$ for some $z$, with $i,j \in N_F$ or $i,j \in N_{SC}$, we set $\ell_{ij}=1$ (resp.\ $\ell_{ij}=0$) if there exists a directed path from $i$ to $j$ (resp.\ from $j$ to $i$) in the DAG. If neither exists, we assign either $0$ or $1$ to $\ell_{ij}$ only ensuring that $\ell_{ij}+\ell_{ji}=1$. It is easy to see that a solution to (\ref{eq:MILP_Inter2}) thus defined satisfies constraints (\ref{eq:milpInter2.const1}) through (\ref{eq:milpInter2.constr3}). Now consider constraint (\ref{eq:milpInter2.constrend}) and suppose that $i \in N_F, j \in N_{F},$ and $k \in N_{SC}$ wlog (other cases can be treated in the same way). If $\ell_{jk}=\ell_{ki}=1$, then this means that $\tilde{\ell}_{jk}=\tilde{\ell_{ki}}=1$ and the DAG defined by $(\tilde{\ell}_{ij})_{ij}$ contains a directed path from $j$ to $i$. This implies that $\ell_{ji}=1$ and $\ell_{ij}=0$ and, thus, constraint (\ref{eq:milpInter2.constrend}) holds. If either $\ell_{jk}$ or $\ell_{ki} \neq 1$, then the constraint trivially holds. This shows the implication.

Now, suppose we have a feasible solution to (\ref{eq:MILP_Inter2}). For each bag $X_z'$, we can draw a directed graph with nodes in $X_z'$ and an edge from node $i$ to node $j$ in $X_z'$ if $\ell_{ij}=1$. As the constraint (\ref{eq:milpInter2.constrend}) precludes cycles from forming, these graphs are all DAGs. Thus, each bag gives rise to a partial ordering of the nodes it contains. One can then consider a partial ordering on all nodes (as each node appears in at least one bag) obtained by taking the union of the partial orders across bags. Indeed, if two nodes are ordered in a specific way in one bag, the presence of the $\{\ell_{ij}\}$ will enforce the same ordering in any other bag where they both appear. Then, from the order-extension principle, one can define a total order on all nodes in $N_F \cup N_{SC}$ which is consistent with this partial order. That is, one can extend $\{\ell_{i,j}\}_{i,j \in N' \cup X'_z, z \in T'}$ to a sequence $\{\ell_{ij}\}_{i, j \in N'}$ in such a way that this latter set represents a DAG. By taking as our solution to (\ref{eq:MILP}) the appropriate subset $\{\ell_{ij},\ell_{ji}\}_{i \in N_F,j \in N_{SC}}$ of the aforementioned sequence combined with the $\{s_i\}$ given by $(\ref{eq:MILP_Inter2})$, we obtain a feasible solution to (\ref{eq:MILP}) as the no-cycle constraint of (\ref{eq:milp.constr4}) is guaranteed to hold by acyclicity of the DAG.
\Halmos
\endproof

\emph{Formulation of the Final BiLP.}

\begin{subequations} \label{eq:MILP_New}
\begin{alignat}{4}
&\min && \quad  &&\sum_{i \in N_F} w_i s_i &&   \notag \\
&\text{s.t.} &&  &&\forall i \in N_F: &&\sum_{j \in J_z^i} r_{ji} \ell_{ji} \geq u_{iz},~\forall z \in T_i^G, \quad 0\geq u_{iz}, ~\forall z \in T_i^{\tilde{G}} \label{eq:milpNew.rji1}\\
& && && &&u_{iz}+\tilde{u}_{ic_1(z)}+\tilde{u}_{ic_2(z)} \geq \tilde{u}_{iz}, ~\forall z\neq z_0 \in T_i', \label{eq:milpNew.rji2}\\
& && && &&\tilde{u}_{iz_0^i}+ \tilde{u}_{ic_1(z_0^i)}+ \tilde{u}_{ic_2(z_0^i)} \geq c_i(1-s_i),\label{eq:milpNew.rji3}\\
& && && &&u_{iz}=\sum_{\tau=0}^{n_u} 2^\tau u^\tau_{iz} \text{ and } \tilde{u}_{iz}=\sum_{\tau=0}^{n_u} 2^\tau \tilde{u}^\tau_{iz},~ \forall z \in T_i', \label{eq:milpNew.binu}\\
& &&  &&\forall j \in N_{SC}: &&~ \sum_{i \in I_z^j} \ell_{ij} \geq v_{jz}, ~\forall z \in T_j^G, \quad 0 \geq v_{jz}, ~\forall z \in T_j^{\tilde{G}},  \label{eq:milpNew.theta1}\\
& && && && v_{jz}+\tilde{v}_{jc_1(z)}+\tilde{v}_{jc_2(z)} \geq  \tilde{v}_{jz}, ~\forall z\neq z_0 \in T_j', \label{eq:milpNew.theta2}\\
& && && && \tilde{v}_{jz_0^j}+\tilde{v}_{jc_1(z_0^j)}+ \tilde{v}_{jc_2(z_0^j)} \geq \theta_j-1,\label{eq:milpNew.theta3}\\
& && && &&v_{jz}=\sum_{\tau=0}^{n_v} 2^\tau v^\tau_{jz} \text{ and } \tilde{v}_{jz}=\sum_{\tau=0}^{n_v} 2^\tau \tilde{v}^\tau_{jz},~\forall z \in T'_j, \label{eq:milpNew.binv} \\
&  &&  &&\forall z \in T': &&\ell_{ij}+\ell_{ji}=1,~\forall i,j \in X_z' \cap N', \label{eq:milpNew.equality1}\\
& && && &&\ell_{ij}+\ell_{jk}+\ell_{ki} \leq 2,~\forall i,j,k \in X_z' \cap N', \label{eq:milpNew.constrend}\\
& && && &&s_i \in \{0,1\},~\forall i \in N_F, \quad \ell_{ij} \in \{0,1\},~\forall i,j \in N' \cap X_z', \forall z \in T', \notag\\
& && && &&u_{iz}^\tau \text{ (resp. $\tilde{u}_{iz}^\tau$)} \in \{0,1\}, ~\forall \tau, \forall i \in N_F, \forall z \in T' \text{ (resp.\ $\forall z \in T' \backslash \bigcup_{i \in N_F} z_{i}^0$)}, \notag\\
& && && &&v_{jz}^\tau \text{ (resp. $\tilde{v}_{jz}^\tau$)} \in \{0,1\}, ~\forall \tau, \forall j \in N_{SC},\forall z \in T' \text{ (resp.\ $\forall z \in T' \backslash \bigcup_{j \in N_{SC}} z_{j}^0$)}. \notag
\end{alignat}
\end{subequations}

\proof{Proof of Proposition \ref{prop:eq.MIN.SCTM}.}
First, note that (\ref{eq:MILP_New}) is equivalent to the following integer linear program (ILP):
\begin{subequations} \label{eq:MILP_Inter1}
\begin{alignat}{4}
&\min && \quad  &&\sum_{i \in N_F} w_i s_i &&   \notag \\
&\text{s.t.} &&  &&\forall i \in N_F: &&\sum_{j \in J_z^i} r_{ji} \ell_{ji} \geq u_{iz},~\forall z \in T_i^G, \quad 0\geq u_{iz}, ~\forall z \in T_i^{\tilde{G}}  \label{eq:milpInter1.rji1}\\
& && && &&u_{iz}+\tilde{u}_{ic_1(z)}+\tilde{u}_{ic_2(z)} \geq \tilde{u}_{iz}, ~\forall z\neq z_0 \in T_i', \quad \tilde{u}_{iz_0^i}+ \tilde{u}_{ic_1(z_0^i)}+ \tilde{u}_{ic_2(z_0^i)} \geq c_i(1-s_i),\label{eq:milpInter1.rji2}\\
& &&  &&\forall j \in N_{SC}: &&~ \sum_{i \in I_z^j} \ell_{ij} \geq v_{jz}, ~\forall z \in T_j^G, \quad 0 \geq v_{jz}, ~\forall z \in T_j^{\tilde{G}},  \label{eq:milpInter1.theta1}\\
& && && && v_{jz}+\tilde{v}_{jc_1(z)}+\tilde{v}_{jc_2(z)} \geq  \tilde{v}_{jz}, ~\forall z\neq z_0 \in T_j', \quad \tilde{v}_{jz_0^j}+\tilde{v}_{jc_1(z_0^j)}+ \tilde{v}_{jc_2(z_0^j)} \geq \theta_j-1,\label{eq:milpInter1.theta2}\\
&  &&  &&\forall z \in T': &&\ell_{ij}+\ell_{ji}=1,~\forall i,j \in X_z' \cap N', \label{eq:milpInter1.eq1}\\
& && && &&\ell_{ij}+\ell_{jk}+\ell_{ki} \leq 2,~\forall i,j,k \in X_z' \cap N', \label{eq:milpInter1.constrend}\\
& && && &&s_i \in \{0,1\},~\forall i \in N_F, \quad \ell_{ij} \in \{0,1\},~\forall i,j \in N' \cap X_z', \forall z \in T', \notag\\
& && && &&u_{iz} \text{ (resp. $\tilde{u}_{iz}$)} \in \{0,\ldots,c_{\max}\}  ~\forall i \in N_F,~\forall z \in T' \text{ (resp. $\forall z \in T' \backslash \cup_{i \in N_F} z_{i_0}$)} \notag\\
& && && &&v_{jz} \text{ (resp. $\tilde{v}_{jz}$)} \in \{0,\ldots,\theta_{\max}\} ~\forall j \in N_{SC},~\forall z \in T' \text{ (resp. $\forall z \in T' \backslash \cup_{j \in N_{SC}} z_{j_0}$)}, \notag
\end{alignat}
\end{subequations}
where $c_{\max}=\max \{c_1,\ldots,c_n\}$ and $\theta_{\max}=\max \{\theta_1,\ldots,\theta_m\}-1$. Simply note that $\sum_{\tau=0}^{n_u} u_{iz}^\tau$ (resp. $\tilde{u}_{iz}^\tau$) is the binary formulation of $u_{iz}$ (resp. $\tilde{u}_{iz}$), and $\sum_{\tau=0}^{n_v} v_{jz}^\tau$ (resp. $\tilde{v}_{jz}^\tau$) is the binary formulation of $v_{jz}$ (resp. $\tilde{v}_{jz}$). 

We now show that the ILP (\ref{eq:MILP_Inter1}) is equivalent to (\ref{eq:MILP_Inter2}). By virtue of Proposition \ref{prop:eq.10} and the above, the result follows. Assume that constraints (\ref{eq:milpInter1.rji1}) and (\ref{eq:milpInter1.rji2}) hold. By iteratively using constraint (\ref{eq:milpInter1.rji2}) as we go up the tree $T'_i$ from the leaves to the roots, we obtain that $\sum_{z \in T'_i} u_{iz}\geq c_{i}(1-s_i)$. Now, using (\ref{eq:milpInter1.rji1}), it follows that $\sum_{z \in T_{i}^G} \sum_{j \in J_{z}^i}  r_{ji}\ell_{ji} \geq c_i(1-s_i)$. As the sets $\{J_z^i\}_{z\in T_i^G}$ partition $J_i=\{j~|~(j,i) \in E_{SC,F}\}$, we obtain (\ref{eq:milpInter2.const1}). Conversely, if (\ref{eq:milpInter2.const1}) holds, then we can simply set $u_{iz}=\min \{c_{\max},\sum_{j \in J_z^i} r_{ji} \ell_{ji}\}$ if $z \in T_i^G$, or $u_{iz}=0$ if $z \in T_i^{\tilde{G}}$, and $\tilde{u}_{iz}=\min \{c_{\max}, u_{iz}+\tilde{u}_{ic_1(z)}+\tilde{u}_{ic_2(z)}\}$. Then, (\ref{eq:milpInter1.rji1}) and (\ref{eq:milpInter1.rji2}) hold. 

Likewise, assume that constraints (\ref{eq:milpInter1.theta1}) and (\ref{eq:milpInter1.theta2}) hold. By iteratively using constraint (\ref{eq:milpInter1.theta2}) as we go up the tree $T'_j$ from the leaves to the roots, we obtain that $\sum_{z \in T'_j} v_{jz}\geq \theta_j-1$. Now, using (\ref{eq:milpInter1.theta1}), it follows that $\sum_{z \in T_{j}^G} \sum_{i \in I_z^j}  \ell_{ij} \geq \theta_j-1$. As the sets $\{I_z^j\}_{z\in T_j^G}$ partition $I_j=\{i~|~(j,i) \in E_{SC,F}\}$, we obtain (\ref{eq:milpInter2.constr2}). Conversely, if (\ref{eq:milpInter2.constr2}) holds, then we can simply set $v_{jz}=\min \{\theta_{\max}, \sum_{i \in I_z^j} \ell_{ij}\}$ if $z \in T_j^G$, or $v_{jz}=0$ if $z \in T_j^{\tilde{G}}$, and $\tilde{v}_{jz}=\min \{\theta_{\max}, v_{jz}+\tilde{v}_{jc_1(z)}+\tilde{v}_{jc_2(z)}\}$. Then, (\ref{eq:milpInter1.theta1}) and (\ref{eq:milpInter1.theta2}) hold. As the remaining constraints and the objectives are the same, (\ref{eq:MILP_Inter1}) and (\ref{eq:MILP_Inter2}) are equivalent. 
\Halmos
\endproof

\proof{Proof of Proposition \ref{prop:tw.MILP.New}.}
Recall that $\mathcal{T}'=(T',\{X'_z\}_{z \in T'})$ is the tree decomposition of $G'$, where we assume that each node $z$ has no more than two children. We build a tree decomposition $\mathcal{S}=(S,\{W_z\}_{z \in S})$ from $\mathcal{T}'$, where $S=T'$ and each bag $W_z$ contains variables from (\ref{eq:MILP_New}) instead of nodes. We then show that such a tree decomposition is in fact a valid tree decomposition of the intersection graph of (\ref{eq:MILP_New}) and that its width is upper-bounded by $O(\omega'^2+\omega' \log_2(\vartheta_{\max}))$. This proves the result. 
We now specify the bags $\{W_z\}_{z\in S}$:
\begin{enumerate}
            \item Let $z \in S$ and let $X'_z$ be the corresponding bag of nodes in $T'$. Then, $W_z=\bigcup_{i,j \in N' \cap X'_z}\{\ell_{ji}\}.$
            \item For each $i \in N_{F}$, consider $T'_i$ in $T'$ with an arbitrary root node. For each node $z$ in $T'_i$: 
                \begin{enumerate}
                    \item If $z$ is the root node, add variables $\{u^\tau_{iz}\}_{\tau=0,\ldots,n_u}$ and $ s_i$, as well as the variables also present in its children's nodes, $ \{\tilde{u}^\tau_{ic_1(z)}\}_{\tau=0,\ldots,n_u}$ and $\{\tilde{u}^\tau_{ic_2(z)}\}_{\tau=0,\ldots,n_u}$, to $W_z$.
                    \item If $z$ is not a leaf and not the root, add variables $\{u^\tau_{iz}\}_{\tau=0,\ldots,n_u}$ and $\{\tilde{u}^\tau_{iz}\}_{\tau=0,\ldots,n_u}$, as well as the variables also present in its children's nodes, $\{\tilde{u}^\tau_{ic_1(z)}\}_{\tau=0,\ldots,n_u}$ and $\{\tilde{u}^\tau_{ic_2(z)}\}_{\tau=0,\ldots,n_u}$, to $W_z$.
                    \item If $z$ is a leaf, add variables $ \{u^\tau_{iz}\}_{\tau=0,\ldots,n_u}$ and $\{\tilde{u}^\tau_{iz}\}_{\tau=0,\ldots,n_u}$ to $W_z$.
                \end{enumerate}
            \item Likewise, for each $j \in N_{SC}$, consider $T'_j$. For each node $z$ in $T'_j$: 
                \begin{enumerate}
                \item If $z$ is the root node, add variables $\{v^\tau_{jz}\}_{\tau=0,\ldots,n_v}$, as well as the variables also present in its children's nodes, $\{\tilde{v}^\tau_{jc_1(z)}\}_{\tau=0,\ldots,n_v}$ and $\{\tilde{v}^\tau_{jc_2(z)}\}_{\tau=0,\ldots,n_v}$, to $W_z$.
                \item If $z$ is not a leaf and not the root, add variables $\{v_{jz}\}_{\tau=0,\ldots,n_v}$ and $\{\tilde{v}_{jz}\}_{\tau=0,\ldots,n_v}$, as well as the variables also present in its children's nodes, $\{\tilde{v}^\tau_{jc_1(z)}\}_{\tau=0,\ldots,n_v}$ and $\{\tilde{v}^\tau_{jc_2(z)}\}_{\tau=0,\ldots,n_v}$, to $W_z$.
                \item If $z$ is a leaf, add variables $\{v^\tau_{jz}\}_{\tau=0,\ldots,n_v}$ and $\{\tilde{v}^\tau_{jz}\}_{\tau=0,\ldots,n_v}$ to $W_z$.
                \end{enumerate}
\end{enumerate}

We now show that $\mathcal{S}$ as constructed is a valid tree decomposition for the intersection graph of (\ref{eq:MILP_New}). To do this, we need to prove three points. First, all variables involved in the optimization problem appear in at least one of the bags $\{W_{z}\}_{z \in S}$. This is straightforward to check. Second, if a variable appears in two distinct bags, then it appears in all bags in-between. We proceed by groups of variables. The variables $s_i, \{u^\tau_{iz}\}_{\tau}, \{v^\tau_{jz}\}_\tau$ each only appear in one bag, thus this trivially holds for them. The variables $\{\tilde{u}^\tau_{iz}\}_\tau$ and $\{\tilde{v}^\tau_{jz}\}_\tau$ appear in two bags, however, these are parent/children combinations, so the property holds. This leaves variables $\ell_{ji}$. Suppose that $\ell_{ji}$ appears in bag $W_{z_1}$ and $W_{z_2}$ and that there is at least one bag between $W_{z_1}$ and $W_{z_2}$. The assumption implies that $j,i \in X_{z_1}'$ and $j,i \in X_{z_2}'$. As $\mathcal{T}$ is a tree decomposition, it follows that $j$ and $i$ appear in all bags between $X_{z_1}'$ and $X_{z_2}'$, thus $\ell_{ji}$ also appears in all bags between $W_{z_1}$ and $W_{z_2}.$ Third, if a group of variables appears in a constraint, then this group appears in at least one bag of $\mathcal{S}$, because the group forms a clique in the intersection graph. We proceed constraint by constraint. For constraint (\ref{eq:milpNew.rji1}), by construction, $j\in J_{z}^i \subseteq X_{z}'$ for $z \in T_i^G$, thus $i,j \in X_z'$, and $\ell_{ji} \in W_{z}.$ Furthermore, from Step 2b in the construction of $\mathcal{S}$, $\{u_{iz}^\tau\}_\tau \in W_z.$ For constraints (\ref{eq:milpNew.rji2}) and (\ref{eq:milpNew.rji3}), this is straightforward from steps 2a and 2b. A similar reasoning applies to constraints (\ref{eq:milpNew.theta1}), (\ref{eq:milpNew.theta2}), and (\ref{eq:milpNew.theta3}). For constraints (\ref{eq:milpNew.constrend}), this follows from step 1. Thus, $\mathcal{S}$ is a valid tree decomposition for the intersection graph of (\ref{eq:MILP_New}). 

We now upper-bound the treewidth by looking at the size of each one of the bags $\{W_z\}_{z \in S}.$ Consider the algorithm to build $\mathcal{S}$ and recall that $\omega'$ is the treewidth of $G'$, and thus the maximum size of $X_z'$ for all $z \in T.$ Step 1 only occurs once and at the end of it, $W_z$ contains at most $\omega'^2$ nodes. Then, for each node $z \in T'$, steps 2-3 happen at most a combined $\omega'$ times as $X'_z$ contains at most $\omega'$ nodes. Thus, $z$ appears in at most $\omega'$ trees $T'_i$ or $T'_j$. During Step 2 (resp.\ Step 3), a maximum of 4 sets of variables are added to the node, with each set having size at most $\max \{n_u, n_v\}$, where $n_u, n_v$ are as defined in (\ref{eq:nu.nv}). Thus, at the end of the construction of $\mathcal{S}$, $W_z$ contains at most $\omega'^2+4\max \{n_u, n_v\} \cdot \omega'=O(\omega'^2+\omega' \log_2(\vartheta_{\max}))$ nodes.
\Halmos
\endproof

\proof{Proof of Theorem \ref{thm:eq.refor}.}
From Proposition \ref{prop:eq.MIN.SCTM}, we have that (\ref{eq:MILP_New}) solves $MIN$-$SCTM$. From Proposition \ref{prop:tw.MILP.New}, the intersection graph of (\ref{eq:MILP_New}) has treewidth at most $O(\omega'^2+\omega' \log_2(\vartheta_{\max})).$ We now count the variables appearing in (\ref{eq:MILP_New}). Recall that $T'$ has at most $4(n+m)$ bags. We have that (\ref{eq:MILP_New}) has $n$ variables $s_i$, at most $4(n+m) \cdot \omega'^{3}$ variables $\ell_{ij}$, at most $3 \cdot n_u \cdot 4(n+m)$ variables $\{u_{iz}^\tau\}$ and $\{\tilde{u}_{iz}^\tau\}$, and at most $3 \cdot n_v \cdot 4(n+m)$ variables $\{v_{jz}^\tau\}$ and $\{\tilde{v}_{jz}^\tau\}$. From Propositions \ref{prop:courcelle} and \ref{prop:binLP}, we obtain the result.
\Halmos
\endproof

\subsection{Proof of Proposition \ref{prop:kappa}} \label{appendix:kappa}

\proof{Proof of Proposition \ref{prop:kappa}.} When $\kappa=0$, (\ref{eq:LP.k}) reads:
\begin{equation*} 
    \begin{aligned}
    \min_{Y_\mathbb{S}}~~~~~~&\sum_i w_i Y_{\{s_i\}}\\
    \text{s.t.}~~~~~~~&Y_{\emptyset}=1\\
    \forall z \in S:~ &Y_{\emptyset} -Y_{\{x_i^s\}} \geq 0, Y_{\{x_i^s\}} \geq 0,~\forall i=1,\ldots,\omega_z\\
    & g_{\emptyset}^l Y_{\emptyset}+ \sum_{i=1}^{\omega_z} g_{\{x_{i}^z\}}^l \cdot Y_{\{x_i^z\}} \geq 0,~\forall l=1,\ldots,l_z,
\end{aligned}
\end{equation*}
By definition of an intersection graph, all constraints of (\ref{eq:MILP_New}) are associated to at least one node $z \in S$. The result follows.
\Halmos
\endproof

\section{A Dynamic Programming-Based FPT Algorithm} \label{appendix:dp}

Corollary \ref{cor:sctm.ltm} shows that if we can solve (\ref{eq:tss.ltm}) for $G'$, we can solve $MIN$-$SCTM$ for $G$. The algorithm introduced by \citet{ben2011treewidth} is an FPT algorithm (with parameter $\omega'$) for target set selection under the LTM, assuming thresholds are bounded.  
Thus, we show how to generalize the algorithm to our setting, taking into consideration the characteristics of the auxiliary graph $G'$: (i) edges are weighted, (ii) only certain types of nodes (the firm-nodes) can be part of the seed set, (iii) the objective is minimizing costs rather than number of nodes in the seed set, and (iii) we are interested in structural assumptions on $G$ more so than $G'.$

The algorithm first computes a tree decomposition of $G'$ with treewidth $\omega'$ and of a specific type:
\begin{definition}[Nice tree decomposition]
Let $G' = (N',E')$. A tree decomposition $\mathcal{T}'=(T',(X'_z)_{z \in Z})$ of $G'$ with treewidth $\omega'$ is nice if and only if $T'$ is rooted at some node $\tilde{z}$, $|X_z'| = \omega'+1\ \forall z \in T'$, and all nodes are of exactly one of the following types: (a) leaf nodes, (b) replace nodes ($z$ has exactly one child, $z_0$, and there are $u, v \in N'$, $u \neq v$, such that $X_z' \setminus X_{z_0}' = \{u\}$ and $X_{z_0}' \setminus X_{z}' = \{v\}$), and (c) join nodes ($z$ has exactly two children, $z_0$ and $z_1$, and $X_z' = X_{z_0}' = X_{z_1}'$).
\end{definition}
An example of a nice tree decomposition is given in Figure \ref{fig:treedec_nice}.

\begin{figure}[t]
    \centering
    \subfloat[Rearranged auxiliary graph $G'$ from Figure \ref{fig:sctm_ltm_2}.]{
        \includegraphics{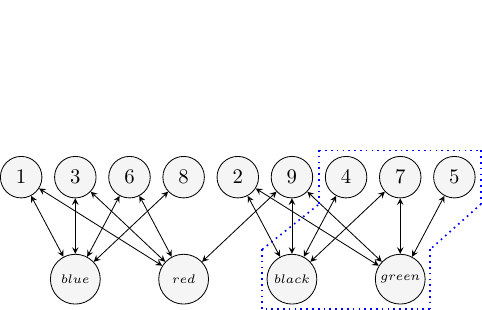}
        \label{fig:nice_a}
        \vphantom{
            \includegraphics{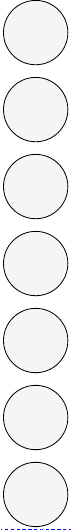}
        }
    }
    \hspace{-0.2cm}
    \subfloat[Nice tree decomposition $T$ of $G'$.]{
        \includegraphics{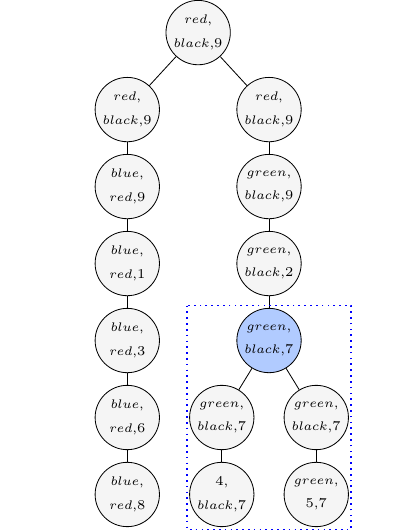}
        \label{fig:nice_b}
    }
    \caption{Example of a nice tree decomposition. Let $r_{ji}=c_i=\theta_j=1$ for all nodes $i$ and hyperedges $e_j$.}\label{fig:treedec_nice}
\end{figure}

Before any formal results, we provide high-level intuition as to how the algorithm works. 
Assume a nice tree decomposition $\mathcal{T}'=(T',(X'_z)_{z \in Z})$ of $G'$ with treewidth $\omega'$. Possibly overlapping subgraphs of $G'$ are created from the leafs of $T'$, each containing at most $\omega'$ nodes. Moving up one level in the tree corresponds to growing the subgraphs constructed at the previous level, either by adding nodes or merging two subgraphs. When we get to the root node, the subgraph considered is the complete graph $G'$. The key property of the subgraphs constructed via this process is that each one only interacts with its complement in $G'$ via a set of \emph{boundary nodes}, which is of size at most $\omega'$. In other words, if there exists an edge between node $i$ in the subgraph and a node outside of the subgraph, then it must be that $i$ belongs to the boundary nodes, and there cannot be more than $\omega$ such nodes. The minimum cost seed set for $G'$ is then obtained recursively using these subgraphs. First, the minimum cost seed sets for the subgraphs corresponding to the leafs are found through enumeration (which is exponential in $\omega'$ but not necessarily in $n$). Then, as we go up the tree, the algorithm uses the previously established minimum cost seed sets for smaller subgraphs to obtain that of the larger subgraph. The computational effort is limited because nodes in the smaller subgraph that are not part of the boundary are not affected by any nodes that may enter the subgraph. Thus, we can focus on nodes in the boundary. As we reach the root node, the relevant subgraph is all of $G'$, allowing us to derive a minimum cost seed set for the full auxiliary graph and, thus, the original graph $G$.

\begin{algorithm2e}[!ht]
 \caption{Algorithmic solution to minimum cost seed set problem. \label{alg:solving}}
 \SingleSpacedXI
\KwData{Auxiliary graph $G'$ with $tw(G') = \omega'$.}
\KwResult{Minimum cost seed set $(S_0')^*$ of $G'$.}
\textbf{Initialization:} Compute a nice tree decomposition $\mathcal{T}'=\left(T',\{X_z'\}_{z\in T'} \right)$ of $G'$, rooted at $\tilde{z}$ with width $\omega'$;
$\tilde{C} \leftarrow \{0,\ldots,\vartheta_{max}\}^{\omega'+1}$; $\tilde{A} \leftarrow \{0,\ldots,\omega'+1\}^{\omega'+1}$; $red \leftarrow 0^{\omega'+1}$\;
\For{$z \in T'$ where $z$ is a leaf node}{
    \For{$\tilde{c} \in \tilde{C}$}{
        \For{$\tilde{a} \in \tilde{A}$}{
            $S_0^z[\tilde{c},\tilde{a}] \leftarrow$ compute minimum cost seed set through enumeration\;
        }
    }
}
\While{$z \in T'$ has not been traversed, but all its child nodes have}{
    \If{$z$ is a replace node with child $z_0$}{
        $G^{z}_0 = (G^z_0,E^z_0) \leftarrow G^z$;
        $E^z_0 \leftarrow E^z_0 \setminus \{(u,i'), (i',u):i' \in X_z'\}$\;
        \For{$\tilde{c} \in \tilde{C}$}{
            $c \leftarrow \tilde{c}$; $c(v) \leftarrow c_v'$\;
            \For{$\tilde{a} \in \tilde{A}$}{
                $\hat{A} \leftarrow \{a \in \tilde{A}:a(i') = \tilde{a}(i')\ \forall i' \neq v \}$; $a^* \leftarrow \arg \min_{a \in \hat{A}} \mathcal{C}\left(S_0^{z_0}[c,a]\right)$\;
                $S_0^z[\tilde{c},\tilde{a}] \leftarrow
                \begin{cases}
                    S_0^{z_0}[c,a^*] & \text{ if } \tilde{c}(u) = 0 \\
                    S_0^{z_0}[c,a^*] \cup \{u\} & \text{ if } \tilde{c}(u) \neq 0, u \text{ is a firm-node} \\
                    \{i' \in N'\ :\ i' \text{ is a firm-node}\} & \text{ otherwise.}
                \end{cases}$\;
            }
        }
        \For{$e_1 = (u,i') \in \{(u,i'):i' \in X_z'\}$ and $e_2 = (i',u)$}{
            $E^z_0 \leftarrow E^z_0 \cup \{e_1,e_2\}$;
            $S_0^{z'} \leftarrow S_0^z$\;
            \For{$\tilde{c} \in \tilde{C}$}{
                $\tilde{c}^u \leftarrow \tilde{c}$;
                $\tilde{c}^{i'} \leftarrow \tilde{c}$;
                $\tilde{c}^u(u) \leftarrow \max\{\tilde{c}(u) - w'_{i',u},0\}$;
                $\tilde{c}^{i'}(i') \leftarrow \max\{\tilde{c}(i') - w'_{u,i'},0\}$\;
                \For{$\tilde{a} \in \tilde{A}$}{
                    $S_0^z[\tilde{c},\tilde{a}] \leftarrow \begin{cases}
                        S_0^{z'}[\tilde{c},\tilde{a}] & \text{ if } \tilde{a}(i') = \tilde{a}(u)\\
                        S_0^{z'}[\tilde{c}^u,\tilde{a}] & \text{ if } \tilde{a}(i') < \tilde{a}(u)\\
                        S_0^{z'}[\tilde{c}^{i'},\tilde{a}] & \text{ if } \tilde{a}(i') > \tilde{a}(u)\\
                    \end{cases}$\;
                }
            }
        }
    }
    \If{$z$ is a join node with children $z_0$ and $z_1$}{
        \For{$\tilde{a} \in \tilde{A}$}{
            \For{$i' \in X_z'$}{
                $red(i') \leftarrow \sum_{\{j' \in X_z\ :\ (j',i') \in E' \text{ and } \tilde{a}(j') < \tilde{a}(i') \}} w_{j',i'}'$\;
            }
            \For{$\tilde{c} \in \tilde{C}$}{
                $(\tilde{c}^{z_0}, \tilde{c}^{z_1}) \leftarrow \arg \min_{\left\lbrace \tilde{c}^{z_0'},\tilde{c}^{z_1'}\ :\ \tilde{c}^{z_0'} + \tilde{c}^{z_1'} = \tilde{c} + red \right\rbrace} \mathcal{C}\left(S_0^{z_0}[\tilde{c}^{z_0'},\tilde{a}] \cup S_0^{z_1}[\tilde{c}^{z_1'},\tilde{a}]\right)$\;
                $S_0^z [\tilde{c},\tilde{a}] \leftarrow S_0^{z_0}[\tilde{c}^{z_0},\tilde{a}] \cup S_0^{z_1}[\tilde{c}^{z_1},\tilde{a}]$\;
            }
        }
    }
}
$(S_0')^* \leftarrow S^{\tilde{z}}_0[c,a]$, with $c^*(i') = c_{i'}'\ \forall i' \in X_{\tilde{z}}'$ and $a = \arg \min_{\tilde{a} \in \tilde{A}} \mathcal{C}\left(S^{\tilde{z}} _0[c^*,\tilde{a}]\right)$\;
\end{algorithm2e}

\begin{lemma} \label{lem:DP_lem1}
Assume an auxiliary graph $G'$ as defined in Section \ref{sec:aux_graph}, as well as the LTM activation process. Algorithm \ref{alg:solving} applied to $G'$ returns a minimum cost seed set.
\end{lemma}

\proof{Proof.}
We introduce some notation.
Consider a subtree in $T'$ rooted at a node $z_0$. Define with $G^{z_0} = (N^{z_0},E^{z_0})$ the subgraph of $G'$ induced by $\bigcup_{z \in subtree(z_0)} X_z'$. Nodes in $X_{z_0}'$ may be connected to nodes in $G'$ outside of $G^{z_0}$, but other nodes in $G^{z_0}$ cannot be connected to those outside. Thus, denote $X_{z_0}'$ as the \emph{boundary} of $z_0$ (it may include unconnected nodes, due to the requirement that $|X_z'| = \omega'+1$). For example, in Figures \ref{fig:nice_b} (resp.\ Figure \ref{fig:nice_a}), the dashed line indicates the subtree rooted at $\{green,black,7\}$ (resp.\ the associated subgraph).
Nodes $black$ and $green$, which are in the boundary, connect to nodes $2$ and $9$ outside of the dashed line. The boundary may include additional nodes such as $7$ due to $|X_z'| = \omega'+1$. Node $4$, on the other hand, is not in the boundary and cannot connect to nodes outside of the dashed line.
Moving upwards in the tree from $z_0$ to a replace node $z_R$, a node $v \in X_{z_0}'$ is replaced by a node $u$ to arrive at the new boundary $X_{z_R}'$. By construction, node $v$ cannot share edges with any of the nodes outside $G^{z_0}$. Say, we move from $\{green,black,7\}$ to $\{green,black,2\}$. In the boundary, Node $7$ is replaced by Node $2$, and Node $7$ does indeed not share edges with nodes outside of the subgraph induced by $black$, $green$, $2$, $4$, and $5$.

Next, let $\tilde{c} \in \{0,\ldots,\vartheta_{\max}\}^{\omega'+1}$ and $\tilde{a} \in \{0,\ldots,\omega'+1\}^{\omega'+1}$ be threshold and activation vectors, where $\vartheta_{\max} = \max \left\lbrace \max_{i =1,\ldots,n} \{ c_{i}\}, \max_{j=1,\ldots,m} \{ \theta_{j}\} -1 \right\rbrace$. We arbitrarily assign a one-to-one mapping from the root node boundary to these vectors, denoting with $\tilde{c}(i')$ (resp.\ $\tilde{a}(i')$) the mapping from $i'$ in the boundary. Define mappings for other boundaries recursively: (i) if $v$ replaces $u$, $v$ is mapped to the same index and the mapping remains unchanged otherwise; and (ii) if a node has two children, boundary and mapping remain unchanged.

Finally, for any $z \in T'$, we define a matrix $S_0^z$ with rows indexed by the vectors $\tilde{c}$ and columns indexed by the vectors $\tilde{a}$. The size of the matrix is $[(\vartheta_{\max}+1) \cdot (\omega'+2)]^{\omega'+1} = [\vartheta_{\max} \cdot \omega']^{O(\omega')}$. Each entry is a set of nodes in $G'$ and matrices will be computed recursively bottom-up. The last matrix to be computed (corresponding to the root node $\tilde{z}$) gives us the minimum cost seed set of $G'$, $(S_0')^*$. 
For any set $S_0^z[\tilde{c},\tilde{a}]$, we define its cost by $\mathcal{C}\left(S_0^z[\tilde{c},\tilde{a}]\right)= \sum_{i' \in S_0^z[\tilde{c},\tilde{a}]} w_i$, where $w_i$ is the cost of adding node $i \in N_F$ to the seed set.

Algorithm \ref{alg:solving} assumes a nice tree decomposition $\mathcal{T}'$ of $G'$. It is based on \citet{ben2011treewidth}, but adapted to reflect the specificities of $G'$.
For each leaf node $z$, each $\tilde{c}$, and each $\tilde{a}$, we compute a ``minimum cost seed set'' (Lines 2--5 in Algorithm \ref{alg:solving}) based on the subgraph $G^z$ and assuming that a node $i' \in X_z'$ (i) has threshold $\tilde{c}(i')$; (ii) can activate only if all nodes $j' \in X_z'$ with $\tilde{a}(j') < \tilde{a}(i')$ are already active and if all $j' \in X_z'$ with $\tilde{a}(j') = \tilde{a}(i')$ activate at the same time; and (iii) can only be part of the seed set if it is a firm-node.
For example, let $z_0 = \{4,black,7\}$, $\tilde{c} = (1,1,0)$, and $\tilde{a} = (2,1,0)$. Node $7$ has a threshold of $0$. It also has the lowest value in the activation vector. Hence, Node $7$ activates for any seed set. Node $black$ has a threshold of $1$. However, there is an edge $(7,black)$ with $w_{7,black}' = 1$ and Node $black$ has a lower value in the activation vector than Node $4$. Once Node $7$ is active, Node $black$ also activates. Thereafter, Node $4$ can also activate. It follows that the minimum cost seed set is $\emptyset$.
Assume, instead, that $\tilde{c} = (1,1,3)$, and $\tilde{a} = (2,1,0)$. Node $7$ still has to activate first. However, its threshold is $3$, so it must be part of the seed set. With Node $7$ in the seed set, the remainder of the activation process is unchanged, so the minimum cost seed set is $\{7\}$.
Finally, assume that $\tilde{c} = (1,1,3)$, and $\tilde{a} = (2,0,1)$. Node $black$ has to activate first, but SC-nodes cannot be added to the seed set. Hence, no such minimum cost seed set exists, and we adopt the convention of equating the minimum cost seed set to the entire set of node-nodes. That is, $S_0^{z_0}[(1,1,3),(2,0,1)] = \{i' \in N'\ :\ i' \text{ is a firm-node} \}$.

Next, we recursively generate $S_0^{z}$ for all other nodes $z \in T'$.
If $z$ is a replace node (Lines 7--19 in Algorithm \ref{alg:solving}) with child $z_0$, the subgraph $G^z$ of $G'$ has exactly one more node than $G^{z_0}$, node $u$. Node $u$ replaces node $v$ in the boundary, mapped to entry $i_v$ in both activation and threshold vectors. For each $\tilde{a}$, define set $\hat{A}= \{ a \in \{0,\ldots,\omega'+1\}^{\omega}: a(i') = \tilde{a}(i')\ \forall i' \neq v \}$ grouping all activation vectors that differ only in entry $i_v$.
For example, if $z_0 = \{4,black,7\}$, $z = \{green,black,7\}$, then $u = green$ and $v=4$. Moreover, if $\tilde{a} = (2,1,0)$, then $\hat{A} = \{ (0,1,0), (1,1,0), (2,1,0), (3,1,0)\}$.
Then, for each $\tilde{c}$, construct intermediary sets $S_0^z[\tilde{c},\tilde{a}]$ to remove dependencies on $v$ and introducing the node $u$ (Lines 9--13 in Algorithm \ref{alg:solving}). We define them thus: (i) if $\tilde{c}(u) = 0$, then $S_0^z[\tilde{c},\tilde{a}] =S_0^{z_0}[c,a]$ where $c(i') = \tilde{c}(i'), i' \neq v$, $c(v) = c'_v$, and $a= \arg \min_{a \in \hat{A}} \mathcal{C}\left(S_0^{z_0}[c,a]\right)$, (ii) if $\tilde{c}(u) > 0$ and $u$ is a firm-node, then $S_0^z[\tilde{c},\tilde{a}] = S_0^{z_0}[c,a] \cup \{u\}$ with $a$ and $c$ as before, (iii) if $\tilde{c}(u) > 0$ and $u$ is a SC-node, then $S_0^z = \{i' \in N'\ :\ i' \text{ is a firm-node}\}$.
In the previous example, if $\tilde{a} = (2,1,0)$ and $\tilde{c} = (0,1,3)$, then $c = (1,1,3)$ and $S_0^z[(2,1,0),(0,1,3)]$ is the least costly of the following seed sets: $\{4,7\}$ (corresponding to $a = (0,1,0)$), $\{4,7\}$ (corresponding to $a = (1,1,0)$), $\{7\}$ (corresponding to $a = (2,1,0)$), $\{7\}$ (corresponding to $a = (3,1,0)$). As $\tilde{c}(green) = 0$, it follows that $S_0^z[(2,1,0),(0,1,3)] = \{7\}$.

Given $\tilde{c}$ and $\tilde{a}$, $S_0^z[\tilde{c},\tilde{a}]$ reflects the minimum cost seed set to activate the subgraph $G^z$ minus any edges between node $u$ and other nodes in the boundary. We iteratively refine the intermediary sets while adding these edges to the subgraph (Lines 14--19 in Algorithm \ref{alg:solving}).
Choose one such edge, $(u,i')$, and the corresponding edge $(i',u)$ and copy the entries of matrix $S_0^z$ to a new matrix $S_0^{z'}$. For any $\tilde{a}$ and $\tilde{c}$, there are three options: (i) if $\tilde{a}(u) < \tilde{a}(i')$, node $u$ has to activate before node $i'$. Once $u$ is active, it contributes $w'_{i',u}$ towards activation of $i'$. Hence, update $S_0^z[\tilde{c},\tilde{a}] = S_0^{z'}[\tilde{c}^{i'},\tilde{a}]$, where $\tilde{c}^{i'}$ is identical to $\tilde{c}$, except in entry $i'$: $\tilde{c}^{i'}(i') = \max\{\tilde{c}(i') - w'_{u,i'},0\}$. (ii) if $\tilde{a}(u) > \tilde{a}(i')$, $u$ has to activate after $i'$. Once $i'$ is active, it contributes a benefit of $w'_{i',u}$ towards activation of $u$. Hence, update $S_0^z[\tilde{c},\tilde{a}] = S_0^{z'}[\tilde{c}^{u},\tilde{a}]$, where $\tilde{c}^{u}$ is identical to $\tilde{c}$, except in entry $u$: $\tilde{c}^{u}(u) = \max\{\tilde{c}(u) - w'_{u,i'},0\}$. (iii) if $\tilde{a}(u) = \tilde{a}(i')$, the nodes have to activate simultaneously, so they cannot influence each other. Thus, $S_0^z[\tilde{c},\tilde{a}] = S_0^{z'}[\tilde{c},\tilde{a}]$.
In the previous example, add the edges $(green,7)$ and $(7,green)$ with $w'_{green,7}=w'_{7,green}=1$, and assume again that $\tilde{a} = (2,1,0)$ and $\tilde{c} = (0,1,3)$. Then, $\tilde{a}(green) = 2 > 0 = \tilde{a}(7)$ and we consider $\tilde{c}^{green} = (0,1,3)$. The latter vector is unchanged because the threshold in the second entry is already $0$, so $S_0^z[(0,1,3),(2,1,0)] = S_0^{z'}[(0,1,3),(2,1,0)]$. If, however, $\tilde{a} = (0,1,2)$, then we need to consider the threshold vector $\tilde{c}^7 = (0,1,2)$ to account for the fact that $7$ always activates after $green$. In this case, $S_0^z[(0,1,3),(0,1,2)] = S_0^{z'}[(0,1,2),(0,1,2)]$.
Revise $S_0^{z'} = S_0^z$ and repeat for all edges.

If $z$ is a join node (Lines 20--26 in Algorithm \ref{alg:solving}) with children $z_0$ and $z_1$, $G^z = G^{z_0} \cup G^{z_1}$ (there are no edges between nodes of the two subgraphs outside the boundary). Fix $\tilde{a}$ and $\tilde{c}$. A node $i' \in X_z'$ may benefit from activations in both subgraphs $G^{z_0}$ and $G^{z_1}$. For example, consider the join node with $z=(red,black,9)$ in Figure \ref{fig:treedec_nice} and take $\tilde{a} = (2,0,1)$ (i.e., $black$ has to activate first, followed by $9$, then $red$) and $\tilde{c}=(1,0,2)$ (i.e., $\tilde c_{red} = 1$, $\tilde c_{black} = 0$, $\tilde c_{9} = 2$). In this case, the minimum cost seed set for the subgraph in $G'$ associated with the left-hand subtree is $\{9\}$: $black$ activates due to its threshold, $red$ activates as $9$ is active and then $1$, $3$ and $6$ activate, followed by $blue$ and finally $8$. Following a similar reasoning, a minimum cost seed set for the subgraph in $G'$ associated to the right-hand subtree is $\{5\}$. However, while $G'$ does activate if $\{5,9\}$ is active, this is not the smallest seed set, which is given by $\{9\}.$
Thus, we construct $S_0^{z}[\tilde{c},\tilde{a}]$ to account for synergies. For each $i' \in X_z'$, we define the following weight reduction which avoids double-counting: $red(i') = \sum_{\{j' \in X_z'\ :\ (j',i')\in E' \text{ and } \tilde{a}(j') < \tilde{a}(i')\}} w'_{i',j'}$. We then take $S_0^z [\tilde{c},\tilde{a}] = S_0^{z_0}[\tilde{c}^{z_0},\tilde{a}] \cup S_0^{z_1}[\tilde{c}^{z_1},\tilde{a}]$, where
\begin{align*}
    (\tilde{c}^{z_0},\tilde{c}^{z_1})=\arg &\min_{f,g \in \{0,\ldots,k-1\}^{\omega}} \mathcal{C}\left(S_0^{z_0}[f,\tilde{a}] \cup S_0^{z_1}[g,\tilde{a}]\right)\\
    &\text{s.t. } f(i') + g(i') = \tilde{c}(i') + red(i') \text{ for all }i' \in X_z'.
\end{align*}
For example, if $X_{z_0}' = X_{z_1}' = \{green,black,7\}$, $\tilde{a} = (2,1,0)$, and $\tilde{c} = (0,1,3)$, then $red(black) = 1$, $red(green) = 1$, and $red(7) = 0$. It follows that $S^z_0[(0,1,3),(2,1,0)]$ is based on the union of seed sets corresponding to the threshold vectors $(\tilde{c}^{z_0},\tilde{c}^{z_1})=\arg \min_{\left\lbrace f, g\ :\ f + g = (1,2,3) \right\rbrace}  \mathcal{C}\left(S_0^{z_0}[f,(2,1,0)] \cup S_0^{z_1}[g,(2,1,0)] \right)$.

Proceed until the root note $\tilde{z}$. To obtain $(S_0')^*$, take the threshold vector $c^*$ corresponding to the actual thresholds, that is $c^*(i') = c'_{i'}\ \forall i' \in X_{\tilde{z}}'$. As the optimal seed set induces an activation sequence $\tilde{a}$, $(S_0')^* = S^{\tilde{z}}_0[c^*,a^*]$, with $a^* = \arg \min_{\tilde{a}} S^{\tilde{z}} _0[c^*,\tilde{a}]$.
\Halmos
\endproof

\begin{lemma} \label{lem:DP_lem2}
Assume an auxiliary graph $G'$ as defined in Section \ref{sec:aux_graph} with treewidth $tw(G') = \omega'$. Algorithm \ref{alg:solving} applied to $G'$ runs in $[ \vartheta_{\max} \cdot \omega' ]^{O(\omega')} (n+m)$ time.
\end{lemma}

\proof{Proof.}
The decomposition of $G'$ into a minimal tree requires time exponential in $\omega'$ but linear time when $\omega'$ is bounded \citep{bodlaender1996linear} and, given an arbitrary tree decomposition, one can always construct a a nice tree decomposition of the same width in linear time. \citep{ben2011treewidth}.
Moreover, $G'$ contains $n+m$ nodes, so a nice tree decomposition of width $\omega'$ exists with at most $(\omega'+1) \cdot (n+m)$ nodes. Hence, the number of entries of any $S_0^z$ is bounded by $[\vartheta_{max}  \omega']^{O(\omega')}$.

The number of computations for each entry of a leaf node is $2^{\omega'+1}$ (Lines 2--5 in Algorithm \ref{alg:solving}). The number of computations for each entry of a replace node is determined by comparing all $\omega'+2$ activation options of $v$ (Lines 8--13 in Algorithm \ref{alg:solving}) and iterating through each of at most $\omega'$ edges that $u$ shares with other nodes in the boundary (Lines 14--19 in Algorithm \ref{alg:solving}). The number of computations for each entry of a join node is determined by comparing combinations of thresholds of the boundaries, which is upper-bounded by a constant factor of $\vartheta_{max}^{\omega'+1}$. It follows that the maximum number of computations required for Algorithm \ref{alg:solving} is in $ \left[\vartheta_{max}  \omega' \right]^{O(\omega')} (n+m)$.
\Halmos
\endproof

We can now state and proof the main result regarding the DP algorithm:

\begin{theorem}\label{theo:gen-solution}
Let $G$ be a hypergraph as defined in Section \ref{sec:SCTM}.
Assume that $tw(G)=\omega$. Then, $MIN$-$SCTM$ can be solved exactly via a dynamic program in 
$[ \vartheta_{\max} \cdot \omega ]^{O(\omega)} (n+m)$
time.
\end{theorem}

\proof{Proof of Theorem \ref{theo:gen-solution}.}
From Lemma \ref{lem:DP_lem1}, we know that Algorithm \ref{alg:solving} provides a minimum cost seed set for auxiliary graph $G'$, assuming the LTM activation process. Meanwhile, Line 13 of the algorithm ensures that the seed set contains only nodes $i \in N_F$ and the assumption on how seed set costs are computed ensures that these nodes are weighted with $w_i$. It follows from Proposition \ref{prop:sctm.eq.ltm} that the seed set identified is also a solution to $MIN$-$SCTM$. Combined with Lemma \ref{lem:DP_lem2}, the fact that constructing $G'$ from $G$ requires a constant factor of $n+m$ operations, and Proposition \ref{prop:courcelle} ($\omega'  \leq \omega +1$), the result follows.
\Halmos
\endproof

From this, we can directly derive a corollary corresponding to Corollary \ref{cor:lp}:
\begin{corollary} \label{cor:dp}
There is a dynamic programming-based FPT algorithm for solving $MIN$-$SCTM$ with parameter $\omega$ and $\vartheta_{\max}$, running in time at most
$2^{O(\omega \log_2(\omega))} \vartheta_{\max}^{O(\omega)} (n+m).$
\end{corollary}

\section{Description of data and measures for numerical experiments} \label{appendix:empirical}

\subsection{Description of supply chain network data sets} \label{appendix:data}

We generate three datasets of supply chain networks, with key statistics summarized in Table \ref{tab:data}.
We discuss the generation processes below. In all cases, we let $c_i = r_{ji} = 1$, $\theta_j = k_j$, and $w_i \sim N(1,0.1)$ for all $i \in N_F, j \in N_{SC}$. Standardization allows us to remove (predictable) effects of costs and benefits while focusing on the effects of graph structures. At the same time, perturbing $w_i$ enables us to obtain unique solutions even in networks with repetitive structures.

\begin{table}[t]
\centering
\caption{Means and standard deviations of key statistics of each data set employed.}
\label{tab:data}
\begin{tabular}{lccc}
\hline 
 \up\down
 & \citet{Willems2008} & Random & Modular \\
\hline
\up
\multirow{2}{*}{$n$} & 152.54 & 24.23 & 29.91 \\
    & (178.38) & (0.77) & (0.38)\\[0.1cm]
\multirow{2}{*}{$m$} & 157.33 & 27.54 & 54.57\\
                    & (337.27) & (10.27) & (19.32) \\[0.1cm]
\multirow{2}{*}{$k$} & 5.34 & 5.0 & 5.0 \\
                        & (1.98) & (0.0) & (0.0) \\[0.1cm]
\multirow{2}{*}{$\omega'$} &  8.96 & 9.62 & 14.85 \\
& (9.02) & (2.72)  & (4.25)\\
\hline
\end{tabular}
\end{table}

\subsubsection*{\citet{Willems2008} networks.}

The original data set represents 38 acyclic networks of companies in 22 industries. While the data contains information about direct connections between nodes, it does not specify which sets of nodes belong to which supply chains (i.e., we do not have access to the hyperedges of $G$). To circumvent this difficulty, we randomly generate the hyperedges of the graph using the observable edges in the following way: We consider all possible paths between nodes in the first and the last tier, that is, nodes with no incoming, respectively, no outgoing edges. We then assume for each path with probability 0.05, 0.25, or 0.5 that it is a supply chain and remove nodes not part of any supply chain. 
We repeat this generation process ten times for each network-probability combination and remove networks with less than 15 remaining nodes, as well as those for which we cannot find a seed set guaranteed to be within 5\% of the minimum cost within two hours using Gurobi. This results in a total of 657 supply chain networks.

\subsubsection*{Random networks.}

We fix five tiers and the following configurations of nodes per tier: $(5,5,5,5,5)$, $(2,6,9,6,2)$, and $(2,2,4,7,10)$, where the value at index $l$ represents the number of nodes at tier $l$. We also fix a number of supply chains $\overline{m} \in \{20,40\}$. Each node and supply chain is randomly assigned a value drawn from the uniform distribution on $(0,1)$.
Then, at each tier, we identify the $h \in \{1,2,\ldots,6\}$ nodes whose assigned value is closest to that of the supply chain.
Of those nodes, we randomly choose one to be part of the supply chain. We repeat the generation process if a supply chain with the same nodes already exists and remove nodes that belong to no supply chain after all have been generated. The instance is discarded if more than 10\% of the generated nodes have been removed.
We repeat the entire process 20 times for each parameter combination and discard networks for which we cannot find a seed set guaranteed to be within 5\% of the minimum cost within two hours using Gurobi. This results in 525 supply chain networks.

\subsubsection*{Random modular networks.}

We fix the node-tier structure $(3,3,3,3,3)$ and generate $\overline{m} \in \{5,10,15,20,25,30\}$ distinct supply chains by randomly selecting one node of each tier to be part of the supply chain. We copy the resulting network and obtain two disconnected (but identical) hypergraphs. We then generate $\overline{m}$ ``connecting" supply chains, that is, supply chains that contain one node of each tier from either of the two originally disconnected hypergraphs. We keep track of each supply chain network generated in this process (one for each additional connecting supply chain). We repeat this process four times for each initial parameter combination. Nodes that do not belong to any supply chains are removed; if more than 10\% of the generated nodes have been removed, the instance is discarded. This results in 434 supply chain networks. Based on how we generate the networks, they have comparatively high treewidth, so Gurobi is frequently unable to identify a useful lower bound. Hence, we take the solutions obtained by Gurobi after three hours of run time and experiment with different improvement heuristics. As we cannot identify a single instance in which an improvement is found, we assume the solutions found using Gurobi are sufficiently close to optimal. Note that these networks are designed to obtain varying degrees of modularity (see Appendix \ref{appendix:measures}) while keeping other structural measures largely constant.

\subsection{The relationship between Jaccard clustering and  IP (\ref{eq:MILP_Inter2})} \label{subsec:cluster.IP}

We make explicit the connection between Jaccard clustering and the integer program given in (\ref{eq:MILP_Inter2}). Consider a firm-node $i$ in $G'$ and let $k$ be a firm-node sharing a supply chain $j \in N_{SC}$ with $i$. As $i$ and $j$ (resp. $j$ and $k$) are neighbors in $G'$, $i$ and $j$ (resp. $j$ and $k$) necessarily appear in the same bag in the tree decomposition of $G'.$ If $i,j$ and $j,k$ appear in different bags, then $j$ must appear in all intermediate bags, by definition of a tree decomposition. Thus, in a setting where $i$ and $k$ share many supply chains, i.e., the numerator of $NS(i,k)$ is high, the minimal tree decomposition will likely place $i$ and $k$, together with all shared supply chains $j$, in the same bag to minimize treewidth. If $i,j,k$ are all in the same bag in the tree decomposition of $G'$, then constraint (\ref{eq:milp.constr4}) now applies to the three nodes, i.e., we need $\ell_{ij}+\ell_{jk}+\ell_{ki}\leq 2.$ For this constraint to hold, at least one of $\ell_{ji}, \ell_{jk}, \ell_{ki}$ needs to be equal to zero. If this is either $\ell_{ji}$ or $\ell_{jk}$, then constraints (\ref{eq:milpInter2.const1})  and (\ref{eq:milpInter2.constr2}) become harder to meet, unless we set $s_i$ or $s_k$ to 1, that is, we add $i$ or $k$ to the seed set. This phenomenon is further exacerbated in the setting where the denominator of $NS(i,k)$ is low. In this case, $i$ and $k$ belong to few supply chains in total, which makes constraint  (\ref{eq:milpInter2.constr2}) harder to satisfy as the sum over all supply chains that $i$ (resp. $k$) belongs to only contains few terms. This, in turn, also pushes the seed set to be larger.

\subsection{Description of network measures and parameters} \label{appendix:measures}

We consider the following previously defined measures: the number of nodes $n$, the number of supply chains $m$, the maximum number of nodes per supply chain $k$, the treewidth of the auxiliary graph $\omega'$, and the clustering metric $J$. In addition, we also define the following measures:

\subsubsection*{Alternative clustering metrics.} First, \emph{projection clustering}. We consider the (non-bipartite) projection of $G'$ onto firm-nodes and use the traditional definition of clustering for graphs. More specifically, we construct a graph $G''$ from the $n$ firm-nodes, with an edge between two firm-nodes if they have at least one supply chain in common. The standard clustering coefficient is defined, for example, in \citet{latapy2008basic}, and we take the average over all nodes of the projection. Second, \emph{projection clustering (weighted)} follows the same principles, but each edge in $G''$ is weighted by the number of supply chains the nodes of the edge shave in common.
Third, \emph{hourglass clustering}. On non-bipartite graphs, the clustering coefficient of a node can equivalently be defined as the number of triangles containing the node divided by the number of triplets containing the node with at least two edges. To extend this idea to bipartite graphs, one can divide the number of fully connected quadruplets by the number of quadruplets with at least three links \citep{latapy2008basic}. We employ this extension and consider both the average ratio across firm-nodes and the total ratio throughout the graph. Finally, \emph{repetition of partners}. This measures how many firm-nodes, on average, a given firm-node shares supply chains with.

\subsubsection*{Modularity.}
A commonly used definition is provided by \citet{newman2006modularity}: given a partition of $n$ nodes of a graph $H$ into $z$ groups $\mathcal{P} = (p_1, \ldots, p_z)$, the modularity of $H$ is
$Q=\frac{1}{2 \eta} \sum_{ij} \left( A_{ij} -\frac{\eta_i \eta_j}{2\eta}\right) \delta(p_i,p_j),$
where $\eta$ is the sum of all edge weights in the graph, $\eta_i$ is the sum of the weights of the edges attached to node $i$, $A$ is the (weighted) adjacency matrix of $H$, and $\delta(p_i,p_j)$ is equal to 1 if $p_i=p_j$ (that is, $i$ and $j$ are in the same community) and 0 otherwise.
This definition has no direct extension for hypergraphs and does not apply to bipartite graphs (such as our auxiliary graph). Hence, we use the weighted projection of $G'$ onto firm-nodes, denoted by $G''$, as in the case of \emph{projection clustering (weighted)} to compute modularity. The graph $G''$ is non-bipartite but keeps most of the relevant information about linkages between nodes. We then use the commonly applied Clauset-Newman-Moore greedy modularity maximization algorithm to find the partition leading to the largest $Q$ \citep{clauset2004finding}.

\subsubsection*{Other relevant measures from the supply chain and network literature.}

For our predictive models, we identify \emph{accessibility} and \emph{interconnectedness} as important measures for supply chain networks \citep{bellamy2014influence}. Accessibility refers to the information centrality of nodes, that is, the length of paths ending at a given node. Interconnectedness, meanwhile, refers to the number of shared relationships between connections of a node. As the measures directly relate to diffusion and are defined for undirected graphs, we seek to apply them to the auxiliary graph $G'$. However, they are not well-defined for bipartite graphs, so we consider them on the (weighted) projection $G''$.

Finally, we compute all measures from \citet[][Table 4]{perera2017network} on a modified version of $G$, where directed edges connect nodes in subsequent tiers sharing a supply chain. The authors summarize key supply chain network measures from the empirical literature. 
We omit clustering and modularity, which are already specified. As we deal with a directed graph, we compute \emph{assortativity} based on all in/out degree combinations. This leads to eleven measures, to which we add the \emph{average number of supply chains of a firm}.

\subsection{Predictive models of the seed set size} \label{appendix:prediction}

We first use a random forest regression to predict the inverse logit function of the percentage of seed nodes in the optimal solution on the measures introduced in Appendix \ref{appendix:measures}. In particular, we select 80\% of the instances of a dataset for training and choose hyperparameters by applying 5-fold cross-validation on 100 randomly chosen combinations. We then evaluate the best model identified by computing the root mean square errors (RMSE) for the remaining 20\% of instances between the actual and predicted percentage of seed nodes. RMSEs are $0.021$, $0.024$, and $0.016$ for the three datasets. Unlike the case of linear regression, the importance of each regressor in a random forest can only be computed indirectly by calculating and weighing the Shapley values of the different decision trees for a subset of data. We use the Python package \emph{shap} for this task \citep{NIPS2017_7062}, finding that $J$ has the highest importance among clustering metrics for all networks and either the highest or second-highest importance among all regressors.

\begin{figure}[t]
    \centering
    \subfloat[Relative weight of $J$'s coefficient]{\label{fig:penalty_exp_2}
        \begin{adjustbox}{clip,trim=0cm 0.3cm 0cm 0cm,max width=\textwidth}
        \includegraphics{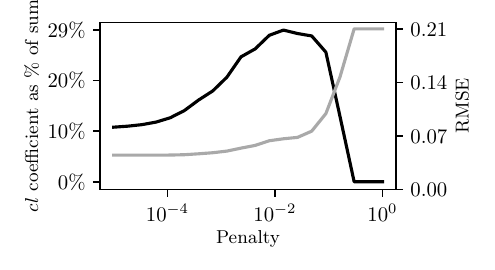}
        \end{adjustbox}
    }
    \hspace{-0.9cm}
    \subfloat[Relative rank of $J$'s coefficient]{\label{fig:penalty_exp_1}
        \begin{adjustbox}{clip,trim=0cm 0.3cm 0cm 0cm,max width=\textwidth}
        \includegraphics{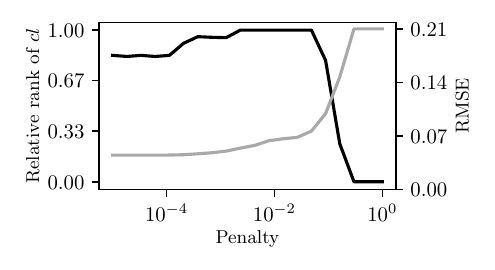}
        \end{adjustbox}
    }
    \caption{Elastic net regression on the \citet{Willems2008} networks with an L1-ratio of 0.8 and varying penalty coefficient. The black (resp.\ gray) lines indicate the value of the left-hand (resp.\ right-hand) axes.}
    \label{fig:penalty_exp}
\end{figure}

To gain more insight into the importance of Jaccard clustering compared to other variables, we introduce an elastic net regression on the same variables, i.e., a regression combining L1-norm (``Lasso'') and L2-norm (``Ridge'') penalty terms, and vary the regression's regularization penalty coefficient on a logarithmic scale. 
Figure~\ref{fig:penalty_exp} depicts the results for the \citet{Willems2008} networks. The results for other datasets are omitted for brevity but show the same patterns. In particular, any model with good explanatory power (before the RMSE increases steeply) puts a high weight on $J$. As we increase the penalty, i.e., requiring the model to have fewer explanatory variables, $J$ becomes more important. In Figure \ref{fig:penalty_exp_2}, we consider the coefficient obtained for each variable during the regression and plot in black the coefficient corresponding to $J$ divided by the sum of all (absolute) coefficients. As can be seen, this ratio increases to a third. In Figure \ref{fig:penalty_exp_1}, we observe that the relative rank (in black) of the absolute value of the coefficient of $J$ increases to 1. The two curves then drop suddenly, but only when we have reached a penalty value that generates a highly biased model (as seen from the ``explosion" of the RMSE in gray). Our results are consistent for all relative weights of L1 and L2, as long as the L1 term is high enough to ensure convergence, i.e., in the range $\frac{L1}{L1+L2} \in [0.5,1.0]$.

\SingleSpacedXI
\renewcommand{\refname}{References for the Appendix}
\putbib[library]
\end{bibunit}
\end{APPENDICES}

\end{document}